\documentclass[11pt]{amsart}

\usepackage{amssymb, amscd}
\usepackage{graphicx}

\newcommand\al{\alpha}

\newcommand\g{\gamma}
\newcommand\Dt{\Delta}

\renewcommand\th{\vartheta}

\newcommand\ld{\lambda}

\newcommand\x{\xi}

\newcommand\Ph{\Phi}
\newcommand\ph{\varphi}

\newcommand\lieg{{\mathfrak g}}
\newcommand\liek{{\mathfrak k}}
\newcommand\lieh{{\mathfrak h}}
\newcommand\liesl{\mathfrak {sl}}
\newcommand\liegl{\mathfrak {gl}}
\newcommand\lieu{{\mathfrak u}}

\newcommand\CC{\mathbb C}

\newcommand\RR{\mathbb R}
\newcommand\ZZ{\mathbb Z}

\newcommand\SL{{\mathrm{SL}}}

\newcommand\GL{{\mathrm{GL}}}
\newcommand\SO{{\mathrm{SO}}}
\newcommand\SU{{\mathrm{SU}}}
\newcommand\U{{\mathrm{U}}}

\newcommand\End{\operatorname{End}}
\newcommand\re{\operatorname{Re}}
\newcommand\im{\operatorname{Im}}
\newcommand\tr{\operatorname{tr}}

\makeatletter 
\newcommand\matc[4]{\left( {#1\@@atop #3}{#2\@@atop #4}\right)}

\newcommand\matr[4]{\left( {\hfill #1\@@atop\hfill #3}{\hfill
#2\@@atop\hfill #4}\right)}

\newcommand\matl[4]{\left( { #1\@@atop #3}{ #2\@@atop\hfill #4}\right)}

\makeatother

\newcommand\widearray[1]{\renewcommand\arraystretch{1.4} \begin{array}{#1}}

\newcommand\lw[1]{\hbox{}_{#1}\!}

\newcommand\vz[1]{\mathchoice{\left\{ #1 \right\}}{\left\{ #1 \right\}}{\{
#1 \}}{\{ #1 \}}}

\newcommand\vzm[2]{\mathchoice{\left\{\, #1 : #2 \,\right\}}{\{\, #1
:\allowbreak #2 \,\}}{\{ #1 :\allowbreak #2 \}} {\{ #1 :\allowbreak #2 \}}}

 \theoremstyle{plain}
\newtheorem{thm}{Theorem}[section]
\newtheorem{lem}[thm]{Lemma}
\newtheorem*{lemsec1}{Lemma \ref{XbetaH}}

\newtheorem*{propsec4}{Proposition \ref{expD1}}
\newtheorem*{prop4.5}{Proposition \ref{expD2}}
\newtheorem*{prop4.6}{Proposition \ref{expE1}}
\newtheorem*{prop4.7}{Proposition \ref{expE2}}
\newtheorem*{lem5.3}{Lemma \ref{dderx1}}
\newtheorem*{lem5.4}{Lemma \ref{ddery1}}
\newtheorem*{lem5.5}{Lemma \ref{ddery2}}
\newtheorem*{lem5.6}{Lemma \ref{dderx2}}
\newtheorem*{lem5.7}{Lemma \ref{ddery1x1}}
\newtheorem*{lem5.8}{Lemma \ref{derx2y1}}
\newtheorem*{lem5.9}{Lemma \ref{derx2y2}}

\newtheorem*{esfneg}{Theorem \ref{esfneg}}
\newtheorem{prop}[thm]{Proposition}
\newtheorem{cor}[thm]{Corollary}
\newtheorem{conj}[thm]{Conjecture}

\theoremstyle{definition}
\newtheorem{defn}[thm]{Definition}

\theoremstyle{remark}

\newcommand\rmk[1]{\medskip\par\noindent{\em #1. }\ignorespaces}

\title[Matrix valued spherical functions]{\sc Matrix valued spherical
functions
associated to the complex projective plane}

\author{ F. A. Gr\"unbaum}
\address{Departament of Mathematics, University of California, Ber\-ke\-ley
CA 94705}
\email{grunbaum@math.berkeley.edu}

\author{I. Pacharoni}
\address{CIEM-FaMAF, Universidad Nacional de C\'or\-do\-ba,
C\'or\-do\-ba~5000, Argentina}
\email{pacharon@mate.uncor.edu}

\author{J. Tirao}
\address{CIEM-FaMAF, Universidad Nacional de C\'or\-do\-ba,
C\'or\-do\-ba~5000, Argentina}
\email{tirao@mate.uncor.edu}

\dedicatory{TO E. H. ZARANTONELLO, teacher and friend.}
\thanks{This paper is partially supported by NSF grants FD9971151 and
1-443964-21160 and by CONICET grant PIP655-98.}

\begin{document}
\begin{abstract} The main purpose of this paper is to compute all
irreducible spherical functions
on $G=\SU(3)$ of arbitrary type $\delta\in \hat K$, where
$K={\mathrm{S}}(\U(2)\times\U(1))\simeq\U(2)$. This is accomplished by
associating to a spherical
function $\Phi$ on $G$ a matrix valued function $H$ on the complex
projective plane $P_2(\CC)=G/K$.
It is well known that there is a fruitful connection between the
hypergeometric function of Euler and
Gauss and the spherical functions of trivial type associated to a rank one
symmetric pair $(G,K)$.
But the relation of spherical functions of types of dimension bigger than
one with classical analysis,
has not been worked out even in the case of an example of a rank one pair.
The entries of $H$ are solutions of two systems of ordinary differential
equations.
There is no ready made approach to such a pair of systems, or even to a
single system of this kind.
In our case the situation is very favorable and the solution to this pair
of systems can be exhibited
explicitely in terms of a special class of  generalized hypergeometric
functions $\lw{p+1}F_p$.
\end{abstract}
\maketitle

\section{Introduction and statement of results}\label{statements}
The complex projective plane can be realized as the homogeneous space
$G/K$, where $G=\SU(3)$ and
$K={\mathrm{S}}(\U(2)\times\U(1))$. We are interested in determining, up to
equivalence, all irreducible spherical functions,
associated to the pair $(G,K)$. If $(V,\pi)$ is a finite dimensional
irreducible representation of $K$ in the equivalence
class $\delta \in \hat K$, a spherical function on $G$ of type $\delta$ is
characterized by
\begin{enumerate}
\item[i)] $\Ph:G\longrightarrow \End(V)$ is analytic.
\item[ii)] $\Ph(k_1gk_2)=\pi(k_1)\Ph(g)\pi(k_2)$, for all $k_1,k_2\in K$,
$g\in G$, and $\Phi(e)=I$.
\item[iii)] $[\Dt_2\Ph ](g)=\ld_2\Ph(g)$, $[\Dt_3\Ph ](g)=\ld_3\Ph(g)$ for
all $g\in G$ and for
some $\ld_2,\ld_3\in \CC$.
\end{enumerate}
Here $\Dt_2$ and $\Dt_3$ are two algebraically independent generators of the
polynomial algebra $D(G)^G$ of all differential operators on $G$ which are
invariant under left and right multiplication by
elements in $G$. A particular choice of these operators is given in
Proposition \ref{defDelta}.

The set $\hat K$ can be identified with the set $\ZZ\times \ZZ_{\geq 0}$.
If $k=\left(\begin{smallmatrix} A&0\\ 0& a
\end{smallmatrix}\right)$, with $A\in \U(2)$ and $a=(\det A)^{-1}$, then
$$\pi(k)=\pi_{n,\ell}(A)=\left(\det A\right )^n\,A^\ell,$$
where $A^\ell$ denotes the $\ell$-symmetric power of $A$, defines an
irreducible representation of $K$ in the class
$(n,\ell)\in \ZZ\times \ZZ_{\geq 0}$.

For $n\geq 0$, the
representation $\pi_{n,\ell}$ of $\U(2)$ extends to a unique holomorphic
multiplicative map of
${\mathrm M}(2,\CC)$ into $\End(V_\pi)$, which we shall still denote by
$\pi_{n,\ell}$. For any
$g\in {\mathrm M}(3,\CC)$, we shall denote by $A(g)$ the left upper
$2\times 2$ block of $g$, i.e.
$$A(g)=\matc{g_{11}}{g_{12}}{g_{21}}{g_{22}}.$$

For any $\pi=\pi_{n,\ell}$ with $n\geq 0$ let $\Phi_\pi:G\longrightarrow
\End(V_\pi)$ be defined by
$$\Phi_\pi(g)=\Phi_{n,\ell}(g)= \pi_{n,\ell}(A(g)).$$
Then, in Theorem \ref{Phipiesf} we prove that $\Phi_\pi$ is a spherical
function of type $(n,\ell)$.
This is the main result of Section \ref{Spherical functions}. These
particular spherical functions will play a crucial role
in the rest of the paper.

Observe that when $n<0$ the function $\Phi_\pi$ above can be defined on the
open set
$${\mathcal A}=\vzm{g\in G}{\det A(g)\neq 0}.$$
The group $G=\SU(3)$ acts in a natural way in the complex projective plane
$P_2(\CC)$. This action is transitive and $K$ is
the isotropy subgroup of the point $(0,0,1)\in P_2(\CC)$. Therefore $
P_2(\CC)= G/ K$. We shall identify the complex plane
$\CC^2$ with the affine plane $\vzm{(x,y,1)\in P_2(\CC)}{(x,y)\in \CC^2}$.

The canonical projection $p:G\longrightarrow P_2(\CC)$ maps the open dense
subset ${\mathcal A}$ onto the affine plane
$\CC^2$. Observe that ${\mathcal A}$ is stable by left and right
multiplication by elements in $K$.

To determine all spherical functions
$\Phi:G\longrightarrow \End(V_\pi)$ of type $\pi=\pi_{n,\ell}$ we use the
function $\Phi_\pi$ in the following way: in the open
set ${\mathcal A}$ we define a function $H$ by
$$H(g)=\Phi(g)\, \Phi_\pi(g)^{-1},$$
where $\Phi$ is supposed to be a spherical function of type $\pi$.
Then $H$ satisfies
\begin{enumerate}
\item [i)] $H(e)=I$.
\item [ii)]$ H(gk)=H(g)$, for all $g\in {\mathcal A}, k\in K$. \item [iii)]
$H(kg)=\pi(k)H(g)\pi(k^{-1})$, for all
$g\in {\mathcal A}, k\in K$.
\end{enumerate}

\noindent Property ii) says that $H$ may be considered as a function on
$\CC^2$.

The fact that $\Phi$ is an eigenfunction of $\Delta_2$ and $\Delta_3$,
makes $H$ into an eigenfunction of certain
differential operators $D$ and $E$ on $\CC^2$. For completeness the
explicit computation of these operators is carried out
fully in Section \ref{reduction}.

In Section \ref{secc1var} we take full advantage of the $K$-orbit structure
of $P_2(\CC)$ combined with property iii) of our
functions $H$. The affine plane $\CC^2$ is $K$-stable and the corresponding
line at infinity $L$ is a $K$-orbit. Moreover the
$K$-orbits in $\CC^2$ are the spheres
$S_r=\vzm{(x,y)\in \CC^2}{|x|^2+|y|^2=r^2}.$ Thus we can take the points
$(r,0)\in S_r$ and $(1,0,0)\in L$ as
representatives of $S_r$ and $L$, respectively. Since
$(M,0,1)=(1,0,\textstyle\frac 1M)\longrightarrow (1,0,0)$ when
$M\rightarrow \infty$, the closed interval $[0,\infty]$ parametrizes the
set of $K$-orbits in $P_2(\CC)$.

Thus there exist ordinary differential operators $\tilde D$ and $\tilde E$
on the open interval $(0,\infty)$ such that
$$(D\,H)(r,0)=(\tilde D\tilde H)(r)\, , \quad (E H)(r,0)=(\tilde E\tilde
H)(r),$$
where $\tilde H(r)=H(r,0)$, $r\in(0,\infty)$.
These operators $\tilde D$ and $\tilde E$ are explicitely given in Theorems
\ref{mainth1} and \ref{mainth2}. We need to
compute a number of second order partial derivatives of the function
$H:\CC^2\longrightarrow \End(V_\pi)$ at the point
$(r,0)$. This detailed computation is broken down in a number of lemmas
included in the Appendix for the benefit of the reader.

Theorems \ref{mainth1} and \ref{mainth2} are given in terms of linear
transformations. The functions $\tilde H$ turn out to
be diagonalizable (Proposition \ref{Hdiagonal}). Thus, in an appropriate
basis of
$V_\pi$ we can write $\tilde H(r)=(h_0(r),\cdots,
h_\ell(r))$. Then we give in Corollaries \ref{sistema} and \ref{sistema2}
the corresponding statements of these theorems in
terms of the scalar functions $h_i$.

In Section \ref{extension} we take into account the behavior of the
function $\tilde H$ associated to a spherical function
$\Phi$ when $r\rightarrow \infty$. The corresponding asymptotic behavior of
$\tilde H$ is given in the first part of
Proposition \ref{extG}. Our strategy to obtain all spherical functions of
$(G,K)$ will be to find all eigenfunctions $\tilde
H$ of $\tilde D$ and $\tilde E$ which satisfy this behavior at infinity.
This is the point of the second part of Proposition
\ref{extG}.

In Section \ref{innerprod} we consider the natural inner product among
continuous maps
from $G$ to $\End(V_\pi)$, which makes $\Delta_2$ and $\Delta_3$ symmetric
(Proposition
\ref{Deltasim}). In Proposition \ref{prodint} we give the explicit
expresion of
this inner product when restricted
to those functions which satisfy $\Phi(k_1gk_2)=\pi(k_1)\Phi(g) \pi(k_2)$.

We also introduce here the variable $t=(1+r^2)^{-1}$ which converts the
operators $\tilde D$ and $\tilde E$
into new operators $D$ and $E$. Then we present a classical analysis
motivation for the inner product
given in Proposition \ref{prodint}.

In Section \ref{autofunciones} we give a complete description of a method
to obtain all $C^\infty$ solutions of $DH=\ld H$.
This construction rests on some remarkable factorizations of higher order
differential operators with Gauss' hypergeometric
operator as one of the factors. This construction revolves around
introducing an integer parameter $0\leq k \leq \ell$ as well
as a parameter $w$ involved in expressing $\ld$ appropriately.

The case $\ell=0$ is analized in detail and the solutions are given
generically in terms of $\lw{2}F_1$.
For higher $\ell$ the solutions are obtained by acting on $\lw{2}F_1$
either by differential operators or by appropriate
shift operators.

As we remark in Section \ref{Spherical functions}, there exists a one to
one correspondence between the set of all
equivalence classes of finite dimensional irreducible representations of
$G$ which contain the representation
$\pi_{n,\ell}$ of $K$, and the set of all equivalence classes of
irreducible spherical functions of $(G,K)$ of type
$(n,\ell)$. This is the starting point we take in Section \ref{parametriz}
to get a parametrization of all these equivalence
classes in terms of all tuples $(p,q,k_1,k_2)\in \ZZ^4$ such that $p+q\geq
k_1\geq q \geq k_2 \geq 0$, with $n=k_1+2k_2-p-2q$,
$\ell=k_1-k_2$. This is the content of Corollary \ref{param2}. Also in
Corollary \ref{expldmu} we give the eigenvalues $\ld$
and $\mu$ which correspond to the differential operators $D$ and $E$ in
terms of these parameters. Moreover, we introduce here
two integer parameters $w$, $k$ subject to the following four inequalities:
$0\leq w$, $0\leq k\leq \ell$, $0\leq w+n+k$,
which give a very convenient parametrization of the irreducible spherical
functions of type $(n,\ell)$.

Section \ref{aeigenfunc} starts with Proposition \ref{Hanali} where the
local behaviour of the functions $H=H(t)$ in the
variable $t=(1+r^2)^{-1}$ at $t=0$ associated to a spherical function
$\Phi$ of type $(n,\ell)$ is established. In particular
it is proved that $H$ is analytic and that has a zero of order at least
$-n-\ell$, when $n+\ell<0$. Thus we aim at getting
joint solutions to $DH=\ld H$ and $EH=\mu H$ in the form of formal power
series $H(t)=\sum_{j\geq 0} t^j H_j$. Since $E$
commutes with $D$ we look at the linear map defined by it on the vector
space $V(\ld)$ of all formal solutions of $DH=\ld H$.
This equation gives rise to a three term recursion relation for the
coefficients $H_j$.

When $n\geq 0$, one sees that the linear map $\eta:H(t)\mapsto H(0)$ is an
isomorphism
of $V(\ld)$ onto $\CC^{\ell+1}$. Using $\eta$, $E$ becomes an
$(\ell+1)\times (\ell+1)$ matrix $L=L(\ld)$. In Corollary
\ref{eigenv}, $\eta$ is used to establish a one to one correspondence
between joint solutions to $DH=\ld H$, $EH=\mu H$ and
the eigenvalues $\mu$ of $L(\ld )$. Moreover there is no joint
eigenfunction given by formal power series unless the pair
$(\ld,\mu)$ lies on the union of $(\ell+1)$ straight lines. For any point
on this curve, if $H_0$ is an eigenvector of
$L(\ld)$ with eigenvalue $\mu$, the recursion relation becomes a two term
recursion and can be solved explicitely in terms of
generalized hypergeometric functions. In particular they converge for
$|t|<1$.

When $n<0$ the isomorphism $\eta$ must be redefined and the corresponding
matrix $L$ is not as simple as the one obtained
when $n\ge0$. This suggest that our choice of the function $H$ associated
to a spherical function is most appropriate only
when $n\ge0$.

In Section \ref{solutions} we illustrate in full detail the results
obtained above in
the cases $\ell=0,1,2$ and $n$ an arbitrary integer. The main ingredient
here is a particular class of generalized
hypergeometric functions among those of the form $\lw{p+1}F_p$. For each
$\ell$, we need to deal with several values of $p$ in
the range $1\le p\le \ell+1$ and thus for $\ell=0$ (and $n=0$) we obtain
the known result involving only $\lw{2}F_1$.

Generalized hypergeometric functions have numerator and denominator
parameters, as well as an independent variable $t$. The
functions that enter in our explicit formulas satisfy the condition that
all but two of the numerator parameters exceed by one
all but one of the denominator parameters. In the classical case, with
$\ell=0$, there is not enough room for this phenomenon
to show up.

We state a fairly explicit conjecture on the dependence of these
hypergeometric functions on its main parameters, namely the
ones that are not related by the shift described above. These two numerator
parameters and the denominator one exhaust the
ones that appear in the case of $\lw{2}F_1$.

In view of the very important role played by all functions
$\lw{p+1}F_p$, $p$ arbitrary, in our description of the entries of the
matrix valued spherical functions, it is worth to note that certain
generalized hypergeometric functions, more precisely the functions
$\lw{3n-4}F_{3n-5}$, have appeared in the expression of matrix entries of
representations of the group U$(n)$. In this case one can consult \cite{KG} as
well as the very systematic treatment in \cite{KV}. This very nice result extends
to the case of arbitrary $n$ the well known results involving Gauss'
function $\lw{2}F_1$ in the case of U$(2)$. It would
be nice to see the functions identified here playing a useful role in other
problems in geometry or physics.

In Section \ref{recur} we observe that our spherical functions satisfy not
only a
differential equation, but also an appropriate recursion relation in the
discrete variable $w$. We display the results in
full in the case $\ell=2$ and $n$ nonnegative.

This three term recursion relation is well known in the classical case,
when $\ell=0$, and it gives the main way to compute in
a numerically stable way the so called classical orthogonal polynomials.
Its importance is not only a practical matter.   It
gives one of the earliest instances of the bispectral property. By repeated
applications of the Darboux process this innocent
looking property is tied up with all sorts of other issues of recent
interest in mathematical physics, including nonlinear
integral evolution equations, $W$-algebras, interesting isomonodromy
deformations, the whole area of random matrix models,
etc. There are even $q$ versions of all of this. Moreover this is not
restricted to the case of one (either spatial or
spectral) variable.

It is natural to wonder about the relation between the matrix valued
functions
constructed in Section 12 and the relatively new theory of matrix valued
orthogonal polynomials. This issue has been addressed in \cite{GPT}. It
suffices to state here that the functions considered in Section 12 do not
satisfy all the conditions in the theory, see for instance \cite{DVA}.

Finally we remark that we have chosen this example as one of the simplest
to analyze among the symmetric spaces of rank one
leading to matrix valued spherical functions. We hope to deal with other
simple examples in the near future. One important
case of scalar valued spherical functions of non trivial type has been
considered in Part I, Chapter 5 of \cite{HS}.
The referee has wondered what is the perspective for doing analogous things
in the case of SU$(n)$, and E. Stein had also suggested to consider this
more
general problem. At this point this appears to be a very interesting
challenge.

One can only speculate that many of the connections that make Gauss'
function a vital part of mathematics at the end of the
twenty century will be shared by its matrix valued version discussed here.

It is a pleasure to thank the referee for a very thorough job. In
particular the referee has redone independently a few of the
computer algebra checks that we had originally done with Maxima
by using very nice Maple worksheets.

\section{Spherical functions}\label{Spherical functions}
Let $G$ be a locally compact unimodular group and let $K$ be a compact
subgroup of $G$. Let $\hat K$ denote the set of all
equivalence classes of complex finite dimensional irreducible
representations of $K$; for each $\delta\in \hat K$, let
$\x_\delta$ denote the character of $\delta$, $d(\delta)$ the degree of
$\delta$, i.e. the dimension of any representation in
the class $\delta$, and $\chi_\delta=d(\delta)\x_\delta$. We shall choose
once and for all the Haar measure $dk$ on
$K$ normalized by
$\int_K dk=1$.

We shall denote by $V$ a finite dimensional vector space over the field
$\CC$ of complex numbers and by $\End(V)$ the space
of all linear transformations of $V$ into $V$. Whenever we shall refer to a
topology on such vector space we shall be talking
about the unique Hausdorff linear topology on it.

By definition a zonal spherical function (\cite{He}) $\ph$ on $G$ is a
continuous complex valued function which satisfies
$\ph(e)=1$ and
\begin{equation}\label{defclasica}
\ph(x)\ph(y)=\int_K \ph(xky)\, dk \qquad \qquad x,y\in G.
\end{equation}

A fruitful generalization of the above concept is given in the following
definition.
\begin{defn}\label{defesf}(\cite{T},\cite{GV})
A spherical function $\Ph$ on $G$ of type $\delta\in \hat K$ is a
continuous function on $G$ with values in $\End(V)$ such
that
\begin{enumerate} \item[i)] $\Ph(e)=I$. ($I$= identity transformation).

\item[ii)] $\Ph(x)\Ph(y)=\int_K \chi_{\delta}(k^{-1})\Ph(xky)\, dk$, for
all $x,y\in G$.
\end{enumerate}
\end{defn}

\begin{prop}\label{propesf}(\cite{T},\cite{GV})  If $\Ph:G\longrightarrow
\End(V)$ is a spherical function of type $\delta$
then:
\begin{enumerate}
\item[i)] $\Ph(kgk')=\Ph(k)\Ph(g)\Ph(k')$, for all $k,k'\in K$, $g\in G$.
\item[ii)] $k\mapsto \Ph(k)$ is a representation of $K$ such that any
irreducible subrepresentation belongs to $\delta$.
\end{enumerate}
\end{prop}

Concerning the definition let us point out that the spherical function
$\Ph$ determines its type univocally
(Proposition \ref{propesf}) and let us say that the number of times that
$\delta$ occurs in the representation
$k\mapsto \Ph(k)$ is called the {\em height} of $\Ph$.

Let $\ph$ be a complex valued continuous solution of the equation
\eqref{defclasica}. If $\ph$ is not identically zero then
$\ph(e)=1$. (cf. \cite{He}, p. 399). This result generalizes in the
following way: we shall say that a function
$\Ph:G\longrightarrow \End(V)$ is {\em irreducible} whenever $\Ph(g)$,
$g\in G$, is an irreducible family of linear
transformations of $V$ into $V$. Then we have

\begin{prop}(\cite{T})
Let $\Ph$ be an $\End(V)$-valued continuous solution of equation ii) in
Definition \ref{defesf}. If $\Ph$ is irreducible then
$\Ph(e)=I$.
\end{prop}

Spherical functions of type $\delta$ arise in a natural way upon
considering representations of $G$.
If $g\mapsto U(g)$ is a continuous representation of $G$ , say on a
complete, locally convex, Hausdorff topological vector
space $E$, then $$P(\delta)=\int_K \chi_\delta(k^{-1})U(k)\, dk$$ is a
continuous projection of $E$ onto
$P(\delta)E=E(\delta)$; $E(\delta)$ consists of those vectors in $E$, the
linear span of whose $K$-orbit is finite dimensional
and splits into irreducible $K$-subrepresentations of type $\delta$.
Whenever $E(\delta)$ is finite dimensional, the function
$\Ph:G\longrightarrow \End(E(\delta))$ defined by $\Ph(g)a=P(\delta)U(g)a$,
$g\in G, a\in E(\delta)$ is a spherical function
of type $\delta$. In fact, if $a\in E(\delta)$ we have
\begin{align*}
\Ph(x)\Ph(y)a&= P(\delta)U(x)P(\delta)U(y)a=\int_K \chi_\delta(k^{-1})
P(\delta)U(x)U(k)U(y)a\, dk\\
&=\left(\int_K\chi_\delta(k^{-1})\Ph(xky)\, dk\right) a. \end{align*}
If the representation $g\mapsto U(g)$ is topologically irreducible (i.e.
$E$ admits no non-trivial
closed $G$-invariant subspace) then the associated spherical function $\Ph$
is also irreducible.

If a spherical function $\Phi$ is associated to a Banach representation of
$G$ then it is quasi-bounded,
in the sense that there exists a semi-norm $\rho$ on $G$ and $M\in\RR$ such
that $\|\Phi(g)\|\le M\rho(g)$
for all $g\in G$. Conversely, if $\Phi$ is an irreducible quasi-bounded
spherical function on $G$, then
it is associated to a topologically irreducible Banach representation of
$G$ (Godement, see \cite{T}). Thus if
$G$ is compact any irreducible spherical function on $G$ is associated to a
Banach representation of $G$,
which is finite dimensional by Peter-Weyl theorem.

{}From now on we shall assume that $G$ is a connected Lie group. Then it is
not difficult to prove that any spherical function $\Ph:G\longrightarrow
\End(V)$ is
differentiable ($C^\infty$), and moreover that it is analytic. Let $D(G)$
denote the algebra of all
left invariant differential operators on $G$ and let $D(G)^K$ denote the
subalgebra of all operators in
$D(G)$ which are invariant under all right
translation by elements in $K$.

\smallskip
In the following proposition $(V,\pi)$ will be a finite dimensional
representation of $K$ such that any irreducible
subrepresentation belongs to the same class $\delta\in\hat K$. \begin{prop}
\label{defeq}(\cite{T},\cite{GV}) A function
$\Ph:G\longrightarrow \End(V)$ is a spherical function of type $\delta$ if
and only if
\begin{enumerate}
\item[i)] $\Ph$ is analytic.
\item[ii)] $\Ph(k_1gk_2)=\pi(k_1)\Ph(g)\pi(k_2)$, for all $k_1,k_2\in K$,
$g\in G$, and
$\Phi(e)=I$.
\item[iii)] $[D\Ph ](g)=\Ph(g)[D\Ph](e)$, for all $D\in D(G)^K$, $g\in G$.
\end{enumerate}
\end{prop}

Let us observe that if $\Ph:G\longrightarrow \End(V)$ is a spherical
function then $\Ph:D\mapsto
[D\Ph](e)$ maps $D(G)^K$ into $\End_K(V)$ ($\End_K(V)$ denotes the space of
all linear maps of $V$
into $V$ which conmutes with $\pi(k)$ for all $k\in K$) defining a finite
dimensional
representation of the associative algebra $D(G)^K$. Moreover the spherical
function is irreducible
if and only if the representation $\Ph: D(G)^K\longrightarrow \End_K (V)$
is surjective. As a
consequence of this we have:

\begin{prop}(\cite{T},\cite{GV}) The following properties are equivalent:
\begin{enumerate}
\item[i)] $D(G)^K$ is commutative.
\item[ii)] Every irreducible spherical function of $(G,K)$ is of height
one. \end{enumerate}
\end{prop}

In this paper the pair $(G,K)$ is $(\SU(3),{\mathrm S}(\U(2)\times
\U(1)))$. Then it is known
that $D(G)^K$ is abelian in this case; moreover $D(G)^K$ is isomorphic to
$D(G)^G\otimes
D(K)^K$ (cf. \cite{Co}, \cite{Kp}), where $D(G)^G$ (resp. $D(K)^K$) denotes
the subalgebra of all
operators in $D(G)$ (resp. $D(K)$) which are invariant under all right
tranlations from $G$ (resp.
$K$). Now a famous theorem of Harish-Chandra says that $D(G)^G$ (resp.
$D(K)^K$) is a polynomial
algebra in two algebraically independent generators $\Delta_2$ and
$\Delta_3$ (resp. $Z$ and
$\Delta_K$).

The first consequence of this is that all irreducible spherical
functions of our pair $(G,K)$
are of height one.

The second consequence is that to find all spherical functions of type
$\delta\in \hat K$ is
equivalent to taking any irreducible representation $(V,\pi)$ of $K$ in the
class $\delta$ and to
determine all analytic functions
$\Ph:G\longrightarrow \End(V)$ such that \begin{enumerate}
\item [(1)] $\Ph(k_1gk_2)=\pi(k_1)\Ph(g)\pi(k_2)$, for all $k_1,k_2\in K$,
$g\in G$.

\smallskip
\item [(2)] $[\Delta_j\Ph](g)=\Ph(g)[\Delta_j\Ph](e)$, $j=2,3$.
\end{enumerate}

In fact, because $Z$ and $\Delta_K$ are in $D(K)^K$ and $\Ph$ satisfies (1)
we have
$$[Z\Ph](g)=\Ph(g)\,\dot\pi(Z)=\Ph(g)[Z\Ph](e)$$ and
$$[\Delta_K\Ph](g)=\Ph(g)\dot\pi(\Delta_K)=\Ph(g)[\Delta_K\Ph](e).$$
Here $\dot\pi:\liek_\CC\longrightarrow \End(V_\pi)$ denotes the derivative
of the representation $\pi$ of $K$.  We also
denote with $\dot\pi$ the representation of $D(K)$ in $\End(V_\pi)$ induced
by $\dot\pi$.

Therefore an analytic function $\Ph$ which satisfies (1) and (2) verifies
conditions i), ii) and
iii) of Proposition \ref{defeq}, and hence it is a spherical function.

\medskip
The group $G=\SU(3)$ consists of all $3\times 3$ unitary matrices of
determinant one. The subgroup
$K={\mathrm S}(\U(2)\times \U(1))$ consists of all unitary matrices of the
form
$$k=\left(\begin{matrix} A&0\\ 0& a \end{matrix}\right)$$ where $A\in
\U(2)$ and $a=(\det A)^{-1}$.

Clearly the map $k\mapsto A$ defines a Lie isomorphism of $K$ onto $\U(2)$.
This isomorphism will be used freely in what
follows. In particular $\hat K$ will be identified with $\hat{\U(2)}$. Let
us recall that the identity representation $\pi_1$
of $\U(2)$ in $\CC^2$, as well as the
$\ell$-symmetric power of it $\pi_\ell:A\mapsto A^\ell$ of dimension
$\ell+1$ are irreducible. Moreover the representations
$\pi_{n,\ell}$ of $\U(2)$ defined by
$$\pi_{n,\ell}(A)=\left(\det A\right )^n\, A^\ell\qquad n\in \ZZ, \ell\in
\ZZ_{\geq 0}$$ is a complete set of
representatives of elements in $\hat{\U(2)}$. Thus $\hat{\U(2)}$ can be
identified with the set $\ZZ\times \ZZ_{\geq 0}$.

Now let us observe that the representation $\pi_{n,\ell}$ extends to a
unique holomorphic representation of $\GL(2,\CC)$  and
that when $n\geq 0$ this extends in turn to a unique holomorphic
multiplicative map of ${\mathrm M}(2,\CC)$ into
$\End(V_\ell)$ which we shall still denote by $\pi_{n,\ell}$.

For any $g\in {\mathrm M}(3,\CC)$, we shall denote by $A(g)$ the left upper
$2\times 2$ block of $g$, ie.
$$A(g)=\matc{g_{11}}{g_{12}}{g_{21}}{g_{22}}.$$

Now we are in a position to introduce a particular spherical function on
$G$ of type $(n,\ell)$ for all $n\geq 0$.  For any
$\pi=\pi_{n,\ell}$ with $n\geq 0$ we define the function
$\Phi_\pi:G\longrightarrow \End(V_\pi)$ in the following way:
$$\Phi_{n,\ell}(g)=\Phi_\pi(g)= \pi(A(g)). $$

Observe that when $n<0$ the above formula defines a function $\Phi_\pi$ on
the set ${\mathcal A}=\vzm{g\in G}{\det A(g)\neq
0}$. Although in this case $\Ph_\pi$ is no longer spherical it will have a
useful role to play.

To state the following lemma we need some notation. Let $E_{i,j}$ denote
the $3\times 3$ matrix with $1$ in the
$(i,j)$-place and $0$ elsewhere. We define the elements $X_\beta$ and
$X_\g$ in the complexified Lie algebra $\liesl(3,\CC)$ of $\SU(3)$
by: $X_\beta=E_{2,3}$ and $X_\g=E_{1,3}$.

\begin{lem} \label{holomorfia}The function $\Phi_\pi$ has the following
properties:

\begin{enumerate}
\item[i)] $\Phi_\pi(k)=\pi(k)$, for all $k\in K$.

\item[ii)] $\Phi_\pi(k_1gk_2)=\Phi_\pi(k_1)\Phi_\pi(g)\Phi_\pi(k_2)$, for
all $k_1,k_2\in K,\,g\in \mathcal A$.

\item[iii)] $X_\beta(\Phi_\pi)(g)=X_\g(\Phi_\pi)(g)=0$, for all $g\in
\mathcal A$. \end{enumerate}
\end{lem}
\begin{proof}
i) is obvious and ii) follows directly from the definition of $\Phi_\pi$
upon observing that $A(gk)=A(g)A(k)$ and
$A(kg)=A(k)A(g)$ for all $g\in G$, $k\in K$.

iii) The function $\Phi_\pi$ extends to a unique holomorphic function
$$\Phi_\pi: \SL(3,\CC)\longrightarrow \End(V).$$
Therefore for $g\in G$ we have
$$ X_\beta \left(\Phi_\pi\right)(g)={\frac{d}{d\,t}\Big|}_{t=0}\Phi_\pi(g
\,\text{exp} tX_\beta)=0,$$
because $A(g \,\text{exp}\,tX_\beta)=A(g)$ for all $t$. In a similar way it
follows that $X_\g(\Phi_\pi)=0$.
\end{proof}

Before proving the next theorem we need some more preparation.

\begin{lem}\label{menores} Let $g\in \SU(3)$ and let $g_{(i|j)}$ be the
$2\times 2$  matrix obtained from $g$ deleting the
$i$-th row and the $j$-th column. Then
$$\overline g_{ij}=(-1)^{i+j}\det(g_{(i|j)}).$$ \end{lem}

The set ${\mathcal A}=\vzm{g\in G}{\det A(g)\neq 0}$ is an open dense
subset of $G$. Observe that by Lemma
\ref{menores}, ${\mathcal A}$ can also be described as the set of all $g\in
G$ such that $g_{33}\neq 0$ and that ${\mathcal
A}$ is stable by left and right multiplication by elements in $K$.

The group $G=\SU(3)$ acts in a natural way in the complex projective plane
$P_2(\CC)$. This action is transitive and $K$ is
the isotropy subgroup of the point $(0,0,1)\in P_2(\CC)$. Therefore $$
P_2(\CC)\simeq G/ K.$$ Moreover the $G$-action on
$P_2(\CC)$ corresponds to the action induced by left multiplication on
$G/K$. We shall identify the complex plane $\CC^2$
with the affine plane $\vzm{(x,y,1)\in P_2(\CC)}{(x,y)\in \CC^2}$ under the
map $ (x,y)\mapsto (x,y,1)$.

The projection map $p:G\longrightarrow P_2(\CC)$ defined by $p(g)=g\cdot
(0,0,1)$ maps the open set ${\mathcal A}$ onto the
affine plane $\CC^2$.

Let us consider on $\CC^2$ the real linear coordinates $(x_1,x_2,y_1,y_2)$
defined by: $x_1(x,y)+i\,x_2(x,y)=x$ and
$y_1(x,y)+i\,y_2(x,y)=y$ for all $x,y\in \CC^2$. We also introduce the
following usual notation: \begin{equation}
\frac{\partial }{\partial x}=\frac 12 \left( \frac{\partial }{\partial
x_1}- i\,\frac{\partial }{\partial x_2}\right)\, ,
\quad \frac{\partial }{\partial y}=\frac 12 \left( \frac{\partial
}{\partial y_1}- \,i\frac{\partial }{\partial y_2}\right)
\end{equation}

\

The next lemma will be proved in Section \ref{generalidades}.

\begin{lem}\label{XbetaH} Given $H\in C^\infty (\CC^2)$ denote also by
$H\in C^\infty ({\mathcal A})$ the function defined by
$H(g)=H(p(g))$, $g\in {\mathcal A}$. Then we have $$(X_\beta H)(g)= {
\frac{\overline g_{11}}{g_{33}^2}}\; \frac{\partial
H}{\partial y}-{\frac{\overline g_{21}}{g_{33}^2}}\; \frac{\partial
H}{\partial x}\, , \qquad (X_\g H)(g)=- { \frac{\overline
g_{12}}{g_{33}^2}}\; \frac{\partial H}{\partial y}+{\frac{\overline
g_{22}}{g_{33}^2}}\; \frac{\partial H}{\partial x}.$$
\end{lem}

\begin{cor}\label{antiholom}
A function $H\in C^\infty (\CC^2)$ is antiholomorphic if and only if
$(X_\beta H)(g)=(X_\g H)(g)=0$ for all $g\in {\mathcal
A}.$
\end{cor}
\begin{proof}
The Cauchy Riemann equations say that $H$ is antiholomorphic precisely when
$ \frac{\partial H}{\partial x}=
\frac{\partial H}{\partial y}=0$. Now the corollary follows from Lemma
\ref{XbetaH} upon observing that for $g\in \mathcal A$
the matrix $$\frac{1}{g_{33}^2} \matc{\overline g_{11} }{-\overline g_{21}
}{-\overline g_{12} }{\overline g_{22} }$$ is non
singular. \end{proof}

\begin{thm}\label{Phipiesf} For $n\geq 0$ and $\pi=\pi_{n,\ell}$, the
function $\Phi_\pi$ is an irreducible
spherical function of type $\pi$.
\end{thm}
\begin{proof} We shall prove that $\Ph_\pi$ satisfies the integral equation
$$\Ph_\pi(g)\Ph_\pi(a)=\int_K\chi_{\pi}(k^{-1})\Ph_\pi(gka)\, dk,$$
for all $g,a\in G$.

We fix $g\in G$ and we consider the function $F:{\mathcal A}\longrightarrow
\End(V_\pi)$ defined by
$$F(a)=\left(\int_K \chi_{\pi}(k^{-1})\Ph_\pi(gka)\, dk \right)
\Ph_\pi(a)^{-1}.$$
Then we have $F(ak)=F(a)$, for all $k\in K$. Therefore we may consider $F$
as a function defined on the affine plane $\CC^2$
in $P_2(\CC)$. By Lemma \ref{holomorfia} iii) we have $X_\beta F=X_\g F=0$.
Thus $F:\CC^2\longrightarrow \End(V_\pi)$ is an
antiholomorphic function (Corollary \ref{antiholom}).

We fix $a\in {\mathcal A}$ and let $p(a)=(x,y,1)$. Now, let us consider the
function $f(w)=F(\overline w x,\overline w y)$,
$w\in\CC$. Then $f$ is a holomorphic function on $\CC$ with values in
$\End(V_\pi)$.

Let $c(\th)=\left(\begin{smallmatrix}e^{-i\th /3}&0&0 \\0&e^{-i\th /3}&0\\
0&0&e^{2i\th /3}\end{smallmatrix}
\right)\in K$. Then $p(c(\th)a)=\left(e^{-i\th}x ,\, e^{-i\th }y,\, 1
\right)$ and
$$f(e^{i\th})=F(e^{-i\th}x, e^{-i\th}y)=F(c(\th)a).$$

We have
\begin{equation*}
\begin{split}
F(c(\th)a)
&= \left(\int_K \chi_{\pi}(k^{-1})\Ph_\pi(gk c(\th)a)\, dk \right)
\Ph_\pi(c(\th)a)^{-1} \displaybreak[0]\\ & =\left(\int_K
\chi_{\pi}(c(\th)k^{-1})\Ph_\pi(gka)\,dk\right) \Ph_\pi(a)^{-1}
\pi(c(\th))^{-1}.
\end{split}
\end{equation*}

Now we notice that $c(\th)$ is in the center of $K$, thus $\pi(c(\th))$ is
a scalar (Schur's lemma) and
$\chi_{\pi}(c(\th)k)=\pi(c(\th))\chi_{\pi}(k)$. Therefore
$$f(e^{i\th})=F(c(\th)a)=F(a),$$
$$F(x,y)=f(1)=f(0)=F(0,0),$$
for all $(x,y)\in \CC^2$. Therefore $F(a)=F(e)$ for all $a\in {\mathcal
A}$, hence
$$\left(\int_K \chi_{\pi}(k^{-1})\Ph_\pi(gka)\, dk \right) \Ph_\pi(a)^{-1}
= \int_K \chi_{\pi}(k^{-1})\Ph_\pi(gk)\, dk=
\Ph_\pi(g)$$ since $\int_K \chi_{\pi}(k^{-1})\pi(k)\, dk=I$ (orthogonality
relations). So, we have proved that
$$\Ph_\pi(g)\Ph_\pi(a)=\int_K \chi_{\pi}(k^{-1})\Ph_\pi(gka)\,dk,$$ for any
$g\in G$ and all $a\in{\mathcal A}$.
Since ${\mathcal A}$ is dense in $G$ and $\Ph_\pi$ is continuous the
theorem is proved.
\end{proof}

\begin{prop}
Given any pair $(G,K)$ and a spherical function $\Ph$ on $G$ of type
$\pi\in \hat K$ the function $\Ph^*:G\longrightarrow
\End{V^*}$ defined by $\Ph^*(g)=\Ph(g^{-1})^T$ is spherical of type
$\pi^*$, where $\pi^*$ denotes the contragradient
representation of $\pi$.
\end{prop}
\begin{proof}
The assertion follows directly from Definition \ref{defesf}. \end{proof}

\begin{cor}\label{pi*} If $\Ph$ is a spherical function of the pair
$(\SU(3),\U(2))$ of type $(n,\ell)$, then $\Ph^*$ is
spherical of type $(-n-\ell,\ell)$.
\end{cor}
\begin{proof} The only thing we have to observe is that $\pi_{n,\ell} ^*$
is equivalent to $\pi_{-n-\ell,\ell}$.
\end{proof}

\medskip
Therefore to determine all spherical functions of type $(n,\ell)$ for all
$n\in \ZZ$ it is enough to determine all spherical
functions of type $(n,\ell)$ with $n\geq -\frac \ell 2$. Observe that for
$n\geq 0$ or $n\leq -\ell$ we have exhibited a
particular spherical function of type $(n,\ell)$: precisely $\Ph_{n,\ell}$
or $\Ph^*_{-n-\ell,\ell}$, respectively.

It is useful to keep in mind the symmetry $n\mapsto -n-\ell$ and the
partition of the integers in the form
$$\ZZ=\vzm{n}{n\geq 0}\cup \vzm{n}{-\ell<n< 0}\cup \vzm{n}{n\leq -\ell}.$$
The spherical functions $\Phi_{n,\ell}$ correspond to taking an irreducible
representation $U$ of $G$ and projecting to the
subrepresentation of $K$ generated by a highest weight vector of $U$,
precisely when $n\ge0$. Moreover the spherical
functions
$\Phi_{n,\ell}^*$ correspond to taking an irreducible representation $U$ of
$G$ and projecting to the subrepresentation of $K$
generated by a lowest weight vector of $U$, exactly when $n\le -\ell$. In
Section \ref{aeigenfunc} this partition of $\ZZ$ will
appear again. Notice that the exceptional interval $\vzm{n}{-\ell<n< 0}$
ocurrs only when $\ell\ge2$.

To exhibit explicitely a spherical function
$\Psi_{n,\ell}$ on
$G$ of type
$(n,\ell)$ for all
$n<0$, we need to introduce the following notation. If $T$ and $S$ are
linear operators on a finite dimensional vector space
$V$ we shall denote by
$T^{\ell-i}\cdot S^i$ the linear map induced by $\ell-i$ factors equal to
$T$ and $i$ factors equal to $S$ acting on the space
of symmetric tensors of rank
$\ell$. Also for $g\in G$ we let $a(g)=g_{33}$ and
$$B(g)=\matc{g_{13}g_{31}\,}{\,g_{13}g_{32}}{g_{23}g_{31}\,}
{\,g_{23}g_{32}}.$$ The proof of the following theorem will be based on
Proposition \ref{defeq} and will be given at the end
of Section \ref{secc1var}. This way of proving is needed because the
function $\Psi_{n,\ell}$ does not satisfy a property
like the one given in Lemma \ref{holomorfia} iii) which allows one to
establish Theorem \ref{Phipiesf}.

\begin{thm}\label{esfneg}
For $n\leq0$
$$\Psi_{n,\ell}(g)=\sum_{0\leq i \leq \min\vz{-n,\ell}} \binom{-n}{i}
a(g)^{-n-i}A(g)^{\ell-i}\cdot B(g)^i$$
is an irreducible spherical function of type $(n,\ell)$. \end{thm}
\noindent Notice that $\Psi_{0,\ell}= \Phi_{0,\ell}$.

\section{Preliminaries}\label{generalidades} The Lie algebra of $G$ is
$\lieg=\vzm{X\in
{\liegl}(3,\CC)} {X=-\overline X^T,\,\tr X=0 }.$ Its complexification is
${\mathfrak g}_{\CC}={\liesl}(3,\CC)$. The Lie
algebra $\liek$ of $K$ can be identified with ${\lieu}(2)$ and its
complexification $\liek_\CC$ with ${\liegl}(2,\CC)$.

The following matrices form a basis of $\lieg$.

$$H_1=\left[ \begin{smallmatrix}
i&0&0 \\ 0&-i&0\\0&0&0\end{smallmatrix}\right], \quad H_2=\left[
\begin{smallmatrix}
i&0&0 \\ 0&i&0\\0&0&-2i\end{smallmatrix}\right],\quad Y_1=\left[
\begin{smallmatrix}
0&1&0 \\ -1&0&0\\0&0&0\end{smallmatrix}\right], \quad Y_2=\left[
\begin{smallmatrix}
0&i&0 \\ i&0&0\\0&0&0\end{smallmatrix}\right],\;$$

$$\;Y_3=\left[ \begin{smallmatrix}
0&0&1 \\ 0&0&0\\-1&0&0\end{smallmatrix}\right], \quad Y_4=\left[
\begin{smallmatrix}
0&0&i \\ 0&0&0\\i&0&0\end{smallmatrix}\right] ,\qquad Y_5=\left[
\begin{smallmatrix}
0&0&0 \\ 0&0&1\\0&-1&0\end{smallmatrix}\right], \quad Y_6=\left[
\begin{smallmatrix}
0&0&0 \\ 0&0&i\\0&i&0\end{smallmatrix}\right]. $$

Let $\lieh$ be the Cartan subalgebra of $\lieg_\CC$ of all diagonal
matrices. The corresponding root space structure is  given
by

$$X_\al=\left[ \begin{smallmatrix}
0&1&0 \\ 0&0&0\\0&0&0\end{smallmatrix}\right], \quad X_{-\al}=\left[
\begin{smallmatrix}
0&0&0 \\ 1&0&0\\0&0&0\end{smallmatrix}\right],\quad H_\al=\left[
\begin{smallmatrix}
1&0&0 \\ 0&-1&0\\0&0&0\end{smallmatrix}\right], $$

$$X_\beta=\left[ \begin{smallmatrix}
0&0&0\\0&0&1 \\0&0&0\end{smallmatrix}\right], \quad X_{-\beta}=\left[
\begin{smallmatrix}
0&0&0 \\ 0&0&0\\0&1&0\end{smallmatrix}\right],\quad H_\beta=\left[
\begin{smallmatrix}
0&0&0 \\ 0&1&0\\0&0&-1\end{smallmatrix}\right], $$

$$X_\gamma=\left[ \begin{smallmatrix}
0&0&1\\0&0&0 \\0&0&0\end{smallmatrix}\right], \quad X_{-\gamma}=\left[
\begin{smallmatrix}
0&0&0 \\ 0&0&0\\1&0&0\end{smallmatrix}\right],\quad H_\gamma=\left[
\begin{smallmatrix}
1&0&0 \\ 0&0&0\\0&0&-1\end{smallmatrix}\right], $$
where $$\al(x_1E_{11}+x_2E_{22}+x_3E_{33})=x_1-x_2,$$
$$\beta(x_1E_{11}+x_2E_{22}+x_3E_{33})=x_2-x_3,$$
$$\g(x_1E_{11}+x_2E_{22}+x_3E_{33})=x_1-x_3.$$

\noindent We have
\begin{align*}
X_\al&=\textstyle\frac 12 (Y_1-iY_2),\quad & X_\beta&=\textstyle\frac 12
(Y_5-iY_6) ,\quad  & X_\g&=\textstyle\frac 12
(Y_3-iY_4), \\ X_{-\al}&=-\textstyle\frac 12 (Y_1+iY_2), \quad &
X_{-\beta}&=-\textstyle\frac 12 (Y_5+iY_6), \quad &
X_{-\g}&=-\textstyle\frac 12 (Y_3+iY_4). \end{align*}

\noindent Let $Z=H_\al+2H_\beta$, $\tilde H_1=2H_\al+H_\beta$ and $\tilde
H_2=H_\beta-H_\al$.

\

Now we shall prove the following lemma stated in Section 2.

\begin{lemsec1}
Given $H\in C^\infty (\CC^2)$ denote also by $H\in C^\infty ({\mathcal A})$
the function defined by $H(g)=H(p(g))$, $g\in
{\mathcal A}$. Then we have $$(X_\beta H)(g)= { \frac{\overline
g_{11}}{g_{33}^2}}\; \frac{\partial H}{\partial
y}-{\frac{\overline g_{21}}{g_{33}^2}}\; \frac{\partial H}{\partial x}\, ,
\qquad (X_\g H)(g)=- { \frac{\overline
g_{12}}{g_{33}^2}}\; \frac{\partial H}{\partial y}+{\frac{\overline
g_{22}}{g_{33}^2}}\; \frac{\partial H}{\partial x}.$$
\end{lemsec1}
\begin{proof}
$2X_\beta=Y_5-iY_6$ and $2X_\g=Y_3-iY_4$. We have \begin{align*} p(g\exp
tY_5)&= \left(
\frac{g_{12}\tan t+g_{13}}{g_{32}\tan t +g_{33}}\, , \, \frac{g_{22}\tan
t+g_{23}}{g_{32}\tan t +g_{33}}\, ,\,1 \right)\\
&\displaybreak[0]\\ & =(u(t), \, v(t), \, 1)=(u_1(t)+iu_2(t),\,
v_1(t)+iv_2(t),\,1). \end{align*} If we put $\dot
u_j=\left(\frac{d u_j}{dt}\right)_{t=0}$ and $\dot v_j=\left(\frac{d
v_j}{dt}\right)_{t=0}$ for $j=1,2$, then
\begin{align*}
Y_5(H)(g)&=
\left(\frac{d}{dt} H \left(p(g\exp t Y_5 )\right) \right)_{t=0}
=H_{x_1}\dot u_1+H_{x_2}\dot u_2+H_{y_1}\dot v_1 + H_{y_2}\dot
v_2. \end{align*} Also by Lemma \ref{menores} we have,
$$ \dot u= -\frac{\overline g_{21}}{g_{33}^2 }\, ,\quad \dot
v=\frac{\overline g_{11}}{g_{33}^2 }. $$

\noindent Similarly, for $Y_6$ let
$$p(g\exp tY_6) =(\tilde u(t), \tilde v(t), \, 1)=(\tilde u_1(t)+i\tilde
u_2(t), \tilde v_1(t)+i\tilde v_2(t),1).$$ Then
\begin{align*} Y_6(H)(g)& = H_{x_1}\dot {\tilde u}_1+H_{x_2}\dot {\tilde
u}_2+H_{y_1}\dot {\tilde v}_1+H_{y_2}\dot {\tilde
v}_2,
\end{align*}
where
$$ \dot {\tilde u}= -\frac{i\,\overline g_{21}}{g_{33}^2 }\, ,\quad \dot
{\tilde v}=\frac{i\,\overline g_{11}}{g_{33}^2 }.$$

\noindent If $u(t)=u_1(t)+iu_2(t)$ then $\dot u_1=\re(\dot u)$ and $\dot
u_2=\im(\dot u).$
Therefore
\begin{equation}\begin{split}
2X_\beta(&H)(g)
\\&= H_{x_1}(\dot u_1-i\dot {\tilde u}_1) +H_{x_2}(\dot u_2-i\dot {\tilde
u}_2)
+H_{y_1}(\dot v_1-i\dot {\tilde v}_1)
+H_{y_2}(\dot v_2-i\dot {\tilde v}_2)\displaybreak[0]\\ &=
-H_{x_1}\,\frac{\overline g_{21}}{g_{33}^2} +iH_{x_2}\,
\frac{\overline g_{21}}{g_{33}^2} + H_{y_1}\,\frac{\overline
g_{11}}{g_{33}^2} -iH_{y_2}\,\frac{\overline g_{11}}{g_{33}^2}\\
&=-2 \frac{\overline g_{21}}{g_{33}^2}\frac{\partial H}{\partial x}
+2\frac{\overline g_{11}}{g_{33}^2}\frac{\partial
H}{\partial y}. \end{split}\end{equation} We proceed in the same way with
$(X_\g H)(g)$ and complete the proof of the
lemma.
\end{proof}

\smallskip
\begin{prop} \label{defDelta}
$D(G)^G$ as a polynomial algebra is generated by \begin{align*}
\Delta_2&= -H_\al^2-\textstyle\frac 13
Z^2-2H_\al-2Z-4X_{-\al}X_\al-4X_{-\beta}X_{\beta}-4X_{-\g}X_{\g} \\
\text{and}\quad\;&\\
\Delta_3&= \textstyle\frac 89 H_\al^3-\frac 89 H_\beta^3+\frac 43
H_\al^2H_\beta-\frac 43 H_\al H_\beta^2+ 8 H_\al^2 +  4H_\al
H_\beta +16H_\al +8H_\beta \\ &\textstyle + 4X_{-\al}X_{\al}H_\al +
8X_{-\al}X_{\al}H_\beta +24 X_{-\al}X_{\al}+\textstyle 12
\left(X_{-\beta}X_{\beta} + X_{-\gamma}X_\gamma\right)
\displaybreak[0]\\
&\,-4
\textstyle X_{-\beta}X_{\beta}\,\tilde H_1 - 4X_{-\gamma}X_{\gamma}\,\tilde
H_2+12X_{-\beta}X_\gamma X_{-\al}+12X_{-\gamma}X_{\beta}X_{\al}.
\end{align*}
\end{prop}
\begin{proof} Since $G$ is a connected Lie group, to verify that $\Delta_2$
and
$\Delta_3$ are in $D(G)^G$ it is enough to check that they are of weight
zero and that
$\mathrm{ad }(X_\alpha)(\Delta_j)=\mathrm{ad }(X_\beta)(\Delta_j)=0$ for
$j=2,3$.
This can be easily accomplished.
To prove that 1, $\Delta_2$ and $\Delta_3$ are algebraically independent
generators
of $D(G)^G$ one can use the well known Harish-Chandra isomorphism
$\xi:D(G)^G\longrightarrow S(\lieh)^W$ where $S(\lieh)^W$
denotes the Weyl group invariants in the symmetric algebra of $\lieh$. (See
\cite{Wa}, Section 3.2.2). We have
$$ \xi(\Delta_2)=-H_\al^2-\textstyle \frac 13 Z^2+4,\qquad \xi(\Delta_3) =
-\textstyle \frac 19 Z^3+ H_\al^2Z+Z^2+3 H_\al -36.$$
If we put $X_1=\frac {\sqrt 3} 3 Z$ and $X_2=H_\al$ then one can verify
that $S(\lieh)^W$ is
generated by the algebraically independent elements 1, $ X_1^2+X_2^2$, and
$X_1(X_1^2-3X_2)$.
(See Proposition 4 in \cite{T2}). This completes the proof of the
proposition.
\end{proof}

\

We write the operators $\Delta_2$ and $\Delta_3$ in the form
$$\Delta_2=\Delta_{2,K}+\tilde\Delta_2,\qquad \Delta_3
=\Delta_{3,K}+\tilde\Delta_3$$
where
\begin{align*}
\Delta_{2,K}& = -H_\al^2-\frac 13 Z^2-2H_\al-2Z-4X_{-\al}X_\al \in
D(K)^K,\displaybreak[0]\\
\tilde\Delta_2&=-4(X_{-\beta}X_{\beta}+X_{-\g}X_{\g})\in
D(G)^K,\displaybreak[0]\\
\Delta_{3,K}&= \textstyle\frac 89 H_\al^3-\frac 89 H_\beta^3+\frac 43
H_\al^2H_\beta-\frac 43 H_\al H_\beta^2+  8 H_\al^2 +
4H_\al H_\beta +16 H_\al +8H_\beta \\ &\quad\textstyle +
4X_{-\al}X_{\al}H_\al
+ 8X_{-\al}X_{\al}H_\beta +24
X_{-\al}X_{\al}\in D(K)^K,\displaybreak[0]\\
\tilde\Delta_3&= \textstyle 12 \left(X_{-\beta}X_{\beta}+
X_{-\gamma}X_\gamma\right)- \textstyle 4X_{-\beta}X_{\beta}\,  \tilde
H_1 - 4X_{-\gamma}X_{\gamma}\,\tilde H_2\displaybreak[0]\\ &\quad +12
X_{-\beta}X_\gamma X_{-\al}+12 X_{-\gamma}X_{\beta}X_{\al}
\in D(G)^K.
\end{align*}

\section{Reduction to $P_2(\CC)$}\label{reduction}

We want to determine all spherical functions $\Phi:G\longrightarrow
\End(V_\pi)$ of type $\pi=\pi_{n,\ell}$. For this purpose we
use the function $\Phi_\pi\in C^{\infty}({\mathcal A})\otimes \End(V_\pi)$
defined by $\Phi_\pi(g)=\pi(A(g))$, in the
following way: in the open set ${\mathcal A}$ we define a function $H$ by
$$H(g)=\Phi(g)\, \Phi_\pi(g)^{-1},$$
where $\Phi$ is supposed to be a spherical function of type $\pi$.
\noindent Then $H$ satisfies
\begin{enumerate}
\item [i)] $H(e)=I$.
\item [ii)]$ H(gk)=H(g)$, for all $g\in {\mathcal A}, k\in K$. \item [iii)]
$H(kg)=\pi(k)H(g)\pi(k^{-1})$, for all
$g\in {\mathcal A}, k\in K$.
\end{enumerate}

\smallskip
The projection map $p:G\longrightarrow P_2(\CC)$ defined before maps the
open set ${\mathcal A}$ onto the affine plane
$\CC^2=\vzm{(x,y,1)\in P_2(\CC)}{(x,y)\in \CC^2}$. Thus ii) says that $H$
may be considered as a function on $\CC^2$.

The fact that $\Phi$ is an eigenfunction of $\Delta_2$ and $\Delta_3$,
makes $H$ into an eigenfunction of certain
differential operators on $\CC^2$, to be determined now.

In the open set ${\mathcal A} \subset G$ let us consider the following
differential operators. For
 $H\in C^{\infty}({\mathcal A})\otimes \End(V_\pi) $ let
$$D(H)=D_1(H)+D_2(H),\qquad \quad E(H)=E_1(H)+E_2(H),$$ where
\begin{align*}
D_1(H)=&-4 (X_{-\beta}X_{\beta}+X_{-\g}X_{\g})(H),\displaybreak[0]\\
D_2(H)=&-4\left(X_{\beta}(H)X_{-\beta}
(\Phi_\pi)\Phi_\pi^{-1} +X_{\g}(H)X_{-\g}(\Phi_\pi)\Phi_\pi^{-1}
\right),\displaybreak[0]\\ &\\
E_1(H)=& -4\textstyle
\left(X_{-\beta}X_{\beta}(H)\right)\,\Ph_\pi\, \dot\pi(\tilde
H_1)\Ph_\pi^{-1} - 4\left(X_{-\gamma}X_{\gamma}(H)\right)\,
\Ph_\pi\,\dot\pi(\tilde H_2)\Ph_\pi^{-1} \\
+ &12
\left(X_{-\beta}X_{\g}(H)\right)\,\Ph_\pi \dot\pi(X_{-\al})\Ph_\pi^{-1}
+12\left(X_{-\gamma}X_{\beta}(H)\right)\,\Ph_\pi\dot\pi\left( X_\al
\right)\Ph_\pi^{-1},\displaybreak[0]
\\ & \\
E_2(H)=& -4X_{\beta}(H)X_{-\beta}(\Ph_\pi)\, \dot\pi(\tilde H_1)\Ph_\pi^{-1}
-4X_{\gamma}(H)X_{-\gamma}(\Ph_\pi)
\,\dot\pi(\tilde H_2)\Ph_\pi^{-1}
\\&
+12X_{\g}(H)X_{-\beta}(\Ph_\pi)
\dot\pi(X_{-\al})\Ph_\pi^{-1}
+12X_{\beta}(H)X_{-\gamma}(\Ph_\pi)
\dot\pi(X_\al)\Ph_\pi^{-1}.
\end{align*}

\begin{prop}\label{lambdamu}
For $H\in C^{\infty}({\mathcal
A})\otimes
\End(V_\pi)$ right invariant under $K$, the function $\Phi=H\Phi_\pi$ is an
eigenfunction of $\Delta_2$ and  $\Delta_3$ on
${\mathcal A}$ if and only if $H$ is an eigenfunction of $D$ and $E$.
Moreover, if $\tilde\lambda$, $\tilde\mu$, $\lambda$ and
$\mu$ denote respectively the corresponding eigenvalues of $\Delta_2$,
$\Delta_3$, $D$ and $E$ then
$$\lambda=\tilde\lambda-\dot\pi(\Delta_{2,K}),\quad\mu=\tilde\mu+\textstyle
3\lambda-
\dot\pi(\Delta_{3,K}).$$\end{prop}
\begin{proof}
If $X\in \liek$ then $X(H)=0$. So, $(XY)(H\Phi_\pi)= H\,(XY)(\Phi_\pi)$,
for all $X,Y\in \liek$. More generally
$\Delta(H\Phi_\pi)=H\,
\Delta(\Phi_\pi)=H\Phi_\pi\,\dot\pi\left(\Delta\right)$, for all $\Delta
\in D(K)$.

\noindent Since $X_{\beta}(\Phi_\pi)=X_{\g}(\Phi_\pi)=0$ we obtain
\begin{align*}
(X_{-\beta}&X_{\beta}+X_{-\g}X_{\g})(H\Phi_\pi)\\
&=(X_{-\beta}X_{\beta}+X_{-\g}X_{\g})(H) \Phi_\pi + X_{\beta}(H)
X_{-\beta}(\Phi_\pi )+ X_{\g}(H) X_{-\g}(\Phi_\pi ).
\end{align*}

\noindent In this way
\begin{equation}\label{DeltaD}
\Delta_2(H\Phi_\pi)=(H\Phi_\pi)\dot\pi\left(\Delta_{2,K}\right)+
D(H)\Phi_\pi.
\end{equation}

\noindent Similarly for $\Delta_3$ we have
\begin{equation}\label{DeltaE}
\begin{split}
\Delta_3(H\Ph_\pi)
=& (H\Ph_\pi) \dot\pi\left(\Delta_{3,K}\right) -
4\left(X_{-\beta}X_{\beta}(H\Ph_\pi)\right)\,
\dot\pi(\tilde H_1) \displaybreak[0]\\
&- 4\left(X_{-\gamma}X_{\gamma}(H\Ph_\pi)\right)\,\dot\pi(\tilde
H_2)+12\left( X_{-\beta}
X_{\g}(H\Ph_\pi)\right)\, \dot\pi(X_{-\al})\\& +12
\left(X_{-\gamma}X_{\beta}(H\Ph_\pi)\right)\, \dot\pi\textstyle(X_{\al}) +
12
\left(X_{-\beta}X_{\beta}+ X_{-\gamma}X_\gamma\right)(H\Ph_\pi) \\ = &
(H\Ph_\pi)\dot\pi\left(\Delta_{3,K}\right) + E(H)
\Phi_\pi-\textstyle 3D(H)\Phi_\pi.
\end{split}
\end{equation}

By Schur's lemma $\dot\pi\left(\Delta_{2,K}\right)$ and
$\dot\pi\left(\Delta_{3,K}\right)$ are scalar because  $\Delta_{2,K}$
and $\Delta_{3,K}\in D(K)^K$. Now it is clear that
$\Delta_2(H\Phi_\pi)=\tilde \ld (H\Phi_\pi)$ and
$\Delta_3(H\Phi_\pi)=\tilde \mu (H\Phi_\pi)$ if and only if $D(H)=\ld H$
and $E(H)=\mu H$ with
$$\ld=\tilde\ld-\dot\pi\left(\Delta_{2,K}\right)\quad \text{ and }\; \quad
\mu=\tilde
\mu-\dot\pi\left(\Delta_{3,K}\right)+\textstyle 3\ld .$$ \end{proof}

Given $H\in C^\infty(\CC^2)\otimes\End(V_\pi)$ we shall also denote by $H\in
C^\infty({\mathcal A})\otimes\End(V_\pi)$ the function defined by
$H(g)=H(p(g))$, $g\in
\mathcal A$. Moreover, if $F$ is a linear endomorphism of
$C^{\infty}({\mathcal A})\otimes\End(V_\pi)$ which preserves the subspace
$C^{\infty}({\mathcal A})^K\otimes\End(V_\pi)$ of all functions which are
right invariant by elements in $K$, then we shall also denote by $F$ the
endomorphism of $C^\infty(\CC^2)\otimes\End(V_\pi)$ which satisfies
$F(H)(p(g))=F(H)(g)$, $g\in\mathcal A$, $H\in
C^\infty(\CC^2)\otimes\End(V_\pi)$.

\begin{lem} The differential operators $D_j$ and $E_j$ $(j=1,2)$ introduced
above, define differential operators $D_j$ and $E_j$ acting on
$C^\infty(\CC^2)\otimes\End(V_\pi)$.
\end{lem}
\begin{proof}
The only thing we really need to proving is that $D_j$ and $E_j$ $(j=1,2)$
preserve the subspace $C^{\infty}({\mathcal A})^K\otimes\End(V_\pi)$.

It is easy to see that $D_1=-4(X_{-\beta}X_\beta+X_{-\gamma}X_\gamma)\in
D(G)^K$.
Then $D_1$ preserves $C^{\infty}({\mathcal A})^K\otimes\End(V_\pi)$. {}From
\eqref{DeltaD}
we get
$$D(H)=\Delta_2(H\Phi_\pi)\Phi_\pi^{-1}-\dot\pi(\Delta_{2,K})H.$$
Since $\Delta_2(H\Phi_\pi)\Phi_\pi^{-1}\in C^{\infty}({\mathcal
A})^K\otimes\End(V_\pi)$ it
follows that $D$ and therefore $D_2$ preserve $C^{\infty}({\mathcal
A})^K\otimes\End(V_\pi)$.

Let us now check that $E_2$ has the same property. Since $K$ is connected
this is equivalent to verifying that for any $X\in\liek$ and all
$H\in C^{\infty}({\mathcal A})^K\otimes\End(V_\pi)$ we have $XE_2(H)=0$,
which in turns amounts to prove that
\begin{align*}
&-XX_{\beta}(H)X_{-\beta}(\Ph_\pi)\dot\pi(\tilde H_1)-
X_{\beta}(H)XX_{-\beta}(\Ph_\pi)\dot\pi(\tilde H_1)\displaybreak[0]\\
&+X_{\beta}(H)X_{-\beta}(\Ph_\pi)\dot\pi(\tilde H_1)\dot\pi(X)-
XX_{\gamma}(H)X_{-\gamma}(\Ph_\pi)\dot\pi(\tilde H_2)\displaybreak[0]\\
&-X_{\gamma}(H)XX_{-\gamma}(\Ph_\pi)\dot\pi(\tilde H_2)+
X_{\gamma}(H)X_{-\gamma}(\Ph_\pi)\dot\pi(\tilde
H_2)\dot\pi(X)\displaybreak[0]\\
&+3XX_{\g}(H)X_{-\beta}(\Ph_\pi)\dot\pi(X_{-\al})+
3X_{\g}(H)XX_{-\beta}(\Ph_\pi)\dot\pi(X_{-\al})\displaybreak[0]\\
&-3X_{\g}(H)X_{-\beta}(\Ph_\pi)\dot\pi(X_{-\al})\dot\pi(X)+
3XX_{\beta}(H)X_{-\gamma}(\Ph_\pi)\dot\pi(X_\al)\displaybreak[0]\\
&+3X_{\beta}(H)XX_{-\gamma}(\Ph_\pi)\dot\pi(X_\al)-
3X_{\beta}(H)X_{-\gamma}(\Ph_\pi)\dot\pi(X_\al)\dot\pi(X)=0.
\end{align*}

Because $X(H)=0$ and $X(\Phi_\pi)=\Phi_\pi\dot\pi(X)$, this is also
equivalent
to showing that
\begin{equation}\label{reduccion}
\begin{split}
&-[X,X_{\beta}](H)X_{-\beta}(\Ph_\pi)\dot\pi(\tilde H_1)-
X_{\beta}(H)[X,X_{-\beta}](\Ph_\pi)\dot\pi(\tilde H_1)\displaybreak[0]\\
&+X_{\beta}(H)X_{-\beta}(\Ph_\pi)\dot\pi([\tilde H_1,X])-
[X,X_{\gamma}](H)X_{-\gamma}(\Ph_\pi)\dot\pi(\tilde H_2)\displaybreak[0]\\
&-X_{\gamma}(H)[X,X_{-\gamma}](\Ph_\pi)\dot\pi(\tilde H_2)+
X_{\gamma}(H)X_{-\gamma}(\Ph_\pi)\dot\pi([\tilde H_2,X])\displaybreak[0]\\
&+3[X,X_{\g}](H)X_{-\beta}(\Ph_\pi)\dot\pi(X_{-\al})+
3X_{\g}(H)[X,X_{-\beta}](\Ph_\pi)\dot\pi(X_{-\al})\displaybreak[0]\\
&-3X_{\g}(H)X_{-\beta}(\Ph_\pi)\dot\pi([X_{-\al},X])
+3[X,X_{\beta}](H)X_{-\gamma}(\Ph_\pi)\dot\pi(X_\al)\displaybreak[0]\\
&+3X_{\beta}(H)[X,X_{-\gamma}](\Ph_\pi)\dot\pi(X_\al)-
3X_{\beta}(H)X_{-\gamma}(\Ph_\pi)\dot\pi([X_\al,X])=0.
\end{split}
\end{equation}

\noindent If we put $X=T\in\lieh$ in \eqref{reduccion} and use
$\gamma(T)=\alpha(T)+\beta(T)$ we get
\begin{align*}
&-\beta(T)X_{\beta}(H)X_{-\beta}(\Ph_\pi)\dot\pi(\tilde H_1)+
\beta(T)X_{\beta}(H)X_{-\beta}(\Ph_\pi)\dot\pi(\tilde H_1)\displaybreak[0]\\
&-\gamma(T)X_{\gamma}(H)X_{-\gamma}(\Ph_\pi)\dot\pi(\tilde H_2)+
\gamma(T)X_{\gamma}(H)X_{-\gamma}(\Ph_\pi)\dot\pi(\tilde
H_2)\displaybreak[0]\\
&+3\gamma(T)X_{\g}(H)X_{-\beta}(\Ph_\pi)\dot\pi(X_{-\al})-
3\beta(T)X_{\g}(H)X_{-\beta}(\Ph_\pi)\dot\pi(X_{-\al})\displaybreak[0]\\
&-3\alpha(T)X_{\g}(H)X_{-\beta}(\Ph_\pi)\dot\pi(X_{-\al})+
3\beta(T)X_{\beta}(H)X_{-\gamma}(\Ph_\pi)\dot\pi(X_\al)\displaybreak[0]\\
&-3\gamma(T)X_{\beta}(H)X_{-\gamma}(\Ph_\pi)\dot\pi(X_\al)
+3\alpha(T)X_\beta(H)X_{-\gamma}(\Phi_\pi)\dot\pi(X_\alpha)=0
\end{align*}

\noindent If we substitute $X=X_\alpha$ in \eqref{reduccion} and use
$[X_\alpha,X_\beta]=X_\gamma$,
$[X_\alpha,X_{-\gamma}]=-X_{-\beta}$, $[X_\alpha,X_{-\alpha}]=H_\alpha$,
$[X_\alpha,X_\gamma]=[X_\alpha,X_{-\beta}]=0$, we obtain
\begin{align*}
&-X_{\gamma}(H)X_{-\beta}(\Ph_\pi)\dot\pi(\tilde H_1)+
\alpha(\tilde
H_1)X_{\beta}(H)X_{-\beta}(\Ph_\pi)\dot\pi(X_\alpha)\displaybreak[0]\\
&+X_{\gamma}(H)X_{-\beta}(\Ph_\pi)\dot\pi(\tilde H_2)+
\alpha(\tilde
H_2)X_{\gamma}(H)X_{-\gamma}(\Ph_\pi)\dot\pi(X_\alpha)\displaybreak[0]\\
&+3X_{\g}(H)X_{-\beta}(\Ph_\pi)\dot\pi(H_\al)+
3X_{\g}(H)X_{-\g}(\Ph_\pi)\dot\pi(X_\al)\displaybreak[0]\\
&-3X_{\beta}(H)X_{-\beta}(\Ph_\pi)\dot\pi(X_\al)=0
\end{align*}
because $3H_\alpha+\tilde H_2-\tilde H_1=0$, $\alpha(\tilde H_1)=3$ and
$\alpha(\tilde H_2)=-3$.
Now the representation theory of GL$(2,\CC)$ tell us that $XE_2(H)=0$ for
all
$X\in\liek$.

{}From \eqref{DeltaE} we get
$$E(H)=\Delta_3(H\Phi_\pi)\Phi_\pi^{-1}-\big(\dot\pi(\Delta_{3,K})-\textstyle
3D\big)(H).$$
Since $\Delta_3(H\Phi_\pi)\Phi_\pi^{-1}\in C^{\infty}({\mathcal
A})^K\otimes\End(V_\pi)$ it
follows that $E$ and therefore $E_1$ preserve $C^{\infty}({\mathcal
A})^K\otimes\End(V_\pi)$.
This finishes the proof of the lemma.
\end{proof}

\medskip
For $g\in \mathcal A$ we put $p(g)=(x,y,1)$. Then
\begin{equation}\label{expxy}
x=g_{13}g_{33}^{-1}\, , \qquad \,y=g_{23}g_{33}^{-1}\, , \qquad \,
|g_{33}|^{-2}=1+|x|^2+|y|^2.
\end{equation}

The proofs of Propositions \ref{expD1}, \ref{expD2}, \ref{expE1} and
\ref{expE2} require simple but lengthy computations. We
include them for completeness in an Appendix at the end of the paper.

\begin{prop}\label{expD1} For $H\in C^\infty(\CC^2)\otimes \End(V_\pi)$ we
have
\begin{align*}
D_1(&H)(x,y)=\\
(&1+|x|^2+|y|^2) \left( \vphantom{(1+|y|^2)}
\,(H_{x_1x_1}+H_{x_2x_2})(1+|x|^2)+(H_{y_1y_1}+H_{y_2y_2})(1+|y|^2)\right.
\displaybreak[0]\\ & \left. \vphantom{(1+|y|^2)}
\;+2(H_{y_1x_1}+H_{y_2x_2})\,\re(x\overline y)
+2(H_{y_1x_2}-H_{y_2x_1})\,\im(x\overline y) \right).
\end{align*}
\end{prop}

\begin{prop}\label{expD2} For $H\in C^\infty(\CC^2)\otimes \End(V_\pi)$ we
have
\begin{align*}
D_2(&H)(x,y)=\\
&- 4 \,\frac{\partial H}{\partial x}\, \dot \pi\begin{pmatrix}
x(1+|x|^2) & x^2{\overline y}\\ &\\ y(1+|x|^2)& x|y|^2
\end{pmatrix} - 4\,\frac{\partial H}{\partial y}\,
\dot \pi\begin{pmatrix} y|x|^2 & x(1+|y|^2) \\ &\\ y^2{\overline x}&
y(1+|y|^2) \end{pmatrix}.
\end{align*}
\end{prop}

\begin{prop}\label{expE1}
For $H\in C^\infty(\CC^2)\otimes \End(V_\pi)$ we have
\begin{align*}
E_1(H)&(x,y)=\\
&(1+|x|^2+|y|^2) \left[ \vphantom{\dot\pi
\matc{-(1+|x|^2)}{-3x\overline y}{0}{2(1+|x|^2)}} \, (H_{x_1x_1}
+H_{x_2x_2})\;\dot\pi \matc{-(1+|x|^2)}{-3x\overline
y}{0}{2(1+|x|^2)}\right.\\
&\\&\left. +(H_{y_1y_1}+H_{y_2y_2})\;\dot\pi
\matc{2(1+|y|^2)}{0}{-3\overline x y}{-(1+|y|^2)}\right.
\displaybreak[0]\\
&\\ & \left.
+(H_{y_1x_1}+H_{y_2x_2})\;\dot\pi \matc{2x\overline y-\overline
xy}{-3(1+|y|^2)}
{-3(1+|x|^2)}{2\overline x y-x \overline y}\right.\\ &\\ &\left.
+i\,(H_{x_2y_1}-H_{x_1y_2})\;\dot\pi \matc{\,-(2x\overline y
+\overline xy)\;}{\;-3\,(1+|y|^2)} {3\,(1+|x|^2)\;}{\;\,2\overline x y+x
\overline y} \right].
\end{align*}
\end{prop}

\begin{prop}\label{expE2}
For $H\in C^\infty(\CC^2)\otimes \End(V_\pi)$ we have
\begin{align*} &E_2(H)(x,y)=\\
& 4\frac{\partial H}{\partial x}\!\!\left(\!\dot\pi\!\matc{0}{x}{0}{y}
\dot\pi
\!\matc{-2x\overline y }{0}{3(1+|x|^2)}{\;x\overline y}
+ \dot\pi\!\matc{x}{0}{y}{0} \dot\pi\! \matc{1+|x|^2}{3x\overline
y}{0}{-2(1+|x|^2)}\!\right)
\displaybreak[0]\\ &\\
&+4\frac{\partial H}{\partial y}\!\!\left(\!\dot\pi\!\matc{0}{x}{0}{y}
\dot\pi\! \matc{-2(1+|y|^2)}{0} {3y\overline x}
{1+|y|^2}+ \dot\pi\!\matc{x}{0}{y}{0}\dot\pi\!\matc{y\overline x
}{\,3(1+|y|^2)}{0}{-2y\overline x}\!\right). \end{align*}
\end{prop}

\section{Reduction to one variable }\label{secc1var}
We are interested in considering the differential operators $D$ and $E$
applied to a function $H\in C^\infty(\CC^2)\otimes
\End(V_\pi)$ such that $H(kp)=\pi(k)H(p)\pi(k)^{-1}$, for all $k\in K$ and
$p$ in the affine complex plane $\CC^2$. This
property of $H$ allows us to find ordinary differential operators $\tilde
D$ and $\tilde E$ defined on the interval
$(0,\infty)$ such that $$ (D\,H)(r,0)=(\tilde D\tilde H)(r),\qquad
(E\,H)(r,0)=(\tilde E\tilde H)(r),$$ where $\tilde
H(r)=H(r,0)$. We also define differential operators $\tilde D_1$, $\tilde
D_2$, $\tilde E_1$ and $\tilde E_2$ in the same way,
that is $(D_1 H)(r,0)=(\tilde D_1\tilde H)(r),\; (D_2 H)(r,0)=(\tilde
D_2\tilde H)(r),\;(E_1 H)(r,0)=(\tilde E_1\tilde H)(r),
\; (E_2 H)(r,0)=(\tilde E_2\tilde H)(r).$

\medskip
The goal of this section is to prove the following theorems.
\begin{thm}\label{mainth1}
We have
\begin{align*}
\tilde D(\tilde H)(r)&=(1+r^2)^2\frac{d^2 \tilde H}{dr^2}
+\frac{(1+r^2)}{r}\biggl(3+r^2-2r^2\dot\pi(H_\g)\biggr)
\frac{d\tilde H}{dr}\displaybreak[0]
\\ & \quad
+\frac {(1+r^2)}{r^2}\left( \dot\pi(J)^2\,\tilde H(r)\,+\tilde
H(r)\dot\pi(J)^2 -2\dot\pi(J)\tilde H(r)\dot\pi(J)\right)\\ &
\quad -\frac {(1+r^2)}{r^2}\left( \dot\pi(T)^2\,\tilde H(r)\,+\tilde
H(r)\dot\pi(T)^2 -2\dot\pi(T)\tilde H(r)\dot\pi(T)\right)
\displaybreak[0] \\ &\quad -4\dot\pi(X_{-\al})\tilde
H(r)\dot\pi(X_{\al})+4\tilde H(r)\dot\pi(X_{-\al})\dot\pi(X_{\al}),
\end{align*}
where $J=\begin{pmatrix}{0}&{1}\\{-1}&{0}\end{pmatrix}$ and
$T=\begin{pmatrix}{0}&{1}\\{1}&{0}\end{pmatrix}$.
\end{thm}

\begin{thm} \label{mainth2} We have
\begin{align*}
&\tilde E(\tilde H)(r)\\&
 =(1+r^2)^2 \,\frac{d^2\tilde H}{d\,r^2}\dot\pi(\tilde H_2)
+\frac{(1+r^2)^2}{r}\,\frac{d\tilde H}{dr}\;\dot\pi(\tilde H_2)
+\frac{2(1+r^2)}{r} \frac{d \tilde H}{d r}\;\dot\pi (\tilde H_1)
\displaybreak[0]
\\&\quad
+ 2r(1+r^2)\,\frac{d\tilde H}{dr}\left(3\, \dot\pi(X_{\al})\dot\pi
(X_{-\al})-\dot\pi(H_\g) \dot\pi(\tilde H_2) \right)\\ &
\quad +\frac{6(1+r^2)^2 }{r}\left( \dot\pi (X_\al)\frac{d\tilde
H}{d\,r}-\frac{d\tilde H}{dr}\dot\pi(X_\al)\right)\; \dot\pi
(X_{-\al})\\&\quad -\frac{6(1+r^2) }{r}\left( \dot\pi
(X_{-\al})\frac{d\tilde H}{d\,r}-\frac{d\tilde
H}{dr}\dot\pi(X_{-\al})\right)\;\dot\pi (X_{\al})\displaybreak[0]\\ & \quad
+ \frac{(1+r^2)}{r^2}\left(\dot\pi(J)^2\, \tilde
H(r)+ \tilde H(r)\, \dot\pi(J)^2-2\dot\pi(J)\tilde H(r)\,
\dot\pi(J)\right)\; \dot\pi (\tilde H_1)\displaybreak[0]
\\&\quad
- \frac{(1+r^2)}{r^2}\left(\dot\pi(T)^2\, \tilde H(r)+ \tilde H(r)\,
\dot\pi(T)^2-2\dot\pi(T)\tilde H(r)\, \dot\pi(T)\right)\;
\dot\pi (\tilde H_1) \displaybreak[0] \\& \quad
+ 4 \left(\dot\pi(X_{-\al})\, \tilde H(r)- \tilde H(r)\,
\dot\pi(X_{-\al})\right)\,\left(-\textstyle\dot\pi(X_\al) \dot\pi
(\tilde H_1) +3\, \dot\pi(H_\g)\dot\pi(X_\al)\right). \end{align*}
\end{thm}

The proof of these theorems will be a direct consequence of Propositions
\ref{D_1H}, \ref{D_2H}, \ref{E_1H} and \ref{E_2H}.
To reach this goal we need some preparatory material.

\smallskip
Since we are interested in considering
functions $H$ on $P_2(\CC)$ such that $H(kp)=\pi(k)H(p)\pi(k)^{-1}$, for
all $k\in K$ and $p\in P_2(\CC)$ we need to see the
$K$-orbit structure of $P_2(\CC)$. The affine plane $\CC^2$ is $K$-stable
and the corresponding line at infinity
$L=\vzm{(x,y,0)}{x,y\in \CC^2}$ in $P_2(\CC)$ is a $K$-orbit. Moreover the
$K$-orbits in the affine plane correspond to the
spheres
$$S_r=\vzm{(x,y)\in \CC^2}{|x|^2+|y|^2=r^2}, \qquad r\in \RR_{\geq 0}$$ in
$\CC^2$.
In each affine orbit we may choose as a representative the point $(r,0)\in
\CC^2$, with $r\geq 0$. In fact given
$(x,y)\in \CC^2$, there exists $k \in K$ such that $k \, (r,0)=(x,y)$ where
$r^2=|x|^2+|y|^2$. Also we may choose  $(1,0,0)$
as a representative of the line at infinity. Since
$(M,0,1)=(1,0,\textstyle\frac 1M)\longrightarrow (1,0,0)$ when
$M\rightarrow \infty$, the closed interval $[0,\infty]$ parametrizes the
set of $K$-orbits in $P_2(\CC)$. This orbit structure
is depicted in the figure below.
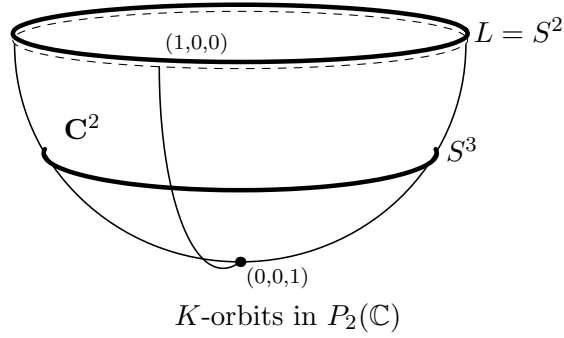
\vspace{1.0cm}
\begin{center}
\includegraphics{proyectivo.1} \\
$K$-orbits in $P_2(\CC)$
\end{center}

\medskip
Turning to the proofs of Propositions \ref{D_1H}, \ref{D_2H}, \ref{E_1H}
and \ref{E_2H} and refering to Propositions
\ref{expD1}, \ref{expD2}, \ref{expE1} and \ref{expE2} we need to compute a
number of first and second order partial
derivatives of the function $H$ at the point $(r,0)\in\CC^2$. These are
given
in Lemmas 5.3 through 5.9 whose proofs are
included in the Appendix.

\begin{lem}\label{dderx1}
At $(r,0)\in \CC^2$ we have
$$H_{x_1}(r,0)=\frac{d\tilde H}{dr}(r) \quad \text{ and }\quad
H_{x_1x_1}(r,0)=\frac{d^2\tilde H}{dr^2}(r).$$
\end{lem}

\begin{lem}\label{ddery1}
At $(r,0)\in \CC^2$ we have
\begin{equation*}
H_{y_1}(r,0)= -{\frac 1r}\left( \dot\pi(J)\tilde H(r)-\tilde
H(r)\dot\pi(J)\right)
\end{equation*}
and
\begin{equation*}
\begin{split}
H_{y_1y_1}(r,0)&= {\frac 1r} \frac{d \tilde H}{d r} +
\frac{1}{r^2}\left(\dot\pi(J)^2\, \tilde H(r)+ \tilde H(r)\,
\dot\pi(J)^2\right) -\frac 2{r^2}\dot\pi(J)\tilde H(r)\, \dot\pi(J),
\end{split}
\end{equation*}
where $J=\matc{0}{1}{-1}{0}$.
\end{lem}

\begin{lem}\label{ddery2}
At $(r,0)\in \CC^2$ we have
\begin{equation*}
H_{y_2}(r,0)= \frac{i}{r}\left( \dot\pi(T)\tilde H(r)-\tilde
H(r)\dot\pi(T)\right)
\end{equation*}
and
\begin{equation*}
\begin{split}
H_{y_2y_2}(r,0)&= {\frac 1r} \frac{d \tilde H}{d r}
-\frac{1}{r^2}\left(\dot\pi(T)^2\, \tilde H(r)+ \tilde H(r)\,
\dot\pi(T)^2\right) +\frac 2{r^2}\dot\pi(T)\tilde
H(r)\, \dot\pi(T),
\end{split}
\end{equation*}
where $T=\matc{0}{1}{1}{0}$.
\end{lem}

\begin{lem}\label{dderx2}
At $(r,0)\in \CC^2$ we have
\begin{equation*}
H_{x_2}(r,0)= {\frac ir}\left( \dot\pi(H_\al)\tilde H(r)-\tilde
H(r)\dot\pi(H_\al)\right)
\end{equation*}
and
\begin{equation*}
\begin{split}
H_{x_2x_2}(r,0)&= {\frac 1r} \frac{d \tilde H}{d r} -
\frac{1}{r^2}\left(\dot\pi(H_\al)^2\, \tilde H(r)+ \tilde H(r)\,
\dot\pi(H_\al)^2\right)\displaybreak[0]\\ & \quad +\frac
2{r^2}\dot\pi(H_\al)\,\tilde H(r)\, \dot\pi(H_\al). \end{split}
\end{equation*}
\end{lem}

\begin{lem}\label{ddery1x1}
At $(r,0)\in \CC^2$ we have
\begin{equation*}
H_{x_1y_1}(r,0)=\! \frac {-1}r \Big(\dot\pi(J)\frac{d \tilde H}{d
r}\!- \!\frac{d \tilde H}{d r}\dot\pi(J)\Big)\!
 + \!\frac{1}{r^2}\left(\dot\pi(J)\tilde H(r)-\tilde H(r)\dot\pi(J)\right),
\end{equation*}
\begin{equation*}
H_{x_1y_2}(r,0)= {\frac ir}\Big(\dot\pi(T)\frac{d \tilde H}{d
r}-\frac{d \tilde H}{d r}\dot\pi(T)\Big) -
\frac{i}{r^2}\left( \dot\pi(T)\tilde H(r)-\tilde H(r)\dot\pi(T)\right).
\end{equation*}
\end{lem}

\begin{lem}\label{derx2y1}
At $(r,0)\in \CC^2$ we have
\begin{equation*}
\begin{split}
H_{y_1x_2}(r,0)=& -
\frac{i}{2r^2}\Big(\dot\pi(H_\al)\,\dot\pi(J)
+\dot\pi(J)\,\dot\pi(H_\al)\Big)\tilde H(r)\displaybreak[0] \\
& - \frac{i}{2r^2}\tilde H(r)\,\Big(\dot\pi(H_\al)\,\dot\pi(J)
+\dot\pi(J)\,\dot\pi(H_\al)\Big)\displaybreak[0]\\ &
+\frac{i}{r^2}\left(\dot\pi(H_\al)\,\tilde H(r)\,\dot\pi(J)
+\dot\pi(J)\,\tilde H(r)\,\dot\pi(H_\al)\right).
\end{split}
\end{equation*}
\end{lem}

\begin{lem}\label{derx2y2}
At $(r,0)\in \CC^2$ we have
\begin{equation*}
\begin{split}
H_{y_2x_2}(r,0)=& -
\frac{1}{2r^2}\Big(\dot\pi(H_\al)\,\dot\pi(T)
+\dot\pi(T)\,\dot\pi(H_\al)\Big)\tilde H(r)\displaybreak[0] \\ & -
\frac{1}{2r^2}\tilde H(r)\,\Big(\dot\pi(H_\al)\,\dot\pi(T)
+\dot\pi(T)\,\dot\pi(H_\al)\Big)\displaybreak[0]\\ &
+\frac{1}{r^2}\left(\dot\pi(H_\al)\,\tilde H(r)\,\dot\pi(T)
+\dot\pi(T)\,\tilde H(r)\,\dot\pi(H_\al)\right).
\end{split}
\end{equation*}
\end{lem}

\medskip
\begin{prop}\label{Hdiagonal}
The function $\tilde H$ is diagonalizable.
\end{prop}
\begin{proof}
Let $M=\vzm{m_a}{|a|=1}$ where
$$m_a=\left(\begin{smallmatrix} a& & \\& a^{-2}& \\&& a
\end{smallmatrix}\right).$$
The subgroup $M$ fixes the points $(x,0,1)$ in the projective plane.
Moreover $H(kg)=\pi(k)H(g)\pi(k^{-1})$, therefore $\tilde
H(r)=\pi(m_a)\tilde H(r)\pi(m_a^{-1})$. Now since all finite dimensional
irreducible $K$-modules are multiplicity free as
$M$-modules, $\tilde H(r)$, $r\ge 0$, and $\dot\pi(H_\alpha)$ diagonalize
simultaneously.
\end{proof}

\begin{prop}\label{D_1H} We have
\begin{align*}
\tilde D_1(\tilde H)(r)&=(1+r^2)^2\frac{d^2\tilde H}{dr^2} +
(1+r^2)(3+r^2)\,\frac{1}{r}\frac{d \tilde H}{d r}
\displaybreak[0]\\ & +\frac {(1+r^2)}{r^2}\left( \dot\pi(J)^2\,\tilde
H(r)\,+\tilde H(r)\dot\pi(J)^2 -2\dot\pi(J)\tilde
H(r)\dot\pi(J)\right) \displaybreak[0]\\ & -\frac {(1+r^2)}{r^2}\left(
\dot\pi(T)^2\,\tilde H(r)\,+\tilde H(r)\dot\pi(T)^2
-2\dot\pi(T)\tilde H(r)\dot\pi(T)\right).
\end{align*}
\end{prop}
\begin{proof}
By Proposition \ref{expD1} we have
\begin{align*}
D_1(H)(r,0)=\Bigl(H_{x_1x_1}&(r,0)+H_{x_2x_2}(r,0)\Bigr)\,(1+r^2)^2
\displaybreak[0]\\
&+\Bigl(H_{y_1y_1}(r,0)+H_{y_2y_2}(r,0)\Bigr)\,(1+r^2).
\end{align*}
Using the lemmas above and the fact that $\tilde H(r)$ and $\dot\pi(H_\al)$
commute because they are simultaneously
diagonalizable, the proposition follows.
\end{proof}

\begin{prop}\label{D_2H} We have
\begin{align*}
\tilde D_2(\tilde H)(r)&=-2r(1+r^2)\frac{d \tilde H}{dr}\, \dot\pi
(H_\g)\displaybreak[0]\\
&\quad
-4\dot\pi(X_{-\al})\tilde H(r)\dot\pi(X_{\al})+4\tilde
H(r)\dot\pi(X_{-\al})\dot\pi(X_{\al}).
\end{align*}
\end{prop}
\begin{proof}
By Proposition \ref{expD2} we have
$$D_2(H)(r,0)=- 4\left( \frac{\partial H}{\partial x}(r,0)\,{\textstyle
\dot \pi \matc{r(1+r^2)}{0}{0}{0} }+\,
\frac{\partial H}{\partial y}(r,0)\, \textstyle\dot
\pi\matc{0}{r}{0}{0}\right).$$

\noindent Now using the lemmas above we compute $\frac{\partial H}{\partial
x}(r,0)$ and $\frac{\partial
H}{\partial y}(r,0)$.
\begin{equation}\label{partialx}
\begin{split}
\frac{\partial H }{\partial x}(r,0)&=\frac 12 \frac{\partial H}{\partial
x_1}(r,0)- \frac i2
\frac{\partial H }{\partial x_2}(r,0)\displaybreak[0]\\ &= \frac 12
\frac{d \tilde H}{dr}+
\frac{1}{2r}\dot\pi(H_\al)\,\tilde H(r) -\frac{1}{2r} \tilde
H(r)\,\dot\pi(H_\al) = \frac 12 \frac{d \tilde H}{dr},
\end{split}
\end{equation}
\begin{equation}\label{partialy}
\begin{split}
\frac{\partial H }{\partial y}(r,0)&=\frac 12 \frac{\partial H}{\partial
y_1}(r,0)- \frac i2
\frac{\partial H }{\partial y_2}(r,0)\displaybreak[0]\\ &= -
\frac{1}{2r}\dot\pi(J)\,\tilde H(r)+\frac{1}{2r} \tilde
H(r)\,\dot\pi(J) +\frac{1}{2r}\dot\pi(T)\,\tilde H(r) -\frac{1}{2r} \tilde
H(r)\,\dot\pi(T) \displaybreak[0]\\ &= \frac{1}{r}
\dot\pi(X_{-\al})\,\tilde H(r)- \frac{1}{r} \tilde H(r)\,\dot\pi(X_{-\al}).
\end{split}
\end{equation}
The proposition follows.
\end{proof}

\medskip

\begin{prop}\label{E_1H}
\begin{equation*}
\begin{split}
\tilde E_1&(\tilde H)(r)\\&
=(1+r^2)^2 \,\frac{d^2\tilde H}{d\,r^2}\dot\pi(\tilde H_2)
+\frac{(1+r^2)^2}{r}\,\frac{d\tilde H}{dr}\;\dot\pi(\tilde H_2)
+\frac{2(1+r^2)}{r} \frac{d \tilde H}{d r}\;\dot\pi (\tilde H_1)
\displaybreak[0]
\\& \quad
+\frac{6(1+r^2)^2 }{r}\Big(\dot\pi (X_\al)\frac{d\tilde
H}{d\,r}-\frac{d\tilde H}{dr}\dot\pi(X_\al)\Big)\;  \dot\pi
(X_{-\al})\\ &\quad -\frac{6(1+r^2) }{r}\Big(\dot\pi
(X_{-\al})\frac{d\tilde H}{d\,r}-\frac{d\tilde
H}{dr}\dot\pi(X_{-\al})\Big)\;\dot\pi (X_{\al})\displaybreak[0]\\ & \quad
+ \frac{(1+r^2)}{r^2}\left(\dot\pi(J)^2\, \tilde
H(r)+ \tilde H(r)\, \dot\pi(J)^2-2\dot\pi(J)\tilde H(r)\,
\dot\pi(J)\right)\; \dot\pi (\tilde H_1)\displaybreak[0]
\\&\quad
-\frac{(1+r^2)}{r^2}\left(\dot\pi(T)^2\, \tilde H(r)+ \tilde H(r)\,
\dot\pi(T)^2-2\dot\pi(T)\tilde H(r)\,\dot\pi(T)\right)\;
\dot\pi (\tilde H_1).
\end{split}
\end{equation*}
\end{prop}
\begin{proof}
By Proposition \ref{expE1} we have
\begin{align*}
E_1(H)(r,0)=&(1+r^2)^2 \Bigl(H_{x_1x_1}(r,0)+H_{x_2x_2}(r,0)\Bigr)
\;\dot\pi(\tilde H_2)\\
&+(1+r^2) \Bigl(H_{y_1y_1}(r,0)+H_{y_2y_2}(r,0)\Bigr)
\;\dot\pi(\tilde H_1)
\displaybreak[0]\\
&-3(1+r^2)\Bigl(H_{y_1x_1}(r,0)+H_{y_2x_2}(r,0)\Bigr)\;\dot\pi
{\textstyle\matc{0}{1}{1+r^2}{0}}
\\ &
+ 3 i (1+r^2)\Bigl(H_{x_2y_1}(r,0)-H_{x_1y_2}(r,0)\Bigr)\;\dot\pi
{\textstyle
\matc{0}{-1}{1+r^2}{0}}.
\end{align*}

\noindent Since, as noticed before, $\tilde H(r)$ and $\dot\pi(H_\al)$
commute, from Lemmas \ref{derx2y1} and
\ref{derx2y2} we get
\begin{align*}
H_{x_2y_1}(r,0)=&
\frac{i}{2r^2}\Bigl(\dot\pi(J)\,\dot\pi(H_\al) -\dot\pi(H_\al)\,\dot\pi(J)
\Bigr)\tilde H(r)\displaybreak[0] \\
& - \frac{i}{2r^2}\tilde
H(r)\,\Bigl(\dot\pi(J)\,\dot\pi(H_\al)-\dot\pi(H_\al)\,\dot\pi(J) \Bigr)
\displaybreak[0] \\
=& \frac{i}{r^2}\tilde H(r)\dot\pi(T)
-\frac{i}{r^2}\dot\pi(T) \tilde H(r).
\end{align*}
\begin{align*}
H_{y_2x_2}(r,0)=&
\frac{1}{2r^2}\Bigl(\dot\pi(T)\,\dot\pi(H_\al) -\dot\pi(H_\al)\,\dot\pi(T)
\Bigr)\tilde H(r)\displaybreak[0] \\
& - \frac{1}{2r^2}\tilde
H(r)\,\Bigl(\dot\pi(T)\,\dot\pi(H_\al)-\dot\pi(H_\al)\,\dot\pi(T)
\Bigr)\displaybreak[0] \\
=& \frac{1}{r^2}\tilde H(r)\dot\pi(J)
-\frac{1}{r^2}\dot\pi(J) \tilde H(r).
\end{align*}

\noindent In this way we obtain
\begin{align}
\tilde E_1&(\tilde H)(r)\nonumber \\
&=(1+r^2)^2 \,\frac{d^2\tilde H}{d\,r^2}\dot\pi(\tilde H_2)
+\frac{(1+r^2)^2}{r}\,\frac{d\tilde H}{dr}\;\dot\pi(\tilde H_2)
+\frac{2(1+r^2)}{r} \frac{d \tilde H}{d r}\;\dot\pi (\tilde H_1)
\displaybreak[0] \nonumber
\\& \quad
+\frac{3(1+r^2) }{r}\Big(\dot\pi (J)\frac{d\tilde H}{d\,r}-\frac{d\tilde
H}{dr}\dot\pi(J)\Big)\;
\dot\pi {\textstyle \matc {0}{1}{1+r^2}{0}} \displaybreak[0] \label{aaa1}\\
&\quad
+\frac{3(1+r^2) }{r}\Big(\dot\pi (T)\frac{d\tilde H}{d\,r}-\frac{d\tilde
H}{dr}\dot\pi(T)\Big)
\;\dot\pi {\textstyle \matc
{0}{-1}{1+r^2}{0}}\displaybreak[0]\label{aaa2}\\ & \quad +
\frac{(1+r^2)}{r^2}\left(\dot\pi(J)^2\, \tilde H(r)+
\tilde H(r)\, \dot\pi(J)^2-2\dot\pi(J)\tilde H(r)\, \dot\pi(J)\right)\;
\dot\pi (\tilde H_1)\displaybreak[0] \nonumber \\
&\quad -\frac{(1+r^2)}{r^2}\left(\dot\pi(T)^2\, \tilde H(r)+ \tilde H(r)\,
\dot\pi(T)^2-2\dot\pi(T)\tilde H(r)\,
\dot\pi(T)\right)\; \dot\pi (\tilde H_1)\nonumber.
\end{align}
We observe that combining lines \eqref{aaa1} and \eqref{aaa2} one gets
\begin{align*}
\eqref{aaa1} &+\eqref{aaa2}\\ &=
\frac{3(1+r^2)^2 }{r} \dot\pi (J+T)\frac{d\tilde H}{d\,r} \dot\pi(X_{-\al})
+\frac{3(1+r^2) }{r} \dot\pi (J-T)
\frac{d\tilde H}{d\,r} \dot\pi(X_{\al})\displaybreak[0]\\
&\quad-\frac{3(1+r^2)^2}{r}\frac{d\tilde H}{d\,r} \dot\pi (J+T)
\dot\pi(X_{-\al})
-\frac{3(1+r^2)}{r} \frac{d\tilde H}{d\,r}\dot\pi (J-T) \dot\pi(X_{\al}).
\end{align*}
By using that $\dot\pi(T+J)=2\dot\pi(X_\al)$ and
$\dot\pi(T-J)=2\dot\pi(X_{-\al})$ the proposition follows.
\end{proof}

\begin{prop}\label{E_2H}
\begin{align*}
&\tilde E_2(\tilde H)(r)=2r(1+r^2)\,\frac{d\tilde
H}{dr}\left(3\, \dot\pi(X_{\al})\dot\pi (X_{-\al})
-\dot\pi(H_\g) \dot\pi(\tilde H_2) \right)
\\& \quad
+ 4\left(\dot\pi(X_{-\al})\, \tilde H(r)- \tilde H(r)\,
\dot\pi(X_{-\al})\right)\,\left(-\textstyle\dot\pi(X_\al) \dot\pi
(\tilde H_1) +3\, \dot\pi(H_\g)\dot\pi(X_\al)\right). \end{align*}
\end{prop}
\begin{proof}
By Proposition \ref{expE2} we have
\begin{align*}
E_2&(H)(r,0)=\\
& 4\frac{\partial H}{\partial x}(r,0)\!\left( \dot\pi\!\matc{0}{r}{0}{0}
\dot\pi\!
\matc{0}{0}{3(1+r^2)}{0}
+ \dot\pi\!\matc{r}{0}{0}{0} \dot\pi\!\matc{1+r^2}{0}{0}{-2(1+r^2)}\!
\right)
\displaybreak[0]\\ &\\
& +4\frac{\partial H}{\partial y}(r,0)\!\left( \dot\pi\!\matc{0}{r}{0}{0}
\dot\pi\! \matc{-2}{\,0}{0}{\,1}
+ \dot\pi\!\matc{r}{0}{0}{0}\dot\pi\!\matc{0}{3}{0}{0}\!\right).
\end{align*} By
\eqref{partialx} and \eqref{partialy} we have
$$\frac{\partial H}{\partial x}(r,0)=\frac 12 \frac{d H}{dr}\, ,\qquad
\quad \frac{\partial H}{\partial y}(r,0)= \frac 1r
\bigl( \dot\pi(X_{-\al})H(r)- H(r)\dot\pi(X_{-\al})\bigr).$$ Now the
proposition follows easily.
\end{proof}

\medskip

The Theorems \ref{mainth1} and \ref{mainth2} are given in terms of linear
transformations. Now we will give the corresponding statements in terms of
matrices by choosing an appropriate basis.

If $\pi=\pi_{n,\ell}$ it is well known (see \cite{Hu}, p.32) that there
exists a basis $\vz{v_i}_{i=0}^\ell$ of $V_\pi$ such
that \begin{align*} &\dot\pi(H_\al)v_i=(\ell-2i)v_i, &\displaybreak[0]\\
&\dot\pi(X_\al)v_i=(\ell-i+1)v_{i-1}, \quad &
(v_{-1}=0), \displaybreak[0]\\ &\dot\pi(X_{-\al})v_i=(i+1)v_{i+1}, \quad &
(v_{\ell+1}=0).\\ \end{align*} Since we are dealing
with a representation of $\U(2)$ these relations have to be supplemented
with
$$\dot\pi(Z)v_i=(2n+\ell)v_i.$$
This follows from
\begin{align*}
\dot\pi(Z)&=\left( \frac{d}{dt}
\pi{\textstyle\matc{e^t}{0}{0}{e^t}}\right)_{t=0}=\left( \frac{d}{dt}
e^{2nt}\pi_\ell{\textstyle\matc{e^t}{0}{0}{e^t}}
\right)_{t=0} = (2n+\ell)I. \end{align*}

\medskip
\noindent We introduce the functions $h_i(r)$ by means of the relations
$$\tilde H(r)v_i=h_i(r)v_i.$$

\begin{cor}\label{sistema}
The function $\tilde H(r)=(h_0(r), \cdots, h_\ell(r) )$, $(r>0)$, satisfies
$(\tilde D\tilde H)(r)=\ld \tilde H(r)$ if and
only if \begin{align*} &(1+ r^2)^2
h_i''+\textstyle\frac{(1+r^2)}r\left(3+r^2-2r^2(n+\ell-i)
\right)h_i'-4i(\ell-i+1)(h_{i-1}-h_i)\displaybreak[0]\\ &
+\textstyle\frac{4(1+r^2)}{r^2}\bigl( (i+1)(\ell-i)(h_{i+1}-h_i)
+i(\ell-i+1)(h_{i-1}-h_{i}) \bigr) =\ld h_i
\end{align*}
for all $i=0,\cdots, \ell$.
\end{cor}
\begin{proof}
We have $J=X_{\al}-X_{-\al}$, $T=X_{\al}+X_{-\al}$ and $H_\g=\frac 12
Z+\frac 12 H_\al.$
Then
\begin{align*}
\dot\pi(J)^2\tilde H(r)v_i=&\;(\ell-i+1)(\ell-i+2)h_i
v_{i-2}\displaybreak[0]\\
-&\left((\ell-i+1)i+(i+1)(\ell-i)\right) h_i v_i +(i+1)(i+2)h_i v_{i+2},
\displaybreak[0]\\
\tilde H(r)\dot\pi(J)^2 v_i=&\;(\ell-i+1)(\ell-i+2)h_{i-2}
v_{i-2}\displaybreak[0]\\
-\,((\ell-i&+1)i+(i+1)(\ell-i)) h_i
v_i+(i+1)(i+2)h_{i+2}v_{i+2},\displaybreak[0]\\
\dot\pi(J)\tilde H(r)\dot\pi(J) v_i=&\;(\ell-i+1)(\ell-i+2)h_{i-1}
v_{i-2}+(i+1)(i+2)h_{i+1} v_{i+2}\displaybreak[0]\\
-&\left((\ell-i+1)ih_{i-1}+(i+1)(\ell-i)h_{i+1}\right)v_i,\displaybreak[0]\\
\dot\pi(T)^2\tilde H(r)v_i=&\;(\ell-i+1)(\ell-i+2)h_i
v_{i-2}\displaybreak[0]\\
+&\left((\ell-i+1)i+(i+1)(\ell-i)\right) h_i v_i +(i+1)(i+2)h_i v_{i+2},
\displaybreak[0]\\
\tilde H(r)\dot\pi(T)^2 v_i=&\;(\ell-i+1)(\ell-i+2)h_{i-2}
v_{i-2}\displaybreak[0]\\
+\,((\ell-i&+1)i+(i+1)(\ell-i)) h_i v_i+(i+1)(i+2)h_{i+2}
v_{i+2},\displaybreak[0]\\
\dot\pi(T)\tilde H(r)\dot\pi(T) v_i=&\;(\ell-i+1)(\ell-i+2)h_{i-1}
v_{i-2}+(i+1)(i+2)h_{i+1} v_{i+2}\displaybreak[0]\\
+&\left((\ell-i+1)ih_{i-1}+(i+1)(\ell-i)h_{i+1}\right)v_i,
\end{align*}

\smallskip
\noindent $\dot\pi(X_{-\al})\tilde H(r)\dot\pi(X_\al)v_i= i (\ell-i+1)
h_{i-1}v_i$,

\smallskip
\noindent $\tilde H(r)\dot\pi(X_{-\al})\dot\pi(X_\al)v_i= i (\ell-i+1)h_i
v_i$.

\smallskip
\noindent Therefore,
\begin{align*}
\bigl(\tilde D\tilde H\bigr)(r)v_i&= \Bigl(
(1+r^2)^2h_i''(r)+\textstyle\frac{(1+r^2)}r\left(3+r^2-2r^2(n+\ell-i)
\right)h_i'(r)
\displaybreak[0]\\ & \quad+\textstyle\frac{4(1+r^2)}{r^2}\bigl(
(i+1)(\ell-i)h_{i+1}(r)
+i(\ell-i+1)h_{i-1}(r)\displaybreak[0]\\ &
\quad- \left( i(\ell-i+1)+(i+1)(\ell-i)\right) h_{i}(r)
\bigr)\displaybreak[0]\\ &\quad -
4i(\ell-i+1)(h_{i-1}(r)-h_i(r))\Bigr) v_i.
\end{align*}
The proof is finished.
\end{proof}

\medskip
\begin{cor}\label{sistema2}
The function $\tilde H(r)=(h_0(r), \cdots, h_\ell(r) )$, $(r>0)$, satisfies
$(\tilde E\tilde H)(r)=\mu \tilde H(r)$ if and
only if
\begin{align*} (&n-\ell+3i)(1+ r^2)^2 h_i''+
6(i+1)(\ell-i)\,\textstyle\frac{(1+r^2)^2}{r}\,h'_{i+1}\displaybreak[0]\\
&-6i(\ell-i+1)\,\textstyle\frac{(1+r^2)}{r}\,
h'_{i-1}+\textstyle(n-\ell+3i)\frac{(1+r^2)}r\left(3+ r^2-2r^2(n+\ell-i)
\right)h_i'\displaybreak[0]\\ &
+4(n+2\ell-3i)\textstyle\frac{(1+r^2)}{r^2}\bigl(
(i+1)(\ell-i)(h_{i+1}-h_i) +i(\ell-i+1)(h_{i-1}-h_i)
\bigr)\displaybreak[0]\\ & + 4i(\ell-i+1)(2n+\ell+3)(h_{i-1}-h_i)
=\mu h_i
\end{align*}
for all $i=0,\cdots, \ell$.
\end{cor}
\begin{proof}
We have $\tilde H_2=\frac 12 Z-\frac 32 H_\al$, $\tilde H_1=\frac 12
Z+\frac 32 H_\al$, $\tilde
H_\g=\frac 12 Z+\frac 12 H_\al$. By using the computations in the proof of
Corollary \ref{sistema}
we obtain:
\begin{align*}
&\bigl(\tilde E\tilde H\bigr)(r)v_i= \Bigl((n-\ell+3i)(1+ r^2)^2 h_i'' \\
& \qquad + (n-\ell+3i) {\textstyle \frac{(1+r^2)^2}{r}} h_i'+ 2
(n+2\ell-3i) {\textstyle
\frac{(1+r^2)}{r}} h_i'\\ &\qquad +2 r(1+r^2) h_i'\bigl( 3(i+1)(\ell
-i)-(n+\ell-i)(n-\ell+3i)\bigr)\\ &\qquad
+6\textstyle\frac{(1+r^2)^2}{r} (i+1)(\ell-i)\,(h'_{i+1}-h_i')
-6\textstyle\frac{(1+r^2)}{r}
\,i(\ell-i+1)( h'_{i-1}-h_i')\\ & \qquad
+4(n+2\ell-3i)\textstyle\frac{(1+r^2)}{r^2}(i+1)(\ell-i)(h_{i+1}-h_i) \\
&\qquad
+4(n+2\ell-3i)\textstyle\frac{(1+r^2)}{r^2}i(\ell-i+1)(h_{i-1}-h_i)
\displaybreak[0]\\
& \qquad +
4i(\ell-i+1)(2n+\ell+3)(h_{i-1}-h_i) \Bigr)v_i.
\end{align*}
\end{proof}

\begin{cor}\label{scalar}
If $\tilde H(r)=(h_0(r), \cdots, h_\ell(r) )$, $(r>0)$, satisfies $(\tilde
D\tilde H)(r)=\ld \tilde H(r)$ and extends to a
$C^\infty$ function on $r\ge0$, then $h_0(0)=\cdots = h_\ell(0)$. \end{cor}
\begin{proof}
{}From Corollary \ref{sistema}, by multiplying by $r^2$ and
evaluating at $r=0$, we get
$$(i+1)(\ell-i)(h_{i+1}(0)-h_i(0))+i(\ell-i+1)(h_{i-1}(0)-h_{i}(0))=0$$
for $i=0,1,\dots,\ell$. Now the assertion follows by
induction starting at $i=0$ and ending at $i=\ell-1$.
\end{proof}

\medskip
We take the opportunity to sketch a proof of Theorem \ref{esfneg} based on
Corollaries \ref{sistema} and \ref{sistema2}.
\begin{esfneg}
For $n\leq 0$
$$\Psi_{n,\ell}(g)=\sum_{0\leq i \leq \min\vz{-n,\ell}} \binom{-n}{i}
a(g)^{-n-i}A(g)^{\ell-i}\cdot B(g)^i$$
is an irreducible spherical function of type $(n,\ell)$.
\end{esfneg}
\begin{proof}
Let $\pi=\pi_{n,\ell}$ and let
$$g(r)=\left(\begin{matrix} \frac{1}{(1+r^2)^{1/2}}& 0&
\frac{r}{(1+r^2)^{1/2}}\\ 0&1&0\\ \frac{-r}{(1+r^2)^{1/2}}& 0&
\frac{1}{(1+r^2)^{1/2}}\end{matrix}\right).$$ By Propositions \ref{defeq}
and \ref{lambdamu} it is enough to check that the
function $\tilde H(r)=\Psi_{n,\ell}(g(r))\Phi_\pi(g(r))^{-1}$ is an
eigenfunction of the differential operators $\tilde D$  and
$\tilde E$ given respectively in Theorems \ref{mainth1} and \ref{mainth2}.

Let $\{e_1,e_2\}$ be the canonical basis of $\CC^2$ and let
$v_i=e_1^{\ell-i}e_2^i$ be the corresponding basis of the space
of symmetric tensors of rank $\ell$. Then we can look at the coordinate
functions $(h_0,\dots,h_\ell)$ of $\tilde H$
associated to the basis $\{v_i\}$. Then it is not difficult to check that
$$h_i(r)=(1+r^2)^n \sum_{0\leq j\leq \min\vz{-n,\ell-i}} (-1)^j
\binom{-n}{j} \binom{\ell-i}{j} r^{2j}.$$

Now it is straightforward but lengthy to verify that $(h_0,\dots,h_\ell)$
is a simultaneous solution of the systems given in
Corollaries \ref{sistema} and \ref{sistema2} with $\lambda=4n(\ell+2)$ and
$\mu=4n(\ell+2)(n-\ell)$.
\end{proof}

\section{Extension to $G$}\label{extension}

Let us recall where we are: if $\Phi$ is a spherical function on $G$ of
type $\pi=\pi_{n,\ell}$ we have associated to  it a
function $\tilde H(r)$, $r\ge0$ with values in $\End(V_\pi)$ which is an
eigenfunction of the differential operators $\tilde D$
and $\tilde E$ (Theorems \ref{mainth1} and \ref{mainth2}). Now we want to
characterize the behavior of $\tilde H$ when
$r\rightarrow \infty$.

The complement ${\mathcal A}^c$ of ${\mathcal A}$ in $G$ is defined by the
condition $g_{33}=0$. This set is clearly left and
right invariant under $K$.

For any $t\in \RR$ let
$$a(t)= \left(\begin{matrix} \cos t& 0& \sin t\\ 0&1&0\\ -\sin t& 0&\cos
t\end{matrix}\right).$$

Then for all $ -\pi/2 <t< \pi/ 2$, $a(t)\in \mathcal A$ and $$a(\pi/2)=\lim
_{t\rightarrow \pi /2} a(t)
= \left( \begin{smallmatrix} 0& 0& 1\\
0&1&0\\ -1& 0&0\end{smallmatrix}\right)\in {\mathcal A}^c.$$

Now if $g\in {\mathcal A}^c$, $p(g)=(g_{13},g_{23}, 0)\in L$. Since the
line at infinity $L$ is a $K$-orbit,  there exists
$k\in K$ such that $p(g)=k p(a(\pi/2))= p(k a(\pi/2))$. Therefore
$g=ka(\pi/2)k'$ for some $k'\in K$. Thus ${\mathcal
A}^c=Ka(\pi/2)K$.

On ${\mathcal A}$ we have written $\Phi=H\Phi_\pi$. Hence if we put
$A(t)=A(a(t))$, for $-\pi/2 <t< \pi/ 2$ we have
$$\Phi(a(t))=\left(\cos t \right)^n H(a(t)) A(t)^\ell =\left(\cos t
\right)^n \tilde H(\tan t) A(t)^\ell,$$
since $p(a(t))=(\tan t, 0,1)$.

If we make the change of variable $r=\tan t$ and let $t\rightarrow \pi/2$
we obtain
$$\Phi(a(\pi/2))=\lim_{r\rightarrow \infty} (1+r^2)^{-\frac n2}\tilde H(r)
\matc{(1+r^2)^{-\frac 12 }}{\;0}{0}{\;1}^\ell,$$
where the exponent $\ell$ denotes the $\ell$-th symmetric power of the
matrix. If we use the basis $\vz{v_i}_{i=0}^\ell$ of
$V_\pi$ introduced before we have proved the existence of
$$\lim_{r\rightarrow \infty} (1+r^2)^{-\frac
{(n+\ell-i)}{2}}h_i(r)=L_i$$ for $i=0.\cdots ,\ell$. This follows directly
from $\dot \pi(H_\g)v_i=(n+\ell-i)v_i$ and $\exp
tH_\g=\matc{e^t}{0}{0}{1}$.

Given $x\in \CC$ with $|x|=1$ let us consider the element
$$b(x)=\left(\begin{matrix} x& 0& 0\\
0&1&0\\
0& 0&x^{-1}\end{matrix}\right)\in K.$$
Then $b(x)a(\pi/2)=a(\pi/2)b(x^{-1})$.
Now
\begin{align*}
\Phi\left(b(x)a(\textstyle\frac{\pi}2)\right)v_i&=\pi(b(x))\Phi(a(\textstyle
 \frac{\pi}2))v_i=x^{n+\ell-i}L_iv_i,
\displaybreak[0]\\ \Phi(a(\textstyle\frac{\pi}2)b(x^{-1}))v_i&=
\Phi(a(\textstyle\frac{\pi}2))\pi(b(x^{-1}))v_i=
L_ix^{-(n+\ell-i)}v_i. \end{align*} Therefore if $i\neq n+\ell$, then
$L_i=0$. Thus we have proved the direct part of the
following proposition.

\begin{prop}\label{extG}
If $\Phi$ is a spherical function on $G$ of type $\pi=\pi_{n,\ell}$ and
$\tilde H=(h_0, \cdots,
h_\ell)$ is the associated function, then \begin{enumerate} \item[i)]
$\tilde H=\tilde H(r)$ is a $C^\infty$-function for
$0\leq r<\infty$ and $\tilde H(0)= (1, \cdots, 1)$.

\item[ii)] If $i\neq n+\ell$, $\displaystyle\lim_{r\rightarrow
\infty}(1+r^2)^{-(n+\ell-i)/2}h_i(r)=0$.

\noindent If $i=n+\ell$,
$\displaystyle\lim_{r\rightarrow \infty} h_i(r)=L_i$ exists. \end{enumerate}

Conversely if $\tilde H$ is an eigenfunction of $\tilde D$ which satisfies
conditions i) and ii), then there exists a unique
$C^\infty$-eigenfunction $\Phi$ on $G$ of $\Delta_2$ which satisfies
condition ii) of Proposition \ref{defeq} to which $\tilde
H$ is associated. Moreover if $\tilde H$ is also an eigenfunction of
$\tilde E$ then $\Phi$ is spherical of type $(n,\ell)$.
\end{prop}
This proposition will play a crucial role in Section \ref{solutions}. We
will already come back to
it in Section \ref{aeigenfunc}.

\begin{proof}
As we said before the second part is the only one which needs to be proved.
First we want to extend $\tilde H$ to
the whole affine plane $\CC^2$, in such a way that the extended function
$H$ satisfies $H(kq)=\pi(k)H(q)\pi(k^{-1})$  for all
$k\in K$ and all $ q\in \CC^2$.  We observe that for $q=(x,y)\in
\CC^2-\vz{(0,0)}$ the element
$$k(q)=k(x,y)=\left(|x|^2+|y|^2\right)^{-{\textstyle\frac
12}}\matc{x}{-\overline y}{y}{\overline x}$$ lies in $K$ and
$q=k(q)(r,0)$ where $r=\left(|x|^2+|y|^2\right)^{\frac 12}$. Let
$$H(q)=\pi(k(q))\tilde H(r)\pi(k(q)^{-1}).$$ Then
$H$ is a
$C^\infty$ function on $\CC^2-\vz{(0,0)}$. Now we shall see that $H$
extends to a continuous function on $\CC^2$. In fact we
can equip $V_\pi$ with a scalar product such that $\pi(k)$ becomes a
unitary operator for all $k\in K$. If we denote with
$\|\,\|$  the corresponding operator norm on $\End(V_\pi)$ the we have $$\|
H(q)-I\|=\| \pi(k(q))(\tilde H(r)-I)
\pi(k(q)^{-1})\|= \|\tilde H(r)-I\|,$$
and our assertion follows.

Let us check now that $H(uq)=\pi(u)H(q)\pi(u^{-1})$ for all $u\in K, q\in
\CC^2$.
We may assume that $q=(x,y)\neq (0,0)$. Then
$k(uq)(r,0)=uq= uk(q)(r,0)$. Therefore there exists $m\in M$ such that
$k(uq)=uk(q)m$. Thus \begin{align*}
H(uq)&=\pi(k(uq))\tilde H(r)\pi(k(uq)^{-1})\\ &=
\pi(u)\pi(k(q))\pi(m)\tilde H(r)\pi(m^{-1})\pi(k(q)^{-1})\pi(u^{-1})\\ &=
\pi(u)H(q)\pi(u^{-1}),
\end{align*}
since $\pi(m)$ and $\tilde H(r)$ commute (Proposition \ref{Hdiagonal}).

Now we lift the function $H$ to a continuous function $H=H\circ p$ on
$\mathcal A$ with values in $\End(V_\pi)$ which is
$C^\infty$ on ${\mathcal A}-K$. Then the function $\Phi=H\Phi_\pi$ is a
continuous function on $\mathcal A$ which is a
$C^\infty$ eigenfunction of $\Delta_2$ on ${\mathcal A}-K$, satisfies
$\Phi(e)=I$ and $\Phi(k_1gk_2)=\pi(k_1)\Ph(g)\pi(k_2)$
for $k_1,k_2\in K$, $g\in \mathcal A$.

We shall next show that $\Phi$ can be extended to a continuous function on
$G$ with property ii) of Proposition \ref{defeq}.

As was pointed out before ${\mathcal A}^c=Ka(\pi/2)K$ and $a(\pi/2)=
\displaystyle\lim_{t\rightarrow \pi/2} a(t)$. So we put
$$\Phi(k_1a(\pi/2)k_2)=\pi(k_1)\displaystyle\lim_{t\rightarrow \pi/2} (\cos
t)^n\tilde H(\tan t)A(t)^\ell \pi(k_2).$$
First of all we see that by hypothesis $L=\lim_{t\rightarrow \pi/2} (\cos
t)^n\tilde H(\tan t)A(t)^\ell$ exists.
But we still
need to verify that if $k_1a(\pi/2)k_2=h_1a(\pi/2)h_2$ then
$\pi(k_1)L\pi(k_2)=\pi(h_1)L\pi(h_2)$, to have a good definition of
$\Phi$. This is equivalent to proving that if $ka(\pi/2)=a(\pi/2)h$ with
$k,h\in K$ then $\pi(k)L=L\pi(h)$. It is easy to see
that $ka(\pi/2)=a(\pi/2)h$ happens if and only if $k=b(x)m$ and
$h=b(x^{-1})m$, for some $x\in \CC, |x|=1$, and some $m\in M$.
Then $$\pi(k)L=\pi(b(x))\pi(m)L=\pi(b(x))L\pi(m) =
L\pi(b(x^{-1}))\pi(m)=L\pi(h),$$ because the hypothesis ii) of the
proposition says precisely that $\pi(b(x))L=L\pi(b(x^{-1}))$.

So we have a continuous function $\Phi$ on $G$ with values in $\End(V_\pi)$
which satisfies condition ii) of Proposition
\ref{defeq} and which is a $C^\infty$-eigenfunction of $\Delta_2$ on the
open set ${\mathcal A}-K$. Since the complement of
this set ${\mathcal A}^c\cup K$ is of Haar measure zero in $G$ and
$\Delta_2$ is an elliptic differential operator we conclude
(Weyl's theorem, see for instance \cite{Ho}, p.96) that $\Phi$ is a
$C^\infty$-eigenfunction of $\Delta_2$. The proposition is proved.
\end{proof}

\section{The inner product}\label{innerprod}

Given a finite dimensional irreducible representation $\pi=\pi_{n,\ell}$ of
$K$ in the vector space $V_\pi$ let
$(C(G)\otimes \End (V_\pi))^{K\times K}$ be the space of all continuous
functions $\Phi:G\longrightarrow \End(V_\pi)$  such
that $\Phi(k_1gk_2)=\pi(k_1)\Phi(g)\pi(k_2)$ for all $g\in G$, $k_1,k_2\in
K$. Let us equip $V_\pi$ with an inner product such
that $\pi(k)$ becomes unitary for all $k\in K$. Then we introduce an inner
product in the vector space $(C(G)\otimes \End
(V_\pi))^{K\times K}$ by defininig
\begin{equation*}
\langle \Phi,\Psi \rangle =\int_G \tr ( \Phi(g)\Psi(g)^*)\, dg\, ,
\end{equation*}
where $dg$ denote the Haar measure on $G$ normalized by $\int_G dg=1$, and
where $\Psi(g)^*$ denotes the adjoint of $\Psi(g)$
with respect to the inner product in $V_\pi$.

Let us write $\Phi=H\Phi_\pi$, $\Psi=F\Phi_\pi$ on the open set $\mathcal
A$ of $G$ and put
$\tilde H(r)=(h_1(r),\cdots, h_\ell(r))$, $\tilde F(r)=(f_1(r),\cdots,
f_\ell(r))$ as we did in Section \ref{secc1var}.

\begin{prop}\label{prodint}
If $\Phi, \Psi\in \left(C(G)\otimes \End (V_\pi)\right)^{K\times K}$ then
\begin{equation*}
\langle \Phi,\Psi \rangle =4\sum_{i=0}^\ell \int_0^\infty
r^3(1+r^2)^{-(n+\ell+3-i)}h_i(r)\overline{f_i(r)}\, dr.
\end{equation*}
\end{prop}
\begin{proof}
Let us consider the element $H=E_{13}-E_{31}\in \lieg$. Then $H$ is
conjugate by an element of $G$ to $i(H_\al+H_\beta)$,
thus $\text{ad }H$ has $0$ and $\pm i$ as eigenvalues with multiplicity 2
and $\pm 2i$ as eigenvalues with multiplicity 1.

Let $A=\exp \RR H$ be the Lie subgroup of $G$ of all elements of the form
$$a(t)= \exp tH=\left(\begin{matrix} \cos t& 0&
\sin t\\ 0&1&0\\ -\sin t& 0&\cos t\end{matrix}\right)\, ,\qquad t\in \RR.$$

Then the function $D:A\longrightarrow \RR $ defined by $$D(a(t))=\prod_\nu
\sin it \nu(H), $$
where the product runs over the eigenvalues $\nu=i,2i$ of $\text{ad } H$
counted with multiplicities, becomes
$D(a(t))=-\sin^2 t\sin 2t$.

Now Corollary 5.12, p. 191 in \cite{Hee}
establishes that
$$\int_{G/K} f(gK)\,dg_K=2\int_{K/M}\Big(\int_0^{\pi/2}
|D(a(t))|f(ka(t)K)\,dt\Big)\,dk_M\,,$$
where $dg_K$ and $dk_M$ are repectively the invariant measures on $G/K$ and
 $K/M$
normalized by $\int_{G/K} dg_K=\int_{K/M} dk_M=1$.

Since the function $g\mapsto \tr(\Phi(g)\Psi(g)^*)$ is invariant under left
and right
multiplication by elements in $K$, we have
$$\langle \Phi,\Psi\rangle = 2 \int_{0}^{\pi/2} \sin ^2t\sin 2t
\,\tr\left( \Phi(a(t)\Psi(a(t))^*\right)\,dt.$$

If we put $r=\tan t$ for $0<t<\pi/2$ we have $$\tr\left(
\Phi(a(t)\Psi(a(t))^*\right)=\sum_{i=0}^\ell (1+r^2)^{-(n+\ell-i)}
h_i(r)\overline{f_i(r)}.$$ (See Section \ref{secc1var}). Then
\begin{align*}
2\int_{0}^{\pi/2} \sin ^2t\sin 2t&
\,\tr\left( \Phi(a(t)\Psi(a(t))^*\right)\,dt \\ =&
4\sum_{i=0}^\ell\int_0^\infty r^3
(1+r^2)^{-(n+\ell+3-i)}h_i(r)\overline{f_i(r)}\, dr.
\end{align*}
Thus the proposition follows. \end{proof}

\begin{prop}\label{Deltasim}
If $\Phi,\Psi\in\left(C^\infty(G)\otimes \End (V_\pi)\right)^{K\times K}$
then
$$\langle \Delta_j\Phi, \Psi\rangle=\langle \Phi,\Delta_j \Psi\rangle
\qquad j=1,2.$$
\end{prop}
\begin{proof}
If we apply a left invariant vector field $X\in \lieg$ to the function
$g\mapsto\tr(\Phi(g)\Psi(g)^*)$ on $G$ and then we
integrate over $G$ we obtain
$$0= \int_G \tr\left( (X\Phi)(g)\Psi(g)^*\right)\, dg+\int_G \tr\left(
\Phi(g)(X\Psi)(g)^*\right)\, dg.$$ Therefore $\langle
X\Phi,\Psi\rangle=-\langle \Phi,X\Psi\rangle$. Now let
$\tau:\lieg_\CC\longrightarrow \lieg_\CC$ be the conjugation of
$\lieg_\CC$ with respect to the real linear form $\lieg$. Then $-\tau$
extends to a unique antilinear involutive ${}^*$
operator on $D(G)$ such that $\left(D_1D_2\right)^*=D_2^*D_1^*$ for all
$D_1,D_2\in D(G)$. This follows easily from the fact
that the universal enveloping algebra over $\CC$ of $\lieg$ is canonically
isomorphic to $D(G)$. Then it follows that $\langle
D\Phi, \Psi\rangle=\langle \Phi,D^* \Psi\rangle$.

Now $H_\al^*=H_\al$, $H_\beta^*=H_\beta$ , $X_\al^*=X_{-\al}$,
$X_\beta^*=X_{-\beta}$, $X_\g^*=X_{-\g}$. {}From this it is easy
to verify that $\Delta_2^*=\Delta_2$ and $\Delta_3^*=\Delta_3$. This
completes the proof of the proposition.
\end{proof}

\smallskip
We now describe an alternative way to motivate the choice of inner product
made before Proposition \ref{prodint}.

Start from the system in Corollary \ref{sistema}. Introducing
$t=(1+r^2)^{-1}$, $0<t<1$, we get
\begin{equation}\label{D}
\begin{split}
t\,&(1-t) h_i'' -\left((n+\ell-i+3)t-(n+\ell-i+1)\right)
h_i'\displaybreak[0]\\
& +(\ell-i)(i+1)\,\frac{h_{i+1}-h_{i}}{1-t} - i(\ell-i+1) \,t\,
\frac{h_{i}-h_{i-1}}{1-t}={\textstyle\frac \ld 4}h_i
\end{split}
\end{equation}
with $h_i'=\frac{dh_i}{dt}$. If $D$ denotes the operator acting on
the vector $H(t)=(h_1(t),\cdots, h_\ell(t))$
this can be put in the form
$$DH(t)=D_1H(t)+D_2H(t)={\textstyle\frac \ld 4}H(t).$$
Here
$$(D_1H)=\frac 1{(1-t)t^{n+\ell-i}}\frac{d}{dt}\left( (1-t)^2
t^{n+\ell-i+1} \frac d{dt}\right) h_i$$
and
$$(D_2H)(t)=(\ell-i)(i+1)\,\frac{h_{i+1}-h_{i}}{1-t} -
i(\ell-i+1) \,t\, \frac{h_{i}-h_{i-1}}{1-t}.$$
The aim
is to make $D$ into a densely defined symmetric operator in some
appropriate $L^2$ space. A look at
$D_1$ suggests an obvious choice for an inner product among vector
valued functions $H(t)$ which are $C^2$  and
compactly supported inside the open set $(0,1)$, namely
\begin{equation}
\langle H, K\rangle=
\sum_{i=0}^\ell \int_0^1 h_i(t)\overline k_i(t) (1-t)t^{n+\ell-i}\, dt.
\end{equation}
If one were trying to make $D_1$ symmetric, then each term in this
sum could be multiplied by some arbitrary constant
$\mu_i>0$. One check that the requirement $$\langle DH,K \rangle=
\langle H,D K\rangle$$
forces one to make the choice made above. A change of variables
shows now that, up to a multiplicative constant,
this is the inner product introduced before Proposition \ref{prodint}.

\section{ The $C^\infty$ eigenfunctions of $D$}\label{autofunciones}
This section contains material that is crucial in determining explicit
expressions for all the spherical functions associated
to the complex projective plane. We display these results in Section
\ref{solutions} for some small values of $\ell$.

In the process of obtaining these expressions we have uncovered what we
think are interesting properties of eigenfunctions
for the differential operator $D$. In order to keep the length of this
paper within reasonable bounds we face a problem in
deciding what to include here, as well as the level of generality aimed for
in the proof of some of these statements. We have
opted for a rather mixed approach whereby several lines of attack are
described while proofs are given only in some simple
cases. We are confident that all the statements made below hold true, even
if some of the proofs that we have at this time are
lengthy and rather unilluminating. The statements made below are strong
enough that they can be used in Section
\ref{solutions} to prove that we have exhibited all the spherical functions
of the pair $(G,K)$.

Some of these different ways of getting at the solutions might be useful in
potential applications of this material.

In subsection \ref{l=0} we describe all solutions of $DH = \lambda H$ in
the case $\ell = 0$ and $n$
an arbitrary integer. This is just a review of the classical theory dealing
with the indicial
equation and related elementary material. The analysis is powerful enough
to show that we have written down all spherical
functions when $\ell = 0$. In this case, and only in this case, we give a
complete analysis including the non-analytic
solutions of $DH = \lambda H$ for all integer values of $n$. We also give,
for $\ell=0$, the ``limit point, the limit circle"
analysis of H. Weyl that would be needed to discuss several possible
selfadjoint extensions of the symmetric operator $D$.

Then we describe a way of constructing a $2(\ell+1)$ dimensional space of
solutions of $DH = \lambda H$ when $\ell$ is a
non-negative integer. Here $\lambda$ is an arbitrary parameter.

This construction is based on reducing a second order system of $\ell + 1$
equations to a scalar equation of order
$2(\ell+1)$ and in observing, as the first in a string of miracles, that
the corresponding differential operator has a nice
factorization featuring the Gauss hypergeometric differential operator.
This elimination process (from the bottom up) produces
the last component of a (pair of) vector of solutions while the rest of the
vector can then be trivially constructed by
applying explicit differential operators to this last component. We then
see how the entries in these column vectors can be
used to give the first entries in a set of $2\ell$ other vector solutions
of the same system $DH =
\lambda H$. Once again the rest of these vectors are built by applying to
their first component differential operators that
can be read off from the equation. This construction has some of the flavor
of an unrelated ``miracle'' uncovered by
T.~Koornwinder in \cite{Ko}. See his comments at the end of Section ~1.

The remarkable factorizations in Subsection 8.2 have not yet yielded to a
a nice and general proof for all values of $\ell$. They have been
established
by use of computer algebra for $\ell$ in the range $0,1,\dots,10$.

 An interesting feature of this construction is that it expresses the
solution in terms of differential operators acting on
the Gauss function $\lw{2}F_1$. Alternatively one can express the solution
as a linear combination of ``shifted'' versions of
Gauss' function.

We note that the elimination process mentioned above would fail if one were
trying to apply it to determine solutions of $EH
= \mu H$.

We start by recalling that the classical (scalar) equation of Gauss
\begin{equation}\label{eq1}
\left(t(1-t) \left( \frac {d}{dt} \right)^2 + (c-(a+b+1)t) \frac {d}{dt} -
ab\right) h(t) = 0
\end{equation}
has a two-dimensional space of local solutions $\Phi(a,b,c;t)$ which for
generic values of $a,b,c$ has a basis given by
$$\varphi_1(t) = \lw{2}F_1\left(
\begin{smallmatrix} a,\;b \\ c \end{smallmatrix};t \right) \quad \text{ and
} \quad
\varphi_2(t) = t^{1-c} \lw{2}F_1
\left( \begin{smallmatrix} a-c+1,\;b-c+1 \\ 2-c
\end{smallmatrix};t\right).$$
Any equation of the form
\[ \left( t(1-t)  \left(\frac {d}{dt}\right)^2 + (c-(p+1)t) \frac {d}{dt} -
\mu\right) h(t) = 0 \]
can be solved easily in terms of these solutions. It is enough to introduce
a parameter $w$ as one of the solutions of
\[ \mu = -w(w+p) \] and to put $a = -w$, $b = w+p$. Notice that a different
choice of a solution of the quadratic equation in
$w$ amounts to exchanging $a$ and $b$ above and that they enter
symmetrically in \eqref{eq1}.

\subsection{All solutions of $4DH = \lambda H$ for
$\ell=0$.}\label{l=0}
With the strategy just outlined and the  form of $D$
given in (\ref{D}), we have with $\lambda = -4w(w+n+2)$,
\[ \left( t(1-t) \left( \frac {d}{dt} \right)^2 + ((n+1)-(n+3)t) \frac
{d}{dt} + w(w+n+2)\right) h_0(t) = 0. \]
For generic values of the parameters the functions
$$
\varphi_w(t) = \lw{2}F_1 \left( \begin{smallmatrix} -w,w+n+2 \\ n+1
\end{smallmatrix};t\right) \quad \text{and}\quad \psi_w(t)
= t^{-n} {}_2F_1 \left( \begin{smallmatrix} -w-n,w+2 \\ 1-n
\end{smallmatrix};t\right)
$$
are linearly independent solutions. When $n+1$ or $1-n$ is a non-positive
integer these expressions are of no use, and a more
careful analysis is done below.

Since the replacement of $w$ by $-(w+n+2)$ does not change neither the
equation nor the solutions $\varphi_w,\psi_w$,  we will
consider these values as equivalent.

We will show that for each integer $w$ satisfying $w \ge 0$ and $w+n\ge 0$
we have one solution of our equation that
satisfies the conditions of Proposition \ref{extG}, namely
\begin{equation}\label{eq2}
\lim_{t \rightarrow 1} h_0(t) = \mbox{finite and non-zero} \end{equation}
and
\begin{equation}\label{eq3}
\begin{split}
&\mbox{if $n = 0$ then $\lim_{t \rightarrow 0} h_0(t)$
exists,}\displaybreak[0]\\
&\mbox{if $n\ne 0$ then $\lim_{t \rightarrow 0}  h_0(t)t^{n/2} = 0$.}
\end{split}
\end{equation}

\subsubsection{ Case $n = 0$.} The functions $\varphi_w(t)$ and $\psi_w(t)$
coincide
and give us one solution. This is a polynomial for $w = 0,1,2,\dots$. For
other values of $w$ it behaves, close  to $t = 1$,
like a non-zero constant times $(1-t)^{-1}$.

For $w= -1$, $\varphi_w(t)$ is well defined and we get that \[ \varphi_w(t)
= \frac {1}{1-t} \quad \mbox{ and }
\quad (\log t)\varphi_w(t) \]
are linearly independent solutions. Neither one of them nor any linear
combination of them is bounded at both $t = 0$  and $t
= 1$.

For $w \ne -1$, there is a unique solution of the form \[ (\log
t)\varphi_w(t) + \sum_{0\le i} a_it^i. \]
This is linearly independent of $\varphi_w(t)$ and blows up at $t = 0$.
Therefore when $n = 0$ we only have eigenfunctions
meeting the asymptotic conditions (\ref{eq2}) and (\ref{eq3}) for $w =
0,1,2,\dots$ (or the redundant set $w = -2,-3,\dots$).

\subsubsection{ Case $n > 0$.} Now $\varphi_w(t)$ is well defined for all
$w$, and it
is a polynomial if $w = 0,1,2,\dots$. For other values of $w$ (except the
``redundant'' values $w = -n-2,-n-3,\dots)$  it
behaves like a non-zero constant times $(1-t)^{-1}$ around $t = 1$. The
function $\psi_w(t)$ is well defined and linearly
independent of $\varphi_w(t)$ if $w =-n-1,-n,\dots,-2,-1$ and it is then a
polynomial in $t$ times $t^{-n}$. The product
$t^{n/2}\psi_w(t)$ behaves like $t^{-n/2}$, and blows up, around $t = 0$.

If $w \ne -n-1,-n,\dots,-2,-1$, then there is a unique solution of the form
\[ (\log t)\varphi_w(t) + \sum_{-n\le i,i \ne 0}
a_it^i. \] This solution is linearly independent of $\varphi_w(t)$ and its
product with $t^{n/2}$ is unbounded at $t = 0$.
Thus the only values of $w$ where there is an eigenfunction meeting the
conditions \eqref{eq2} and \eqref{eq3} are $w =
0,1,2,\dots$.

\subsubsection{Case $n < 0$} Now $\psi_w(t)$ is always well defined, and it
is a polynomial in $t$ if
$w = -n,-n+1,-n+2,\dots$ (and also if $w= -2,-3,-4,\dots$). For other
values of $w$ it behaves like a non-zero constant  times
$(1-t)^{-1}$ around $t = 1$. The function $\varphi_w(t)$ is well defined if
$w = -1,0,1,\dots,-n-1$, and it is in fact a
polynomial in $t$ which takes the value $1$ at $t = 0$. It gives a solution
that is linearly independent of $\psi_w(t)$, and
the product $t^{n/2}\varphi_w(t)$ is unbounded at $t = 0$.

If $w \ne -1,0,1,\dots,-n-1$, then there is a unique solution of the form
\[ (\log t)\psi_w(t)
+ \sum_{0 \le i,i \ne -n} a_it^i. \] This solution is linearly independent
of $\psi_w(t)$ and its product with $t^{n/2}$  is
unbounded at $t = 0$, since $a_0 \ne 0$.

By going over the results given above we see that for each integer $w$
satisfying $w \ge 0$ and $w+n \ge 0$ we have {\em one}
solution of our equation that satisfies the properties given in Proposition
$6.1$. For values of $w$ that do not satisfy these
conditions (once one takes into account the symmetry between $-w$ and
$w+n+2$ mentioned earlier) there is {\em no} solution of
the equation that satisfies the conditions in Proposition $6.1$.

We stress that these conditions single out those solutions of our equation
that correspond to spherical functions on $G$.
They should not be construed as ``boundary conditions'' in the traditional
sense of specifying a self-adjoint extension of our
densely defined symmetric operator.

In this regard we close the analysis of the case $l=0$ with the observation
that for $n=0$ the end point $t = 0$ gives a  {\em
limit circle} case (in the terminology introduced by H.~Weyl) while $t = 1$
gives a {\em limit point} case. For $n$ non-zero
both end points are of the {\em limit point} type. These statements can
easily be proved from the detailed description of the
solutions given above.

\subsection{Solutions of $4DH =\lambda H$, $\ell > 0$.} Consider
now the
system \eqref{D} given in Section~\ref{innerprod}.
We fix the values of $(\ell,n)$ and $\lambda$. Exploiting the tridiagonal
nature of our system of equations we can eliminate
--from the bottom up-- all components of $H(t)$ and obtain one equation in
$h_\ell(t)$ of order $2\ell+2$, \[ P_{2\ell+2}
\left( \frac{d}{dt} \right) h_\ell(t) = 0 \] where the following
factorization holds
\[ P_{2\ell+2} \left( \frac {d}{dt} \right) = A_{2\ell} \left( \frac
{d}{dt} \right) \left( t(1-t) \left(\frac {d}{dt}
\right)^2 + (n+1 - (n+\ell+3)t) \frac {d}{dt} - {\textstyle\frac\ld
4}\right) \] with
$A_{2\ell}\left( \frac{d}{dt} \right)$ a differential
operator of order $2\ell$. For $\lambda = -4w(w+n+\ell+2)$ and from the
discussion in the introduction of the section, we can
take for $h_\ell(t)$ any local solution in the space
\[ \Phi(-w,w+n+\ell+2,n+1;t). \]
components $h_i(t)$, $i = 0,1,\dots,\ell-1$ of a vector solution $H(t)$ by
applying to $h_\ell(t)$ differential operators  of
orders $2(\ell-i)$, \[ h_i(t) = R_{2(\ell-i)}\left( \frac {d}{dt} \right)
h_\ell(t). \] These operators can be read off from
the system $4DH = \lambda H$. In particular $R_0 \left( \frac{d}{dt}
\right)= I$.

We claim that from this vector solution (actually any vector in this
two-dimensional space) we can readily obtain an extra
set of $\ell$ vector solutions (actually a $2\ell$ dimensional space of
such).

\noindent Proceed as follows: Pick $k = 1,\dots,\ell$ and consider the
system \[ 4DH(t) = (\lambda - 4k(n+k+1))H(t). \]
Eliminate (from the top down) all components of $H(t)$ except $h_0(t)$. At
the last step one gets an equation of order
$2\ell+2$ for $h_0(t)$, \[ M_{2\ell+2} \left( \frac {d}{dt} \right) h_0 = 0
\] where another remarkable factorization holds,
namely, \begin{align*} &M_{2\ell+2} \left( \frac {d}{dt} \right) {\tilde
R}_{2(\ell-k)} \left( \frac {d}{dt} \right) = \\
&B_{4\ell-2k} \left( \frac {d}{dt} \right) \left( (1-t) t \left( \frac
{d}{dt} \right)^2 + ((n+k+1) - (n+k+\ell+3)t) \frac
{d}{dt} - \lambda \right).
\end{align*}
Here $B_{4\ell-2k}\left( \frac {d}{dt} \right)$ is a differential operator
of order $4\ell-2k$ and ${\tilde R}_{2(\ell-k)}
\left( \frac {d}{dt} \right)$ is the operator obtained from the
$R_{2(\ell-k)} \left( \frac {d}{dt} \right)$ introduced above
by replacing $n$ by $n+k$. In consequence we can take for $h_0(t)$ anything
in the space \[ {\tilde R}_{2(\ell-k)} \left(
\frac {d}{dt} \right) \Phi(-w,w+n+k+\ell+2,n+k+1;t) \] where $w$ is now any
root of
\[ \lambda = -4w(w+n+k+\ell+2). \]
Just as before, from $h_0(t)$ one gets the remaining components of a vector
solution $H(t)$ by the recipe
\[ h_k(t) = Q_{2k} \left( \frac {d}{dt} \right) h_0(t). \] In conclusion we
have $H(t)$ satisfying
$DH = 4(-w(w+n+k+\ell+2) - k(n+k+1))H$, and the relation between $w$ and
the original $\lambda$ is given by
\[ \lambda = -4w(w+n+k+\ell+2) - 4k(n+k+1). \] Note that the first
vector $H(t)$, constructed from $h_\ell(t)$,
corresponds to $k = 0$. In a very explicit sense this $h_\ell(t)$ is the
``seed'' for the entire construction.

\smallskip
The discussion above is summarized as follows: For a fixed $(n,\ell)$ we
associate to each integer $k$ satisfying  $0 \le k
\le \ell$ a two-dimensional space of solutions of $4DH = \lambda H$. The
solutions can be constructed fairly explicitly as the
result of applying certain differential operators to solutions of Gauss's
hypergeometric equation. The construction involves
an auxiliary parameter $w$ related to $\lambda$ by
\[ \lambda = -4w(w+n+k+\ell+2) - 4k(n+k+1). \] The choice of one of the
two roots of this quadratic equation  does not
affect the final result, i.e.~the construction is invariant under the
exchange of $-w$ into $w+n+k+\ell+2$.

The parameters $k$ and $w$ will reappear in Section \ref{solutions}. In
that section
we will use the results of Section \ref{parametriz} giving a
parametrization of those
solutions of $4DH =\lambda H$ which correspond to spherical functions. The
parameter $\ld$ will be replaced by $\ld/4$.  For
the convenience of the reader we anticipate here that these conditions will
turn out to be $0 \le k \le \ell$, as above, and
the extra requirement that $w$ or $-(w+n+k+\ell+2)$ can be chosen to
satisfy the conditions $$ 0 \le w \, ,\qquad \quad 0 \le
w+n+k.$$ An instance of these conditions was seen in Subsection \ref{l=0}
in the special case of $\ell=0$.

By using the well known differentiation formula \[ \frac {d}{dt}\,
\lw{2}F_1 \left( \begin{smallmatrix} a\,,\,b \\ c
\end{smallmatrix};t\right) = \textstyle\frac {ab}{c} \,\lw{2}F_1 \left(
\begin{smallmatrix} a+1\,,\,b+1 \\ c+1
\end{smallmatrix};t\right), \] one can re-express the components of
different vector solutions $H(t)$ of $4DH= \lambda H$ in
terms of ``shifted'' version of $\lw{2}F_1$, as seen below.

We give a fairly general formula for analytic solutions, at $t=0$, when $n
\ge 0$ and
then illustrate it in a couple of cases.

For fixed $\ell$, and $n \ge 0$, and any integer $0 \le k \le \ell$ we have
a vector solution $H^k(t)$ with components
\[H_j^k(t) = \sum_{p=-\min(k,j)}^{\ell-\max(k,j)} z(j,p)\, \lw{2}F_1 \left(
\begin{smallmatrix}\ -w+p,w+n+\ell+2+k - \min(j,k) \\
n+1+\ell-j \end{smallmatrix}; t\right). \] Here $\lambda =
-4w(w+n+\ell+2+k) - 4k(n+k+1)$ and $z(j,p)$ are appropriate
coefficients.

We include two examples. If $\ell=1$ and $k=1$ we have, up to a common
scalar multiple and for $n\ge0$
\begin{align*}
H_0^1(t)&= \lw{2}F_1 \left( \begin{smallmatrix}\ -w\,,\,w+n+4 \\ n+2
\end{smallmatrix};t\right), \displaybreak[0] \\  H_1^1(t)
&=\textstyle \frac {-(n+1)}{w+n+3} \biggl( (w+n+2)\, \lw{2}F_1 \left(
\begin{smallmatrix} -w-1\,,\,w+n+3 \\ n+1
\end{smallmatrix};t\right) + \lw{2}F_1 \left( \begin{smallmatrix} -w,w+n+3
\\ n+1 \end{smallmatrix};t\right)\biggr).
\end{align*}
If $\ell=2$ and $k=2$ we have, up to a scalar multiple and for $n\ge0$
\begin{align*}
H_0^1(t) &=\lw{2}F_1 \left( \begin{smallmatrix} -w\,,\,w+n+6 \\ n+3
\end{smallmatrix};t\right), \\
H_1^1(t) &=\textstyle \frac {-(n+2)}{w+n+5} \biggl ( -(w+n+3)\, \lw{2}F_1
\left( \begin{smallmatrix} -w-1\,,\,w+n+5 \\
n+2 \end{smallmatrix};t\right) \\
&\qquad + \lw{2}F_1 \left(
\begin{smallmatrix} -w\,,\,w+n+5 \\
n+2 \end{smallmatrix};t\right)\biggr ) ,\\ H_2^1(t) &=\textstyle \frac
{(n+1)(n+2)}{(w+n+5)} \biggl ( \frac {w+n+3}{2}\, \lw{2}F_1 \left(
\begin{smallmatrix} -w-2\,,\,w+n+4 \\ n+1 \end{smallmatrix};
t\right) \\ &\qquad +\textstyle \frac {(w+n+3)}{w+n+4}\,\lw{2}F_1 \left(
\begin{smallmatrix} -w-1\,,\,w+n+4 \\ n+1
\end{smallmatrix};t\right) + \frac {1}{w+n+4}\, \lw{2}F_1 \left(
\begin{smallmatrix} -w\,,\,w+n+4 \\ n+1
\end{smallmatrix};t\right)\biggr ). \end{align*}

In Section \ref{solutions} we will give more compact expressions for these
solutions in terms of generalized hypergeometric
functions. The relation between those expresions and the ones here can also
be seen by using \cite{FI} and \cite{FW}.

\section{Parametrizations of spherical functions}\label{parametriz}

As we pointed out in Section \ref{Spherical functions}, there exists a one
to one correspondence between the set
of all equivalence classes of finite dimensional irreducible
representations of $G$ which contain the representation
$\pi_{n,\ell}$ of $K$ and the set of all equivalence classes of irreducible
spherical functions of $(G,K)$ of type
$(n,\ell)$.

We need to quote a restriction theorem of finite dimensional irreducible
representations of
$\GL(n,\CC)$ to
$\GL(n-1,\CC)$ which we shall use. The equivalence classes of finite
dimensional irreducible
holomorphic
representations of $\GL(n,\CC)$ are parametrized by the $n$-tuples of
integers
$m_1\geq m_2\cdots \geq m_n$ in the following way: the diagonal subgroup of
$\GL(n,\CC)$ acts on
the highest weight vector, with respect to the Borel subgroup of all upper
triangular matrices, of
the representation $\tau_{(m_1,\cdots, m_n)}$ by $t_1E_{11}+\cdots
+t_nE_{nn}\longmapsto
t_1^{m_1}\cdots t_n^{m_n}$.

If $\GL(n-1,\CC)$ is identified with the subgroup of all $(n-1)\times
(n-1)$ matrices of
$\GL(n,\CC)$ in the
following way $$\GL(n-1,\CC)\simeq\left(\begin{smallmatrix}
\GL(n-1,\CC)&0\\0&1
\end{smallmatrix} \right),$$ then we have (see \cite{Ze}, p.186):

\begin{prop} When we restrict the representation $m_1\geq m_2\cdots \geq
m_n$ of $\GL(n,\CC)$
to $\GL(n-1,\CC)$ it decomposes as the direct sum of the representations
$k_1\geq \cdots \geq
k_{n-1}$ of $\GL(n-1,\CC)$ such that $m_1\geq k_1\geq m_2\geq k_2\cdots
\geq k_{n-1}\geq m_n$, all
of these with multiplicity one. \end{prop}

Any finite dimensional irreducible representation of $G=\SU(3)$ extends to
a unique
holomorphic
irreducible representation of $\SL(3,\CC)$, and all of these are
restrictions of holomorphic
irreducible representations of $\GL(3,\CC)$. Moreover the highest weight of
the
restriction of
$\tau_{(m_1,m_2, m_3)}$ to $\SL(3,\CC)$ is
$\ld=(m_1-m_2)\ld_1+(m_2-m_3)\ld_2$ where
$\ld_1,\ld_2$
are the fundamental weights of $\SL(3,\CC)$ given by $\ld_1(H)=x_1$ and
$\ld_2(H)=-x_3$, for
$H=x_1E_{11}+x_2E_{22} +x_3E_{33}\in \liesl(3,\CC)$.

Now $K=\mathrm{S}(\U(2)\times\U(1))\simeq \U(2)$ and any finite dimensional
irreducible
representation of $\U(2)$ extends to a unique holomorphic irreducible
representation of
$\GL(2,\CC)$. Thus a $\U(2)$-irreducible subrepresentation of
$\tau_{(m_1,m_2, m_3)}$
corresponds
to a $K$-irreducible subrepresentation parametrized by a pair $(n,\ell)\in
\ZZ\times
\ZZ_{\geq 0}$
which we want to determine.

The parameter $\ell$ can be read off the action of the diagonal subgroup of
$\SU(2)\times \vz{1}$
on a heighest weight vector of the subrepresentation $(k_1,k_2)$ of
$\U(2)$. Thus
$t^{k_1}\left(t^{-1} \right)^{k_2}=t^{k_1-k_2}$ gives $\ell=k_1-k_2$.

The parameter $n$ can instead be read off the action of the center of $K$
in the $(k_1,k_2)$ subrepresentation of $\U(2)$. From
$$\left(\begin{smallmatrix} t&0&0\\0&t&0\\0&0&t^{-2}\end{smallmatrix}
\right)= \left(\begin{smallmatrix}
t^3\,&0&0\\0&t^3\,&0\\0&0&1\end{smallmatrix} \right)
\left(\begin{smallmatrix}
t^{-2}&0&0\\0&t^{-2}&0\\0&0&t^{-2}\end{smallmatrix} \right)$$ we get
$t^{2n}t^\ell=t^{3(k_1+k_2)}t^{-2(m_1+m_2+m_3)}$. Therefore
$2n+k_1-k_2=3(k_1+k_2)-2(m_1+m_2+m_3)$ which gives $n=k_1+2k_2-m_1-m_2-m_3$.

Now from the previous discussion we easily get

\begin{prop}\label{param}
When we restrict to $K$ the finite dimensional irreducible representation
of $G$ with highest weight $\ld=p\ld_1+q\ld_2$  it
decomposes as the direct sum of the representations $\pi_{n,\ell}$ of $K$
given by $n=k_1+2k_2-p-2q$, $\ell=k_1-k_2$ with
$p+q\geq k_1\geq q\geq k_2\geq 0$, all of these with multiplicity one.
\end{prop}
\begin{proof}
The irreducible representation of $G$ with highest weight
$\ld=p\ld_1+q\ld_2$ can be realized as the restriction to $G$  of
the representation $\tau_{(m_1,m_2,m_3)}$ of $\GL(3,\CC)$ with $m_1=p+q,
m_2=q, m_3=0$. Thus $\tau_{(m_1,m_2,m_3)}$ restricted
to $\U(2)$ is the direct sum of the representations $(k_1,k_2)$ with
$p+q\geq k_1\geq q\geq k_2\geq 0$, and a fortiori
$\tau_{(m_1,m_2,m_3)}$ restricted to $K$ is the direct sum of
$\pi_{n,\ell}$ with
$n=k_1+2k_2-p-2q$, $\ell=k_1-k_2$. The proposition follows. \end{proof}

\begin{cor}\label{param2}
The equivalence classes of irreducible spherical functions of $(G,K)$ of
type $(n,\ell)$ can be pa\-ra\-me\-tri\-zed by  the
set of all tu\-ples $(p,q,k_1,k_2)\in \ZZ^4$ such that $p+q\geq k_1 \geq
q\geq k_2\geq 0$, with $n=k_1+2k_2-p-2q$ and
$\ell=k_1-k_2$. Moreover if $\Phi_{(p,q,k_1,k_2)}$ is a spherical function
in the class $(p,q,k_1,k_2)$ then
\begin{align*}
\Delta_2\Phi_{(p,q,k_1,k_2)}&= -\textstyle\frac 43\left( p^2+q^2+pq +3p+3q
\right)\Phi_{(p,q,k_1,k_2)}\displaybreak[0]\\
\Delta_3\Phi_{(p,q,k_1,k_2)}&=4\textstyle( \frac 29 p^3-\frac 29
q^3+\frac13p^2q -\frac13 pq^2+
2p^2+pq+4p+2q)\Phi_{(p,q,k_1,k_2)}.
\end{align*}
\end{cor}
\begin{proof}
The only thing that remains to be proved are the last two formulas. Since
$\Delta_2$ and
$\Delta_3$ are in $D(G)^G$ they act as scalars in any finite dimensional
irreducible representation
of $G$. These scalars $\ld(\Delta_2)$ and $\ld(\Delta_3)$ can be computed
by looking at the action
of $\Delta_2$ and $\Delta_3$ on a heighest weight vector. {}From Proposition
\ref{defDelta} it
follows that
\begin{align*}
\ld(\Delta_2)=&-\ld(H_\al)^2-\textstyle\frac 13
\ld(Z)^2-2\ld(H_\al)-2\ld(Z),\displaybreak[0]\\
\ld(\Delta_3)=&\textstyle\frac 89\ld(H_\al)^3-\frac 89\ld(H_\beta)^3+ \frac
43
\ld(H_\al)^2\ld(H_\beta)-\frac 43
\ld(H_\al)\ld(H_\beta)^2\\
&+8\ld(H_\al)^2+4\ld(H_\al)\ld(H_\beta)+16\ld(H_\al)+8\ld(H_\beta).
\end{align*}

Replacing $\ld(H_\al)=p$, $\ld(H_\beta)=q$ and $\ld(Z)=p+2q$ in the above
expressions the proof of
the corollary is finished. \end{proof}

\medskip
Going back to Proposition \ref{lambdamu} we give in the following corollary
the eigenvalues $\lambda$ and $\mu$
corresponding respectively to the differential operators $D$ and $E$ in
terms of the parameters $(n,l,p,q)$.

\begin{cor}\label{expldmu}
We have
\begin{align*}
\lambda=&-\textstyle\frac 43(p^2+q^2+pq+3p+3q)+\frac 43(\ell^2+n^2+\ell
n+3\ell+3n),\displaybreak[0]\\
\mu=&4(\textstyle\frac 29 p^3-\frac 29 q^3+\frac 13p^2q-\frac 13
pq^2+p^2-q^2+p-q)\\
&-4(\textstyle\frac 29 \ell^3-\frac 29 n^3+\frac 13\ell^2n-\frac 13 \ell
n^2+\ell^2-n^2+\ell-n).
\end{align*}
\end{cor}

\medskip
For a given $(n,\ell)$ the four parameters $(p,q,k_1,k_2)$ are subject to
two relations. One can introduce two integer  free
parameters $w$, $k$ by means of the expressions
\begin{align*}
k_1=&\ \ w+n+\ell+k,\\
k_2=&\ \  w+n+k,\\
p=&\ \ w+\ell-k,\\
q=&\ \ w+n+2k,
\end{align*}
which are consistent with the relations in Proposition \ref{param} and it
is easy to see that the four inequalities in  this
proposition are equivalent to the following four conditions: $0\le w$,
$0\le k\le\ell$ and $0\le w+n+k$.
We recognize here the conditions on $w,k$ mentioned in Section
\ref{autofunciones}.

In terms of this new parameters we have
\begin{equation}\label{lambda}
\begin{split}
\lambda=&-4w(w+n+\ell+k+2)-4k(n+k+1),\displaybreak[0]\\ \mu=&\lambda
(n+3k-\ell)-12k(\ell+1-k)(n+k+1).
\end{split}
\end{equation}

\section{The Analytic Eigenfunctions of D and E}\label{aeigenfunc}

We take up the time honored method of formal power series expansions to
look for simultaneous analytic eigenfunctions of the operators $D$ and $E$.
The system
\eqref{D} given in Section \ref{innerprod}  was
obtained from the one given in Corollary \ref{sistema} by making the change
of variables $t=(1+r^2)^{-1}$ and we got
\begin{equation}\label{Dagain}
\begin{split}
t\,&(1-t) h_i'' +\left((n+\ell-i+1)-t (n+\ell-i+3)\right)
h_i'\displaybreak[0]\\
& +(\ell-i)(i+1)\,\frac{h_{i+1}-h_{i}}{1-t} - i(\ell-i+1) \,t\,
\frac{h_{i}-h_{i-1}}{1-t}=\textstyle\frac \ld 4h_i
\end{split}
\end{equation}
Similarly from
Corollary \ref{sistema2} we get the system \begin{equation}\label{E}
\begin{split}
&(n-\ell+3i)t\,(1-t) h_i'' -3(i+1)(\ell-i)h'_{i+1}\displaybreak[0]\\
&+(n-\ell+3i)\left(n+\ell-i+1-t(n+\ell-i+3)\right)
h_i'+3i(\ell-i+1)th'_{i-1} \displaybreak[0]\\ &+ (n+2\ell-3i)(\ell-i)(i+1)\,
\frac{h_{i+1}-h_{i}}{1-t}-(n+2\ell-3i) i(\ell-i+1)\,
\frac{h_{i}-h_{i-1}}{1-t} \displaybreak[0]\\ &-(2n+\ell+3)i(\ell-i+1)
(h_i-h_{i-1})={\textstyle\frac {\mu}{4}}h_i.
\end{split}
\end{equation}
If we change $\lambda/4$ and $\mu/4$
respectively by $\lambda$ and $\mu$ and use matrix notation, both systems
are equivalent to $DH=\lambda H$ and $EH=\mu H$
where the differential operators $D$ and $E$ are defined below.
\begin{align}
DH&=t(1-t)H''+(A_0-tA_1)H'+\frac{1}{1-t}(B_0-tB_1)H.
\label{DD}\displaybreak[0]\\
EH&=t(1-t)MH''+(C_0-tC_1)H'+\frac{1}{1-t}(D_0+tD_1)H. \label{EE}
\end{align}
Here $H$ denotes the column vector
$H=(h_0,\dots,h_\ell)$ and the coefficient matrices are given by:
\begin{align*}
A_0&=\textstyle\sum_{i=0}^\ell(n+\ell-i+1)E_{ii},\displaybreak[0]\\
A_1&=\textstyle\sum_{i=0}^\ell(n+\ell-i+3)E_{ii},\displaybreak[0]\\
B_0&=\textstyle\sum_{i=0}^\ell(i+1)(\ell-i)(E_{i,i+1}-E_{i,i}),\displaybreak
 [0]\\ B_1&=\textstyle\sum_{i=0}^\ell
i(\ell-i+1)(E_{ii}-E_{i,i-1}),\displaybreak[0]\\
M&=\textstyle\sum_{i=0}^\ell(n-\ell+3i)E_{ii},\displaybreak[0]\\
C_0&=\textstyle\sum_{i=0}^\ell(n-\ell+3i)(n+\ell-i+1)E_{ii}-3\sum_{i=0}^\ell
 (i+1)( \ell-i)E_{i,i+1},\displaybreak[0]\\
C_1&=\textstyle\sum_{i=0}^\ell(n-\ell+3i)(n+\ell-i+3)E_{ii}-3\sum_{i=0}^\ell
 i(\ell-i+1)E_{i,i-1},
\displaybreak[0]\\
D_0&=\textstyle\sum_{i=0}^\ell(n+2\ell-3i)(i+1)(\ell-i)(E_{i,i+1}-E_{ii})
\displaybreak[0]\\
 &\quad -3 \textstyle\sum_{i=0}^\ell(n+\ell-i+1)
i(\ell-i+1)(E_{ii}-E_{i,i-1}),\displaybreak[0]\\
D_1&=\textstyle\sum_{i=0}^\ell(2n+\ell+3)i(\ell-i+1)(E_{ii}-E_{i,i-1}).
\end{align*}

\begin{prop}\label{Hanali}
An irreducible spherical function $\Phi$ on $G$ of type $(n,\ell)$
corresponds to a function $H$ that is analytic at $t=0$
and satisfies $DH=\lambda H$ and
$EH=\mu H$ for some $\lambda,\mu\in\CC$. Moreover if $n+\ell<0$ the
function $H(t)$ has a zero of order at least $-n-\ell$ at
$t=0$.
\end{prop}
\begin{proof} In the previous sections we already established that $H$ is
an eigenfunction of $D$ and $E$.

that
$$\tilde H(r)=(1+r^2)^{n/2}\Phi\left(\begin{matrix} \frac 1{(1+r^2)^{1/2}}&
0& \frac r{(1+r^2)^{1/2}}\\ 0&1&0\\
\frac {-r}{(1+r^2)^{1/2}}& 0& \frac
1{(1+r^2)^{1/2}}\end{matrix}\right)\left(\begin{matrix} (1+r^2)^{1/2}& 0\\
0& 1\end{matrix}\right)^{\ell}.$$
The first thing we get from here is that $\tilde H(r)=\tilde H(-r)$. In
fact we have
\begin{align*}
&\left(\begin{matrix} \frac 1{(1+r^2)^{1/2}}& 0& \frac
{-r}{(1+r^2)^{1/2}}\\ 0&1&0\\ \frac r{(1+r^2)^{1/2}}& 0& \frac
1{(1+r^2)^{1/2}}\end{matrix}\right)\displaybreak[0] \\ &\quad =
\left(\begin{matrix} 1& 0& 0\\ 0&-1&0\\ 0& 0&
-1\end{matrix}\right)\left(\begin{matrix}\frac 1{(1+r^2)^{1/2}}& 0& \frac
r{(1+r^2)^{1/2}}\\ 0&1&0\\ \frac
{-r}{(1+r^2)^{1/2}}& 0& \frac 1{(1+r^2)^{1/2}}\end{matrix}\right)
\left(\begin{matrix} 1& 0& 0\\ 0&-1&0\\ 0&0&-1\end{matrix}\right).
\end{align*}
Moreover
$$\pi\left(\begin{matrix} 1& 0& 0\\ 0&-1&0\\ 0&0&-1\end{matrix}\right)
\qquad \text{and }\qquad
\Phi\left(\begin{matrix}\frac 1{(1+r^2)^{1/2}}& 0& \frac r{(1+r^2)^{1/2}}\\
0&1&0\\ \frac {-r}{(1+r^2)^{1/2}}& 0&  \frac
1{(1+r^2)^{1/2}}\end{matrix}\right)$$ commute because they are diagonal
matrices in the basis $\{v_0,\dots,v_{\ell}\}$.
Hence $\tilde H(r)=\tilde H(-r)$.

If we make the change of variable $s=1/r$ then the function $H=H(s)$ is
even and for $s>0$ we have
$$\tilde H(s)=\frac{(1+s^2)^{n/2}}{s^n}\Phi\left(\begin{matrix} \frac
s{(1+s^2)^{1/2}}& 0& \frac 1{(1+s^2)^{1/2}}\\ 0&1&0\\
\frac {-1}{(1+s^2)^{1/2}}& 0& \frac
s{(1+s^2)^{1/2}}\end{matrix}\right)\left(\begin{matrix} \frac
{(1+s^2)^{1/2}}s& 0\\
0& 1\end{matrix}\right)^{\ell}.$$
Thus $h_i(s)$ is even and $s^{n+\ell-i}h_i(s)$ is analytic at $s=0$. Now in
terms of the variable  $t=(1+r^2)^{-1}$ we have
$s^2=t/(1-t)$. Therefore $h_i=h_i(t)$ has a pole at $t=0$ of order at most
$\left[\frac{n+\ell-i}2\right]$ when $n+\ell-i \geq
2$. When $n+\ell-i \leq 1$, $h_i=h_i(t)$ is analytic at $t=0$ and has also
a zero or order at least
$-\left[\frac{n+\ell-i}2\right]$, when $n+\ell-i\leq -1$.

If we put $H(t)=\sum_{j\ge m}\,t^jH_j$, and assume that $H_m\ne 0$, with
column vector
coefficients
$H_j=(H_{0,j},\dots,H_{\ell,j})\in\CC^{\ell+1}$, then we have
\begin{equation}\label{cond1}
H_{i,j}=0 \qquad \text{ for all } j<-\left[\frac {n+\ell-i}2\right].
\end{equation}
On the other hand from $DH=\lambda H$ and (\ref{DD}) we get that the
following  three term recursion relation holds for all
$j$ \begin{equation}\label{recur1}
\begin{split}
&\big[(j-1)(j-2)+(j-1)A_1-B_1+\lambda\big]H_{j-1}\\
&-\big[2j(j-1)+j(A_0+A_1)- B_0+\lambda\big]H_j+(j+1)(j+A_0)H_{j+1}=0.
\end{split}
\end{equation}
\noindent ~From here for $j=m-1$ we obtain
\begin{equation}\label{cond2}
m(m+n+\ell-i)H_{i,m}=0.
\end{equation}
Now assume that $m<0$. {}From \eqref{cond2} we get $H_{i,m}=0$ for all
$i\neq
m+n+\ell$.
For $i=m+n+\ell$ we have
$$\left[\frac{n+\ell-i}2\right]=\left[\frac{-m}2\right]<-m.$$ Then from
\eqref{cond1}  we obtain $H_{i,m}=0$. Thus $H_m=0$,
which is a contradiction. Therefore $m\geq 0$ and $H=H(t)$ is analytic at
$t=0$.

To prove the last assertion suppose that $n+\ell<0$. First observe that
$m\geq -\left[\frac{n+\ell}2\right]>0$.  Now assume
that $m\leq -n-\ell-1$. Then $m+n+\ell -i\leq -1-i<0$ for all $0\leq i\leq
\ell$. Hence \eqref{cond2} implies that $H_m=0$,
which completes the proof of the proposition.
\end{proof}

\medskip
Thus we want to consider now formal power series
$H(t)=\sum_{j=0}^\infty\,t^jH_j$
with column vector coefficients $H_j\in\CC^{\ell+1}$.

A formal solution of the system $DH=\ld H$ and $EH=\mu H$ is then given by
solving the pair of three term recursion  relations
given by \eqref{recur1} and
\begin{equation}\label{recur2}
\begin{split}
&\big[(j-1)(j-2)M+(j-1)C_1+D_1+\mu\big]H_{j-1}+(j+1)(jM+C_0)H_{j+1}\\
&-\big[2j(j-1)M+j(C_0+C_1)-D_0+\mu\big]H_j=0.
\end{split}
\end{equation}

 Let $V(\lambda)$ denote the vector space of all formal power series
$H(t)=\sum_{j=0}^\infty\,t^jH_j$ such that \eqref{recur1}
holds. The differential operators D and E, closely related to $\Delta_2$
and $\Delta_3$, commute. Therefore E restricts to a
linear map of $V(\lambda)$ into itself.

Now we shall first consider the case $n\geq 0$ where we have complete
results and
then we shall make some comments in the case $n<0$.

Let $\eta:V(\lambda)\rightarrow \CC^{\ell+1}$ be the linear map defined by
$\eta:H(t)\mapsto H(0)$ and let $L$ be the
$(\ell+1)\times(\ell+1)$ matrix given by
\begin{equation}\label{L}
L=L(\ld)=D_0-C_0A_0^{-1}(B_0-\lambda).
\end{equation}
Notice that $n\ge0$ insures that
$A_0$ is invertible.

\begin{thm}\label{amainth} If $n\ge0$ then $\eta:V(\lambda)\rightarrow
\CC^{\ell+1}$
is a surjective isomorphism. Moreover, the following is a commutative
diagram $$\begin{CD}
V(\ld) @ >E >>V(\ld) \\ @ V \eta VV @ VV \eta V \\ \CC^{\ell+1} @> L >>
\CC^{\ell+1}
\end{CD}$$
\end{thm}
\begin{proof} The matrix $j+A_0$ is invertible for all $j\ge0$ since
$$\det(j+A_0)=\Pi_{i=0}^\ell(j+n+\ell-i+1).$$ Therefore
the recurrence relation \eqref{recur1} shows that the coefficient vector
$H_0=H(0)$ determines $H(t)$. This proves the first
assertion.

To prove the second statement let $H(t)\in V(\lambda)$. Then from
(\ref{EE}) and (\ref{recur1}) we get
\begin{align*}
\eta(EH(t))&=EH(0)=C_0H'(0)+D_0H(0)=C_0H_1+D_0H_0\\
&=(-C_0A_0^{-1}(B_0-\lambda)+D_0)H_0=LH_0=L(\eta(H(t)).
\end{align*}
\end{proof}

In the next theorem we give the eigenvalues of the matrix $L=L(\lambda)$
and their multiplicities, as well as  their geometric
multiplicities. We have only partial results on the eigenvectors of $L$ and
we do not display them here. In Section
\ref{solutions} we write them down in the cases $\ell=0,1,2$.

\begin{thm}\label{Jordan} For any $(n,\ell)\in\ZZ\times\ZZ_{\ge0}$ the
characteristic polynomial of the matrix  $L(\lambda)$
is given by
$$\det(\mu-L(\lambda))=\prod_{k=0}^\ell(\mu-\mu_k(\lambda)),$$ where
$\mu_k(\lambda)=\lambda(n-\ell+3k)-3k(\ell-k+1)(n+k+1)$.
Moreover, all eigenvalues $\mu_k(\lambda)$ of $L(\lambda)$ have geometric
multiplicity one. In other words all eigenspaces are
one dimensional.
\end{thm}
\begin{proof}
Let us consider the polynomial $p\in\CC[\lambda,\mu]$ defined by
$p(\lambda,\mu)=\det(\mu-L(\lambda))$. For each integer $k$
such that $0\le k\le\ell$ let $\lambda(w,k)=-w(w+n+\ell+k+2)-k(n+k+1)$ and
let
$\mu_k(\lambda)=\lambda(n+3k-\ell)-3k(\ell-k+1)(n+k+1)$.

\noindent Then from Theorem \ref{amainth} and \eqref{lambda} we have
$$p(\lambda(w,k),\mu_k(\lambda(w,k))=0,$$ for all $w\in\ZZ$
such that $0\le w$ and $0\le w+n+k$. Since there are infinitely many such
$w$, the polynomial function
$w\mapsto p(\lambda(w,k),\mu_k(\lambda(w,k)))$ is identically zero on
$\CC$. Hence, for given $k$ $(0\le k\le\ell)$, we have
$p(\lambda,\mu_k(\lambda))=0$ for all $\lambda\in\CC$.

Now $\Lambda=\{\lambda\in\CC:\mu_k(\lambda)=\mu_{k'}(\lambda) \text{ for
some } 0\le k<k'\le\ell\}$
is a finite set, in fact $\mid\Lambda\mid\le\ell(\ell+1)/2$. Since for any
$\lambda$,
$p(\lambda,\mu)$ is a monic polynomial in $\mu$ of degree $\ell+1$, it
follows that
if $\lambda\in\CC-\lambda$ then
\begin{equation}\label{p}
p(\lambda,\mu)=\prod_{k=0}^\ell(\mu-\mu_k(\lambda))
\end{equation}
for all $\mu\in\CC$.

Now it is clear that \eqref{p} holds for all $\lambda$ and all $\mu$,
which completes the proof of the first asertion.

To prove the second statement of the theorem we point out that the matrix
$L=L(\lambda)$ is a four diagonal matrix of the
form $$L=\sum_{i=1}^{\ell}a_iE_{i,i-1}+\sum_{i=0}^{\ell}b_iE_{i,i}
+\sum_{i=0}^{\ell- 1}c_iE_{i,i+1}+
\sum_{i=0}^{\ell-2}d_iE_{i,i+2},$$ with
\begin{align*}
a_i&=3i(\ell-i+1)(n+\ell-i+1),\\
b_i&=\ld(n-\ell+3i)-3(i+1)(\ell-i)(\ell-2i) - 3i(\ell-i+1)(n+\ell-i+1),\\
c_i&= 3(i+1)(\ell -i) \left (\ell-2i
-(n+\ell-i)^{-1}((i+2)(l-i+1) +\ld)\right),\\ d_i&=3(i+1)(i+2) (\ell-i-1)
(\ell-i)(n+\ell-i)^{-1},
\end{align*}
for all $1\le i\le\ell$.

If $LH_\mu=\mu H_\mu$ with $H_\mu=(h_0,\dots,h_\ell)\in\CC^{\ell+1}$ then
for all $0\le
i\le\ell$ we have
$$a_ix_{i-1}+b_ix_i+c_ix_{i+1}+d_ix_{i+2}=\mu x_i,$$ where we must
interpret $a_0=c_\ell=d_\ell=d_{\ell-1}=0$.
{}From these equations we see that $H_\mu$ is determined by $x_\ell$, which
proves that the geometric multiplicity of the
eigenvalues of $L$ is one. \end{proof}

\begin{cor}\label{eigenv}
For $n\ge 0$ the formal solutions $H(t)=\sum_{j=0}^\infty\,t^jH_j$ of our
system $DH=\lambda H$, $EH=\mu H$ are
parametrized, up to a nonzero multiplicative constant, by the eigenvalues
$\mu$ of the matrix $L=L(\lambda)$.
More precisely, if $H_{\mu}\in \CC^{l+1}$ is a $\mu$-eigenvector of $L$
then $$\eta^{-1}(H_\mu)=\sum_{j=0}^\infty t^jH_j$$
is the corresponding solution.
\end{cor}
\noindent Notice that from the definition of $\eta$ it follows that
$H_0=H_\mu$.

\smallskip
This concludes the discussion of the case $n\geq 0$. In this case none of
the components of $H(t)$ vanishes at $t=0$ unless
$H(t)\equiv 0$.

When $n<0$ things are not so simple. The structure of the order of the
zeros of the components of $H(t)$ at $t=0$, leads us
to consider the following cases: $n\leq -\ell$ and
$-\ell <n<0$. A consequence of this structure is that the isomorphism
$\eta$ introduced in the case
$n\geq 0$ must be redefined, as we do below.

The partition of the integers in the form
\begin{equation}\label{partition}
\ZZ=\vzm{n}{n\geq 0}\cup \vzm{n}{-\ell<n< 0}\cup \vzm{n}{n\leq -\ell}
\end{equation}
was already alluded to in the comments following Corollary \ref{pi*}.

When $n<0$, if $H\in V(\lambda)$ then we have $H(t)=\sum_{j\ge
a}^\infty\,t^jH_j$
where $a=\max\vz{0,-n-\ell}$. See Proposition \ref{Hanali}.

If $n+\ell\leq 0$ from \eqref{recur1}, putting $j=a-1$, we obtain
$aiH_{i,a}=0$. Hence $H_{i,a}=0$ for $1\le i\le\ell$.
Similarly a closer look at \eqref{recur1}  reveals that $H_{i,a+k}=0$
for $0\le k\le\ell-1$ and $k+1\le i\le\ell$, and moreover that the map
$\eta:V(\lambda)\rightarrow \CC^{\ell+1}$ defined by
$$\eta(H)=(H_{0,a},H_{1,a+1},\dots,H_{\ell,a+\ell})$$ is a linear
surjective isomorphism.

If $0< n+\ell<\ell$ from \eqref{cond1} we get that $H_{i,0}=0$ for all
$n+\ell+1\le i\le\ell$. But this does not follow from
\eqref{recur1}. So in this case we should consider the vector space
$W(\lambda)=\{H\in V(\lambda):H_{i,0}=0 \, ,\text{\ for
all\ } n+\ell+1\le i\le\ell\}$. Then it is easy to verify that the elements
$H\in W(\lambda)$ satisfy
$$H_{i,k}=0 \text{\ for all $0\le k\le -n-1$ and $n+\ell+k+1\le
i\le\ell$},$$ which is stronger than \eqref{cond1}. In
particular $H_{i,1}=0$ for all $n+\ell+2\le i\le\ell$. {}From this and
\eqref{recur2} one can prove that the differential
operator $E$ restricts to a linear operator on $W(\lambda)$. Moreover from
\eqref{recur1} it follows that the map
$\eta:W(\lambda)\rightarrow \CC^{\ell+1}$ defined by
$$\eta(H)=(H_{0,0},\dots,H_{n+\ell,0},H_{n+\ell+1,1},H_{n+\ell+2,2},\dots,H_
 {\ell,-n})$$ is a linear surjective isomorphism.

Notice that the structure of the order of zeros discussed above, insures
that any analytic function $H(t)$ in $V(\lambda)$ or
$W(\lambda)$, depending on the value of $n$, meets condition ii) of
Proposition \ref{extG}.

Unfortunately, the matrix $L$ which completes the commutative diagram
corresponding to the one given in Theorem
\ref{amainth}, in the two cases going with $n<0$, is not as simple as the
one obtained when $n\ge0$.

\section{ The spherical functions}\label{solutions}
The purpose of this section is to combine results given in previous
sections to give a fairly explicit expression for the
spherical functions associated to the complex projective plane in terms of
a particular class of generalized hypergeometric
functions. We have the following results which have been verified, by
use of computer algebra, for values
of $\ell$ up to ten, but for which we do not have
yet a general proof.

\begin{conj}\label{conj1}
For a given $\ell \ge 0$, the spherical functions corresponding to the pair
$(\ell,n)$ have components that are expressed  in
terms of generalized hypergeometric functions of the form $
\lw{p+2}F_{p+1}$, namely \[
\lw{p+2}F_{p+1} \left( \begin{smallmatrix} a,b,s_1+1,\dots,s_p+1
\vspace{2pt} \\
c,s_1,s_2,\dots,s_p
\end{smallmatrix} ;t \right) = \sum_{j=0}^{\infty} \frac
{(a)_j(b)_j}{j!(c)_j} (1 + d_1j + \dots + d_pj^p)t^j. \]
\end{conj}

In the case $n \ge 0$ we can give a more precise form of our results, namely

\begin{conj}\label{conj2} For $\ell \ge 0$ and  $n\ge 0$ the spherical
functions are given by diagonal matrices corresponding
to $0 \le k \le \ell$. For a given $k$ the entries along the diagonal make
up a vector $H_i^k(t)$, $0\le i \le \ell$, with
\begin{align*}
& H_i^k(t) = H_i^k(0) \times \\
&\lw{\ell\!+\!2-|k\!-\!i|}F_{\ell\!+\!1-|k-i|} \left(\begin{array}{c}
\!-w\!-\min(i,k),w\!+\!n\!+\!\ell\!+\!k+2-\min(i,k),(s_j\!+\!1) \\
n+\ell-i+1,(s_j)
\end{array};t \right).
\end{align*}
The vector $H_i^k(0)$ is an eigenvector of $L(\lambda)$ defined in
\eqref{L}.  The eigenvalues $\lambda,\mu$ of $D$ and $E$ (see \eqref{DD}
and \eqref{EE}) are given by $\lambda = -w(w+n+\ell+k+2) - k(n+k+1)$ and
$\mu =
\lambda(n-\ell+3k) - 3k(\ell-k+1)$. Finally $(s_j)$ denotes a
row vector of $\ell-|i-k|$ denominator parameters and
$(s_j+1)$ the same coefficients shifted up by one.
\end{conj}

\rmk{Remark} By definition of the function $\lw{p+2}F_{p+1}$ on the
left-hand side of Conjecture \ref{conj1} we have
\[ \sum_{j=0}^{\infty} \frac {(a)_j(b)_j}{j!(c)_j} \frac {(s_1+1)_j\dots
(s_p+1)_j}{(s_1)_j\dots (s_p)_j} \,t^j. \]

Notice that the ratio of the $p$ factors involving $s_1,\dots,s_p$ is given
by \[
\frac {(s_1+j)(s_2+j)\dots (s_p+j)}{s_1s_2 \dots s_p} \] and this serves to
define the polynomial $1+d_1j + \dots + d_pj^p$
appearing in Conjecture \ref{conj1}.

The generalized hypergeometric functions $\lw{p+1}F_p$ have not, so far,
played in representation theory and/or harmonic
analysis a role comparable to that of the celebrated case of Gauss and
Euler, namely $\lw{2}F_1$. An exception to this are the
cases of $p=2,3$ , and a few other low values of $p$
as well documented for instance in \cite{AAR}. Another very
interesting exception to the statement above are the cases
of the families $\lw{3n-4}F_{3n-5}$ 
related to matrix entries for representations of
U(n) discussed in \cite{KG} and \cite{KV}.
For arbitrary $p$ the monodromy of the corresponding equation
has been analized in \cite{BH}.

We observe that the special type of generalized hypergeometric functions
that appear here do satisfy differential equations
of lower order than the order $p+2$ that one should expect. This has been
observed for instance in \cite{LVW}. We are grateful
to M. Ismail for this reference.

\medskip
Next we give a guide that should allow the reader to see how the results in
Sections \ref{parametriz} and \ref{aeigenfunc}, as well as Proposition
\ref{extG} can be
used to see that for $\ell=0,1,2,\dots$ we are exhibiting all spherical
functions
of type $(n,\ell)$. We need to consider separately the cases:
$$\text{a)}\quad n\ge0,\quad\quad \text{b)}\quad n\le -\ell,\quad\quad
\text{c)}\quad
-\ell<n<0,$$
as indicated in \eqref{partition}.

The results in Section \ref{aeigenfunc} are complete for the case $n\ge0$
and simpler
than in the other cases, and we start with it.

arbitrary
$\lambda\in\CC$ we should pick $\mu$ as one of the eigenvalues
$\mu_k(\lambda)$,
$0\le k\le\ell$, of the matrix $L(\lambda)$. Then we determine $H_\mu$ as
the properly normalized eigenvector. Starting with
$H_0=H_\mu$ the three term recursion relation (\ref{recur1}) reduces to a
two term recursion, and this leads to our
hypergeometric functions in terms of a parameter $w\in\CC$ related to
$\lambda$ by $\lambda=-w(w+n+\ell+k+2)-k(n+k+1)$. Since
these hypergeometric functions remain bounded at $t=1$ only when they are
polynomials, we are forced to take $w=0,1,2,\dots$.
Using the appropriate version of Gauss' summation formula we see that we
can normalize $H_\mu$ in order to get
$H(1)=(1,\dots,1)$. Recall that we have shown that the corresponding
diagonal matrix is a scalar matrix, \ref{scalar}. Now
Proposition \ref{extG} guarantees that
$H(t)$ in fact corresponds to a spherical function $\Phi$ on $G$ of type
$(n,\ell)$.
By refering to Section \ref{parametriz} one sees that we have a complete
list.

It is remarkable that for the cases b) and c) the same machinery can be
used, although we do not yet have the  corresponding
theoretical statements. The insistence on bounded eigenfunctions forces a
different choice of $w$, which again is consistent
with the conditions $0\le w$ and $0\le w+n+k$ in Section \ref{parametriz}.

\noindent {11.0. \bf The case $\ell=0$.} We see that the complete list of
spherical functions is given (up to a  scalar
multiple) by:

\noindent {\em a) Case } $n \ge 0$. For $w = 0,1,2,\dots$ we have $\lambda
= -w(w+n+2)$,
$\mu=n\ld$ and
\[
\varphi_w(t) = \lw{2}F_1 \left( \begin{smallmatrix} -w,\,\;w+n+2
\vspace{2pt} \\ n+1
\end{smallmatrix} ;t\right).
\]

\noindent {\em b) Case } $n < 0$. For $w = -n,-n+1,\dots$
we have $\lambda = -w(w+n+2)$, $\mu=n\ld$ and \[
\psi_w(t) = t^{-n}\, \lw{2}F_1 \left( \begin{smallmatrix} -w-n,\,\;w+2
\vspace{2pt} \\ 1-n \end{smallmatrix} ;t\right). \]
These results were already given, and proved, in
Section~\ref{autofunciones}. In this case the operator $E$  is proportional
to $D$.

\subsection{ The case $\ell=1$}
The complete list of matrix (or vector) valued spherical functions is given
below.

\smallskip
\noindent {\em a) Case }$n\ge 0$. We have two families of spherical
functions corresponding either to the
choice $\mu_0,\mu_1$ indicated in
Section \ref{aeigenfunc} or the choice $k = 0,1$ alluded to at the end of
Section \ref{autofunciones} and  given in Section
\ref{aeigenfunc}.

\noindent {\em a.0)} For $k=0$, and $w = 0,1,2,\dots$ we have $\lambda =
-w(w+n+3)$, $\mu=\ld (n-1)$ and
\[
H_0 = \left(\begin{array}{c}
1 -\frac{\lambda}{n+1} \vspace{2pt} \\ 1 \end{array} \right), \]
\[
H(t) = \left( \begin{array}{c}
\left(1-\frac{\lambda}{n+1} \right) \lw{3}F_2 \left( \begin{smallmatrix}
-w,\,\;w+n+3,\,\;\lambda-n \vspace{2pt} \\
n+2,\,\;\lambda - n-1
\end{smallmatrix} ;t \right) \vspace{2pt} \\ \lw{2}F_1
\left(\begin{smallmatrix}
-w,\,\;w+n+3 \vspace{2pt} \\
n+1
\end{smallmatrix} ;t \right)
\end{array} \right).
\]

\smallskip
\noindent {\em a.1)} For $k=1$, and $w = 0,1,2,\dots$ we have $\lambda =
-w(w+n+4)-n-2$,
$\mu=(\ld-3) (n+2)$ and
\[ H_0 = \left( \begin{array}{c}
1 \vspace{2pt} \\ -(n+1) \end{array} \right), \]
\[
H(t) = \left( \begin{array}{c}
\lw{2}F_1 \left( \begin{smallmatrix}
-w,\,\;w+n+4 \vspace{2pt} \\
n+2
\end{smallmatrix} ;t\right) \vspace{2pt} \\ -(n+1)\, \lw{3}F_2
\left(\begin{smallmatrix} -w-1,\,\; w+n+3,\,\;\lambda
\vspace{2pt} \\ n+1,\,\;\lambda-1
\end{smallmatrix} ;t\right)
\end{array} \right).
\]

\

\noindent {\em b) Case }$n < 0$.
Again we have two families corresponding to the choice $k=0$ or $1$ made in
the case $n\ge 0$.

\noindent {\em b.0)} For $k=0$, and $w = -n,-n+1,\dots$ we have $\lambda =
-w(w+n+3)$, $\mu=\ld (n-1)$ and \[
H(t) = \left( \begin{array}{c}
nt^{-n-1}\, \lw{3}F_2 \left( \begin{smallmatrix} w+2,\,\;-w-n-1,\,\; a + 1
\vspace{2pt} \\
-n,\,\;a
\end{smallmatrix} ;t\right) \vspace{2pt} \\ t^{-n}\, \lw{2}F_1 \left(
\begin{smallmatrix} w+3,\,\;-w-n \vspace{2pt} \\  1-n
\end{smallmatrix} ;t \right)
\end{array} \right),
\]
with $a=-w(w+n+3) - 2n-2$.

\noindent {\em b.1)} For $k=1$, and $w = -n,-n+1,\dots$ we have $\lambda
=-w(w+n+4) - n-2$, $\mu=(\lambda-3)(n+2)$ and
\[
H(t) = \left( \begin{array}{c}
t^{-n-1}\, \lw{2}F_1 \left( \begin{smallmatrix} w+3,\,\;-w-n-1 \vspace{2pt}
\\
-n
\end{smallmatrix} ;t \right) \vspace{2pt} \\ \frac bn t^{-n} \,\lw{3}F_2
\left( \begin{smallmatrix} w+3,\,\;-w-n-1,\,\;b+1
\vspace{2pt} \\ 1-n,\,\;b
\end{smallmatrix} ;t \right)
\end{array} \right),
\]
with $b=-w(w+n+4) - 2n - 3$.

\subsection{ The case $\ell=2$}
The complete list follows.

\smallskip
\noindent {\em a) Case $n\geq 0$.}

\noindent {\em a.0)} For $k=0$, and $w = 0,1,2,\dots$ we have $\lambda =
-w(w+n+4)$,
$\mu=\ld (n-2)$ and
$$ H_0=\left( \begin{array}{c}
1+\frac {\lambda(\lambda-3(n+1))}{2(n+1)(n+2)} \vspace{2pt} \\ 1 - \frac
{\lambda}{n+1}\vspace{2pt} \\
1 \end{array} \right),$$
\[
H(t) = \left( \begin{array}{c}
\left( 1+\frac {\lambda(\lambda-3(n+1))}{2(n+1)(n+2)} \right) \lw{4}F_3
\left( \begin{smallmatrix}
-w,\,\;w+n+4,\,\;s_1+1,\,\;s_2+1 \vspace{2pt} \\ n+3,\,\;s_1,\,\;s_2
\end{smallmatrix} ;t \right) \vspace{2pt} \\ \left(1 - \frac {\lambda}{n+1}
\right) \lw{3}F_2
\left( \begin{smallmatrix}
-w,\,\; w+n+4,\,\;s_3+1 \vspace{2pt} \\
n+2,\,\;s_3
\end{smallmatrix} ;t \right) \vspace{2pt} \\ \lw{2}F_1 \left(
\begin{smallmatrix}
-w,\,\;w+n+4 \vspace{2pt} \\
n+1
\end{smallmatrix} ;t \right)
\end{array} \right),
\]
with
\begin{align*}
s_1s_2 &= \textstyle\frac {w(w+3)(w+n+1)(w+n+4)}{2} + (n+1)(n+2),
\displaybreak[0]\\
s_1+s_2 &= -w(w+n+4)-n,\displaybreak[0] \\ s_3 &= -\textstyle\frac
{w(w+n+4)}{2} - \frac {n+1}{2}.
\end{align*}

\noindent {\em a.1)} For $k=1$, and $w = 0,1,2,\dots$ we have $\lambda =
-w(w+n+5) - n-2$,
$\mu=\lambda(n+1) - 6(n+2)$ and
$$ H_0=\left( \begin{array}{c}
\frac {\lambda}{(n+1)(n+2)} \vspace{2pt} \\ -\frac {\lambda +2}{2(n+1)}
\vspace{2pt} \\
1 \end{array} \right),$$
\[
H(t) = \left( \begin{array}{c}
\frac {\lambda}{(n+1)(n+2)}\, \lw{3}F_2 \left( \begin{smallmatrix}
-w,\,\;w+n+5,\,\;s_4+1 \vspace{2pt} \\
n+3,\,\;s_4
\end{smallmatrix} ;t \right) \vspace{2pt} \\ -\frac {\lambda +2}{2(n+1)}\,
\lw{4}F_3
\left( \begin{smallmatrix} -w-1,\,\;w+n+4,\,\;s_5+1,\,\;s_6+1 \vspace{2pt}
\\ n+2,\,\;s_5,\,\;s_6
\end{smallmatrix} ;t \right) \vspace{2pt} \\ \lw{3}F_2 \left(
\begin{smallmatrix}
-w-1,\,\;w+n+4,\,\;s_7+1 \vspace{2pt} \\ n+1,\,\;s_7
\end{smallmatrix} ;t \right)
\end{array} \right),
\]
with
\begin{align*}
s_4 &= -\textstyle\frac {w(w+n+5)}{2} - \frac {n+2}{2}, \displaybreak[0]\\
s_5+s_6 &=-\textstyle\frac {w(w+n+5)}{2} -  \frac
{1}{2}, \displaybreak[0] \\ s_5s_6 &= \textstyle\frac
{(w+1)(w+n+4)(w^2+nw+5w+n)}{8}, \displaybreak[0] \\ s_7
&=-\textstyle\frac {(w+1)(w+n+4)}{2}. \end{align*}

\noindent {\em a.2)} For $k=2$, and $w = 0,1,2,\dots$ we have $\lambda =
-w(w+n+6) - 2n - 6 $, $\mu=\lambda (n+4)-6(n+3)$  and
$$ H_0=\left( \begin{array}{c}
\frac {2}{(n+2)(n+1)} \vspace{2pt} \\
-\frac {2}{n+1} \vspace{2pt} \\
1 \end{array} \right),$$
\[
H(t) = \left( \begin{array}{c}
\frac {2}{(n+2)(n+1)}\, \lw{2}F_1 \left( \begin{smallmatrix} -w,\,\;w+n+6
\vspace{2pt} \\
n+3
\end{smallmatrix} ;t \right) \vspace{2pt} \\ -\frac {2}{n+1}\, \lw{3}F_2
\left( \begin{smallmatrix} -w-1,\,\;w+n+5,\,\;s_8+1
\vspace{2pt} \\
n+2,\,\;s_8
\end{smallmatrix} ;t \right) \vspace{2pt} \\ \lw{4}F_3 \left(
\begin{smallmatrix}
-w-2,\,\;w+n+4,\,\;s_9+1,\,\;s_{10}+1 \vspace{2pt} \\
n+1,\,\;s_9,\,\;s_{10} \end{smallmatrix} ;t \right)
\end{array} \right),
\]
with
\begin{align*}
s_8 &= -\textstyle\frac {w(w+n+6)}{2} - \frac {n+5}{2}, \displaybreak[0] \\
s_9s_{10} &= \textstyle\frac
{(w+1)(w+2)(w+n+4)(w+n+5)}{2}, \displaybreak[0]\\ s_9+s_{10} &= -w(w+n+6) -
(n+6).
\end{align*}

\smallskip
\noindent {\em b) Case $n\leq -2$.} We have, as above, three families going
with the choices $k = 0,1,2$.

\noindent {\em b.0)} For $k=0$, and $w = -n,-n+1,\dots$ we have $\lambda =
-w(w+n+4)$,
$\mu=\ld (n-2)$ and
\[
H(t) = \left( \begin{array}{c}
t^{-n-2}\, \lw{4}F_3\left( \begin{smallmatrix}
w+2,\,\;-w-n-2,\,\;s_1+1,\,\;s_2+1 \vspace{2pt} \\ -n-1,\,\;s_1,\,\;s_2
\end{smallmatrix}; t\right) \\
\frac 2{n+1} t^{-n-1}\, \lw{3}F_2 \left( \begin{smallmatrix}
w+3,\,\;-w-n-1,\,\;s_3+1 \\
-n,\,\;s_3 \end{smallmatrix}; t\right) \\ \frac 2{n(n+1)} t^{-n}\,
\lw{2}F_1 \left( \begin{smallmatrix} w+4,\,\;-w-n \\
 -n+1 \end{smallmatrix}; t\right)
\end{array} \right),
\]
with
\begin{align*}
s_1s_2 &= \textstyle\frac {(w+2)(w+3)(w+n+1)(w+n+2)}{2}, \displaybreak[0]\\
s_1+s_2 &=-w(w+n+4)-3n-4 ,\displaybreak[0]\\
s_3 &= -\textstyle\frac {w(w+n+4)}{2} - \frac {3}{2} (n+1). \end{align*}

\noindent {\em b.1)} For $k=1$, and $w = -n-1,n,-n+1,\dots$ we have
$\lambda = -w(w+n+5)-n-2 $, $\mu=\lambda(n+1) - 6(n+2)$
and \[ H(t)= \left( \begin{array}{c}
t^{-n-2}\, \lw{3}F_2 \left( \begin{smallmatrix} w+3,\,\;-w-n-2,\,\;s_1+1 \\
-n-1,\,\;s_1 \end{smallmatrix}; t\right) \\ \frac
{(\lambda-2n)}{2(n+1)} t^{-n-1}\, \lw{4}F_3 \left( \begin{smallmatrix}
w+3,\,\;-w-n-2,\,\;s_2+1,\,\;s_3+1 \\
-n,\,\;s_2,\,\;s_3
\end{smallmatrix}; t\right) \\
\left( \frac {\lambda}{n(n+1)} - \frac {2}{n} \right) t^{-n}\, \lw{3}F_2
\left( \begin{smallmatrix}
w+4,\,\;-w-n-1,\,\;s_4+1 \\
-n+1,\,\;s_4 \end{smallmatrix}; t \right) \end{array} \right),
\]
with $s_1,s_2,s_3,s_4$ given by the expressions \begin{align*}
s_2s_3 &= \textstyle\frac
{(w+3)(w+n+2)(w^2+nw+5w+3n+2)}{8},\displaybreak[0] \\ s_2+s_3 &=
-\textstyle\frac {w(w+n+5)}{2} - \frac {4n+5}{2},
\displaybreak[0]\\ s_1 &= -\textstyle\frac {w(w+n+5)}{2} - \frac
{(3n+6)}{2}, \displaybreak[0]\\
s_4 &= -\textstyle\frac {w(w+n+5)}{2} -\frac {(3n+4)}{2}. \end{align*}

\noindent {\em b.2)} For $k=2$, and $w = -n-2,-n-1,-n\dots$ we have
$\lambda = -w(w+n+6) - 2n - 6 $,
$\mu=\lambda(n+4) - 6(n+3)$ and \[
H(t) =\left( \begin{array}{c}
t^{-n-2}\, \lw{2}F_1 \left( \begin{smallmatrix} w+4,\,\;-w-n-2 \\ -n-1
\end{smallmatrix}; t\right) \\
\left( \frac {\lambda}{n+1} - 1 \right) t^{-n-1}\, \lw{3}F_2 \left(
\begin{smallmatrix} w+4,\,\;-w-n-2,\,\;s_1+1 \\
 -n,\,\;s_1 \end{smallmatrix}; t\right) \\ \left( \frac
{(\lambda-2)(\lambda+1)}{2n(n+1)} - \frac {\lambda+2}{2(n+1)} \right)
t^{-n}\, \lw{4}F_3 \left( \begin{smallmatrix}
w+4,\,\;-w-n-2,\,\;s_2+1,\,\;s_3+1 \\
-n+1,\,\;s_2,\,\;s_3
\end{smallmatrix}; t\right)
\end{array} \right),
\]
with $s_1,s_2,s_3$ given by
\begin{align*}
s_1 &= -\textstyle\frac {w(w+n+6)}{2} - \frac{3n+7}{2}, \displaybreak[0]\\
s_2s_3 &= \textstyle\frac {w(w+n+6)(w(w+n+6)  + 5n
+ 13)}{2} + (3n^2 + 15n + 20), \displaybreak[0]\\ s_2+s_3 &= -w(w+n+6) -
(3n+6).
\end{align*}

\smallskip
\noindent Finally, we come to the exceptional situation that arises when
$-l < n < 0$.

\

\noindent {\em c) Case $n= -1$.} As in all the previous situations there
are three families of spherical functions,
corresponding to the choices $k=0,1,2$. We observe that the expressions for
$\mu_k(\lambda)$ in Section \ref{aeigenfunc} are
still valid when $n < 0$.

\smallskip
\noindent {\em c.0)} For $k=0$, and $w = 1,2,3,\dots$ we have $\lambda =
-w(w+3)$, $\mu=-3\lambda $ and
\[
H(t)=\left( \begin{array}{c}
\lw{4}F_3\left( \begin{smallmatrix}
-w,\,\;w+3,\,\;s_1+1,\,\;s_2+1 \\
2,\,\;s_1,\,\;s_2
\end{smallmatrix} ;t \right) \\
\frac {2}{w(w+3)}\, \lw{3}F_2 \left( \begin{smallmatrix}
-w,\,\;w+3,\,\;s_3+1 \\ 1,\,\;s_3
\end{smallmatrix}; t\right) \\
-\frac 2{w(w+3)} t\, \lw{2}F_1 \left( \begin{smallmatrix} 1-w,\,\;w+4 \\2
\end{smallmatrix}; t\right)
\end{array} \right),
\]
and we have
\begin{align*}
s_1s_2 &= \textstyle\frac {w^2(w+3)^2}{2},\displaybreak[0] \\ s_1+s_2 &=
-w^2 - 3w+1 ,\displaybreak[0]\\ s_3 &= -
\textstyle\frac {w(w+3)}{2}.
\end{align*}

\noindent {\em c.1)} For $k=1$, and $w =0,1,2,\dots$ we have $\lambda =
-w(w+4)-1$, $\mu=-6 $ and
\[
H(t) =\left( \begin{array}{c}
\lw{3}F_2 \left( \begin{smallmatrix}
-w,\,\;w+4,\,\;s_1+1 \\
2,\,\;s_1 \end{smallmatrix} ;t \right) \\ -\frac
{(w^2+4w-1)}{2(w^2+4w+1)}\, \lw{4}F_3 \left( \begin{smallmatrix}
-w-1,\,\;w+3,\,\;s_2+1,\,\;s_3+1 \\ 1,\,\;s_2,\,\;s_3
\end{smallmatrix}; t \right) \\
t\, \lw{3}F_2 \left( \begin{smallmatrix} -w,\,\;w+4,\,\;s_1+1 \\
2,\,\;s_1
\end{smallmatrix}; t\right)
\end{array} \right),
\]
with
\begin{align*}
s_1 &= -\textstyle\frac {w^2+4w+1}{2},\displaybreak[0] \\ s_2s_3
&=\textstyle\frac {(w+1)(w+3)(w^2+4w-1)}{8},
\displaybreak[0] \\ s_2+s_3 &= -\textstyle\frac {w^2+4w+1}{2}. \end{align*}
\rmk{Remark} Note that in this case $H_2(t) = tH_0(t)$.

\smallskip
\noindent {\em c.2)} For $k=2$, and $w =0,1,2,\dots$ we have $\lambda =
-w(w+5)-4$, $\mu=3\lambda -12 $ and
\[
H(t) =\left( \begin{array}{c}
\lw{2}F_1 \left( \begin{smallmatrix}
-w,\,\;w+5 \\
2 \end{smallmatrix}; t \right) \\
- \lw{3}F_2 \left( \begin{smallmatrix}
-w-1,\,\;w+4,\,\;s_1+1 \\
1,\,\;s_1 \end{smallmatrix}; t\right) \\ -\frac 12(w+1)(w+4)t\, \lw{4}F_3
\left( \begin{smallmatrix}
-w-1,\,\;w+4,\,\;s_2+1,\,\;s_3+1 \\ 2,\,\;s_2,\,\;s_3 \end{smallmatrix};
t\right) \end{array} \right),
\]
with
\begin{align*}
s_1 &= -\textstyle\frac {(w+1)(w+4)}{2}, \\ s_2s_3 &= \textstyle\frac
{(w+1)^2(w+4)^2}{2}, \\ s_2+s_3 &= \textstyle  -(w^2 +
5w + 3).
\end{align*}

\section{A matrix valued form of the bispectral property}\label{recur}

In the case of $\ell =0$ when our spherical functions are scalar valued and
reduce to the well
known case of Jacobi polynomials, one has a three-term recursion relation
giving the product of $t$
times $\phi_w(t)$ as a linear combination, with coefficients independent of
$t$, of the functions
$\phi_{w-1}(t)$, $\phi_w(t)$, and $\phi_{w+1}(t)$.

The reasons behind this fact are well understood in terms of tensor
products of representations.
At a more elementary level this is just a consequence of dealing with a
sequence of orthogonal
polynomials.

{}From a different point of view this three-term recursion relation, {\it
along with} the differential equation in $t$ satisfied
by
$\phi_w(t)$, can be seen as a basic instance of the ``bispectral property"
discussed in a more general context, first in
\cite{DG} and subsequently in a variety of related contexts by several
authors. For a good reference showing
connections of this problem with lots of other areas see \cite{HK}. In this
field the basic cases
play a very crucial role, since they give rise to entire families of
nontrivial bispectral
situations by repeated applications of the Darboux process.

We finish this paper by displaying
a rather intriguing {\it matrix
valued} three term recursion relation that we have found for our spherical
functions. We have not
investigated the reasons behind it and we do not claim that this is the
most natural such recursion
in our context. We expect however to return to this problem in the future.
In fact, since the initial version of this paper was completed we have
made substantial progress in establishing the results illustrated below
and the interested reader may want to see \cite{GPT}.

For given nonnegative integers $n$, $\ell$ and $w$ consider the matrix
whose rows are given by the vectors $H(t)$
corresponding to the values $k=0,1,2,\dots ,\ell$ discussed above. Denote
the corresponding matrix
by
\[ \Phi (w,t). \]

\subsection{ The bispectral property}
As a function of $t$, $\Phi(w,t)$ satisfies two differential equations
\[
D \Phi(w,t)^T = \Phi(w,t)^T\Lambda , \quad E \Phi(w,t)^T=\Phi(w,t)^TM \ .
\]
Here $\Lambda$ and $M$ are diagonal matrices with $\Lambda (i,i)=
-w(w+n+i+\ell +1)- (i-1)(n+i)$,
$M(i,i)=\Lambda(i,i)(n-\ell+3i-3)-3(i-1)(\ell-i+2)(n+i)$, $1\leq
i\leq\ell+1$;
$D$ and $E$ are the differential operators introduced in \eqref{DD} and
\eqref{EE}.
Moreover we have

\begin{thm} There exist matrices $A_w$, $B_w$, $C_w$, independent of $t$,
such that
\[
A_w\Phi(w-1,t)+B_w\Phi(w,t)+C_w \Phi(w+1,t) = t\Phi(w,t)\ .
\]
\end{thm}

The matrices $A_w$ and $C_w$ consist of two diagonals each and $B_w$ is
tridiagonal.
Assume, for convenience, that these vectors are normalized in such a way
that for $t=1$ the matrix
$\Phi(w,1)$ consists of all ones.

We display below the formulas in the case $\ell=2$, $n\geq 0$. We leave the
$\ell$ dependence
quite explicit in these formulas to indicate how things go in the general
case. If we are not in
the case $\ell =2$ the multiplicative factor 2 in the off-diagonal elements
of the three matrices
below has to be altered.

The nonzero entries of the matrix $A_w$ are given as follows \begin{align*}
A_w(i,i)&=w(w+\ell+1)(w+n+i-1)(w+n+\ell+i)(w+\ell-i+2)^{-1}\\
&\times(w+n+2i-1)^{-1}(2w+n+\ell+i)^{-1}(2w+n+\ell+i+1)^{-1},
\displaybreak[0 ] \\
A_w(i,i+1)&=2w(w+\ell+1)(w+\ell-i+1)^{-1}(w+\ell-i+2)^{-1}\\
&\times(w+n+2i-1)^{-1}(2w+n+\ell+i+1)^{-1}.
\end{align*}
Notice that if $w=0$ we have $A_0=0$. This is useful since we
have not defined $\Phi(-1,t)$.

The nonzero entries of the matrix $C_w$ are given as follows
\begin{align*}
C_w&(i,i)\\
=&(w+1)(w+\ell+2)(w+n+i)(w+n+\ell+i+1)(w+\ell-i+2)^{-1}\\
&\times(w+n+2i-1)^{-1}(2w+n+\ell+i+1)^{-1}(2w+n+\ell+i+2)^{-1},
\displaybreak[0]\\ C_w&(i+1,i)\\
=&2(w+1)(w+\ell+2)(w+\ell-i+1)^{-1}\\
&\times(w+\ell-i+2)^{-1}(w+n+2i+1)^{-1}(2w+n+\ell+i+2)^{-1}. \end{align*}
Now we come to the matrix $B_w$. This one is tridiagonal and only its
off-diagonal elements have a
pleasant expression. We have
\begin{align*}
 B_w(i+1,i)=&2(w+n+i)(w+n+\ell+i+1)(w+\ell-i+1)^{-1}\\
&\times(w+n+2i)^{-1}(w+n+2i+1)^{-1}(2w+n+\ell+i+2)^{-1},\displaybreak[0]\\
B_w(i-1,i)=&2(w+n+i-1)(w+n+\ell+i)(w+\ell-i+3)^{-1}\\
&\times(w+n+2i-3)^{-1}(w+n+2i-2)^{-1}(2w+n+\ell+i)^{-1}.
\end{align*}
The expressions for the entries in the main diagonal of $B_w$ are given
below. For $B_w(1,1)$ we get
\begin{align*}
B_w(1,1)=&\big(z(2z+n^2+10n+13)+2(n+1)(n+3)(n+4)\big)\\
&\times(w+2)^{-1}(w+n+2)^{-1}(2w+n+3)^{-1}(2w+n+5)^{-1}, \end{align*}
with $z=w(w+n+4)$. For $B_w(2,2)$ we get \begin{align*}
B_w(2,2)=&\big(q(q(2q+(n+2)(n+12))+4(n^3+10n^2+28n+28))\\
&+(n+4)(3n^3+24n^2+56n+56)\big)(w+1)^{-1}(w+3)^{-1}\\
&\times(w+n+2)^{-1}(w+n+4)^{-1}(2w+n+4)^{-1}(2w+n+6)^{-1}, \end{align*}
with $q=w(w+n+5)$. For $B_w(3,3)$ we get \begin{align*}
B_w(3,3)=&\big(y(2y+n^2+10n+29)+2(n+5)(n^2+5n+10)\big)\\
&\times(w+2)^{-1}(w+n+4)^{-1}(2w+n+5)^{-1}(2w+n+7)^{-1}, \end{align*}
with $y=w(w+n+6)$.

Observe that the matrix $C_w$ is invertible if $w\geq 0$ and $n\geq 0$ as
we are assuming in this
section. This allows one to determine $\Phi(w,t)$ from $\Phi(0,t)$ and
$\Phi(1,t)$ for all
$w=2,3,\dots\ .$

\section{Appendix}\label{apendice}
For completeness we include here the proofs of those propositions and lemmas
given, without them, in Sections 4 and 5.
\begin{propsec4} For $H\in C^\infty(\CC^2)\otimes \End(V_\pi)$
we have
\begin{align*}
D_1(&H)(x,y)=\\
&(1+|x|^2+|y|^2) \left( \vphantom{(1+|y|^2)}
\,(H_{x_1x_1}+H_{x_2x_2})(1+|x|^2)+(H_{y_1y_1}+H_{y_2y_2})(1+|y|^2)\right.
\displaybreak[0]\\ & \left. \vphantom{(1+|y|^2)}
+2(H_{y_1x_1}+H_{y_2x_2})\,\re(x\overline y)
+2(H_{y_1x_2}-H_{y_2x_1})\,\im(x\overline y) \right).
\end{align*}
\end{propsec4}
\begin{proof}
We have
$$D_1(H)=-4 (X_{-\beta}X_{\beta}+X_{-\g}X_{\g})(H)
=\left(Y_3^2+Y_4^2+Y_5^2+Y_6^2\right)(H).$$
We begin by calculating $Y_5^2(H)(g)$, for all $g\in {\mathcal A}$. We have
$$Y_5^2(H)(g)= \left(\frac{d}{ds}\,\frac{d}{dt} H\left(p(g\exp
(s+t)Y_5)\right) \right)_{s=t=0}.$$
Let
$$u(s,t)=\frac{g_{12}\tan(s+t)+g_{13}}{g_{32}\tan(s+t)+g_{33}}\quad
\text{and}\quad
v(s,t)=\frac{g_{22}\tan(s+t)+g_{23}}{g_{32}\tan(s+t)+g_{33}}.$$
Then
\begin{equation*}
p(g\exp (s+t)Y_5)=(u(s,t)\, , \, v(s,t)\, , \, 1).
\end{equation*}
Now using Lemma \ref{menores} we obtain
$$\left(\frac{\partial u}{\partial s}\right)_{s=0} =
\frac{g_{12}g_{33}-g_{13}g_{32}}
{\left(g_{32}\sin t+g_{33}\cos t\right)^2 }=\frac{-\,\overline g_{21}}
{\left(g_{32}\sin t+g_{33}\cos t\right)^2 }$$
and
$$\left(\frac{\partial v}{\partial s}\right)_{s=0} =
\frac{g_{22}g_{33}-g_{23}g_{32}} {\left(g_{32}\sin t+g_{33}\cos
t\right)^2 }=\frac{\overline g_{11}} {\left(g_{32}\sin t+g_{33}\cos
t\right)^2 }.$$

\smallskip
\noindent Therefore at $s=t=0$ we have
$$ \frac{\partial u}{\partial s}
=\frac{\partial u}{\partial t}= -\frac{\overline g_{21}}{g_{33}^2 }\,
,\quad \frac{\partial v}{\partial s}
=\frac{\partial v}{\partial t}=\frac{\overline g_{11}}{g_{33}^2 }$$ and $$
\frac{\partial^2 u}{\partial
s\partial t}=\frac{2g_{32}\overline g_{21}}{g_{33}^3}\, ,\quad
\frac{\partial^2 v}{\partial s\partial t}=\frac{-2g_{32}\overline
g_{11}}{g_{33}^3}.$$

\noindent Similarly let
$$\tilde u(s,t)=\frac{i
\,g_{12}\tan(s+t)+g_{13}}{i\,g_{32}\tan(s+t)+g_{33}}\quad \text{and}\quad
\tilde v(s,t)=\frac{i\,g_{22}\tan(s+t)
+g_{23}}{i\,g_{32}\tan(s+t)+g_{33}}.$$
Then we obtain $ p(g\exp (s+t)Y_6)=(\tilde u(s,t)\, , \, \tilde v(s,t)\, ,
\, 1) $. So, at $s=t=0$ we have
$$ \frac{\partial \tilde u}{\partial s}
=\frac{\partial \tilde u}{\partial t}= -\frac{i\,\overline g_{21}}{g_{33}^2
}\, ,\quad
\frac{\partial\tilde v}{\partial s}
=\frac{\partial\tilde v}{\partial t}=\frac{i\,\overline g_{11}}{g_{33}^2 }$$
and
$$\frac{\partial^2 \tilde u}{\partial s\partial t}=-\frac{2g_{32}\overline
g_{21}}{g_{33}^3}\, ,\quad
\frac{\partial^2 \tilde v}{\partial s\partial t}=\frac{2g_{32}\overline
g_{11}}{g_{33}^3}.$$
By using the chain rule we have
\begin{align*}
Y_5^2(H)(g)&=\sum_{i,j=1,2} H_{x_ix_j}\, \frac{\partial u_i}{\partial
s}\,\frac{\partial u_j}{\partial t}
+\sum_{i,j=1,2} H_{y_iy_j}\, \frac{\partial v_i}{\partial
s}\,\frac{\partial v_j}{\partial t}\displaybreak[0]\\& +
\sum_{i,j=1,2} H_{x_iy_j}\,
\left(\frac{\partial u_i}{\partial s}\,\frac{\partial v_j}{\partial
t}+\frac{\partial u_i}{\partial t}\,
\frac{\partial v_j}{\partial s}\right) +\sum_{i=1,2}
H_{x_i}\frac{\partial^2 u_i}{\partial t\partial s}
+H_{y_i}\frac{\partial^2 v_i}{\partial t\partial s} \end{align*} and
\begin{align*} Y_6^2(H)(g)&=\sum_{i,j=1,2} H_{x_ix_j}\,
\frac{\partial \tilde u_i}{\partial s}\,\frac{\partial \tilde u_j}{\partial
t} +\sum_{i,j=1,2} H_{y_iy_j}\, \frac{\partial
\tilde v_i}{\partial s}\,\frac{\partial \tilde v_j}{\partial
t}\displaybreak[0]\\& +
\sum_{i,j=1,2} H_{x_iy_j}\,
\left(\frac{\partial \tilde u_i}{\partial s}\,\frac{\partial \tilde
v_j}{\partial t}+\frac{\partial\tilde u_i}{\partial t}\,
\frac{\partial \tilde v_j}{\partial s}\right)
+\sum_{i=1,2}H_{x_i}\frac{\partial^2 \tilde u_i}{\partial t\partial s}
+H_{y_i}\frac{\partial^2 \tilde v_i} {\partial t\partial s}.
\end{align*}
We note that if $u(t)=u_1(t)+iu_2(t)$ then
$\frac{du_1}{dt}=\re(\frac{du}{dt})$
and $\frac{du_2}{dt}=\im(\frac{du}{dt})$. Therefore
\begin{equation}\label{X-beXbe}
\begin{split}
(Y_5^2&+Y_6^2)(H)(g)=\left(H_{x_1x_1}+H_{x_2x_2}\right)\,
\frac{|g_{21}|^2}{|g_{33}|^4}
+\left(H_{y_1y_1}+H_{y_2y_2}\right)\, \frac{|g_{11}|^2}{|g_{33}|^4}\\&
\displaybreak[0]\\
&\,-\frac{2\re(g_{11}\overline g_{21})}{|g_{33}|^4}\,
\left(H_{x_1y_1}+H_{x_2y_2}\right)
-\frac{2\im(g_{11}\overline g_{21})}{|g_{33}|^4}\,
\left(H_{x_2y_1}-H_{x_1y_2}\right).
\end{split}
\end{equation}

\noindent We proceed in the same way with $Y_3^2+Y_4^2$ and obtain
\begin{equation}\label{X-gaXga}
\begin{split}
(Y_3^2&+Y_4^2)(H)(g)=\left(H_{x_1x_1}+H_{x_2x_2}\right)\,
\frac{|g_{22}|^2}{|g_{33}|^4}
+\left(H_{y_1y_1}+H_{y_2y_2}\right)\,
\frac{|g_{12}|^2}{|g_{33}|^4}\\&\displaybreak[0] \\
&\,-\frac{2\re(g_{12}\overline g_{22})}{|g_{33}|^4}\,
\left(H_{x_1y_1}+H_{x_2y_2}\right) -\frac{2\im(g_{12}\overline
g_{22})}{|g_{33}|^4}\, \left(H_{x_2y_1}-H_{x_1y_2}\right).
\end{split}
\end{equation}
Hence
\begin{align*}
(D_1)&(H)(g)= (Y_3^2+Y_4^2)(H)(g)+(Y_5^2+Y_6^2)(H)(g)\displaybreak[0] \\&
=\left(H_{x_1x_1}+H_{x_2x_2}\right)
\frac{(1-|g_{23}|^2)}{|g_{33}|^4} +\left(H_{y_1y_1}+H_{y_2y_2}\right)
\frac{\left(1-|g_{13}|^2\right)}{|g_{33}|^4}\displaybreak[0] \\ &+2
\left(H_{x_1y_1}+H_{x_2y_2}\right)\frac{\re(g_{13}\overline
g_{23})}{|g_{33}|^4} +2
\left(H_{x_2y_1}-H_{x_1y_2}\right)\frac{\im(g_{13}\overline
g_{23})}{|g_{33}|^4}.
\end{align*}
Now, using \eqref{expxy} the proposition follows. \end{proof}

\medskip
Recall that the open dense subset $\mathcal A$ of $G$ was defined by the
condition $\det(A(g))\neq0$ or equivalently by
$g_{33}\neq0$.

\begin{lem} \label{X-beta} For any $g\in {\mathcal A}$ let
$B=\matc{0}{g_{13}}{0}{g_{23}}$ and
$C=\matc{g_{13}}{0}{g_{23}}{0}$. Then
\begin{enumerate}
\item[i)] $(X_{-\beta} \Phi_\pi)(g)= \dot \pi ( B A(g)^{-1})\, \pi(A(g))$.

\item[ii)] $(X_{-\g} \Phi_\pi)(g)= \dot \pi (C A(g)^{-1})\, \pi(A(g))$.
\end{enumerate}
\end{lem}
\begin{proof}
Upon considering the unique holomorphic representation of $\GL(2,\CC)$
which extends $\pi$ we may write
$$\Phi_\pi(g\text{ exp } tX_{-\beta})= \pi(A(g\text{ exp }
tX_{-\beta})A(g)^{-1})\, \pi(A(g))$$
for $|t|$ small. Moreover we have
$$\left( X_{-\beta}
\Phi_\pi\right)(g)={\frac{d}{d\,t}\Big|}_{t=0}\Phi_\pi(g\text{ exp }
tX_{-\beta}).$$ Therefore
$$\left(X_{-\beta}\Phi_\pi\right)(g)=
\dot\pi\left(\left(\frac{d}{d\,t}A(g\text{ exp } tX_{-\beta})\right)_{t=0}
A(g)^{-1}\right)\,\pi(A(g)).$$
Since
$${\frac{d}{d\,t}\Big|}_{t=0}A(g\text{ exp } tX_{-\beta})=B,$$ i) follows.
In a similar way one establishes ii).
\end{proof}

\begin{prop4.5} For $H\in C^\infty(\CC^2)\otimes \End(V_\pi)$
we have
\begin{align*}
D_2(&H)(x,y)=\\
&- 4 \,\frac{\partial H}{\partial x}\, \dot \pi\begin{pmatrix}
x(1+|x|^2) & x^2{\overline y}\\ &\\ y(1+|x|^2)& x|y|^2
\end{pmatrix} - 4\,\frac{\partial H}{\partial y}\,
\dot \pi\begin{pmatrix} y|x|^2 & x(1+|y|^2) \\ &\\ y^2{\overline x}&
y(1+|y|^2) \end{pmatrix}.
\end{align*}

\end{prop4.5}
\begin{proof} We have
$$D_2(H)=-4\left(X_{\beta}(H)X_{-\beta}(\Phi_\pi)\Phi_\pi^{-1}
+X_{\g}(H)X_{-\g}(\Phi_\pi)\Phi_\pi^{-1} \right).$$
By Lemma \ref{X-beta}, for any $g\in {\mathcal A}$ we have
$$(X_{-\beta}\Phi_\pi)(g)\, \Phi_\pi(g)^{-1}= {\frac 1{\overline
g_{33}}}\,\dot\pi \matc{-g_{13}g_{21}\,}{\,g_{13}g_{11}}
{-g_{23}g_{21}\,}{\,g_{23}g_{11}} $$  and
$$ (X_{-\g}\Phi_\pi)(g)\,
\Phi_\pi(g)^{-1}= {\frac 1{\overline g_{33}}}\,\dot\pi
\matc{g_{13}g_{22}\,}{\;-g_{13}g_{12}}
{g_{23}g_{22}\,}{\;-g_{23}g_{12}}.
$$
Also, by Lemma \ref{XbetaH} we have
$$(X_\beta H)(g)=-{\frac{\overline g_{21}}{g_{33}^2}}\; \frac{\partial
H}{\partial x}+ { \frac{\overline g_{11}}{g_{33}^2}}\;
\frac{\partial H}{\partial y}\, , \qquad (X_\g H)(g)={\frac{\overline
g_{22}}{g_{33}^2}}\; \frac{\partial H}{\partial x}- {
\frac{\overline g_{12}}{g_{33}^2}}\; \frac{\partial H}{\partial y}.$$

\noindent Hence
\begin{align*}
D_2(H)(g)
=&-4\left( {\frac{\overline g_{21}}{g_{33}^2}}\; \frac{\partial H}{\partial
x}-{ \frac{\overline g_{11}}{g_{33}^2}}\;
\frac{\partial H}{\partial y}\right)\; {\frac 1{\overline g_{33}}}\,\dot\pi
\matc{g_{13}g_{21}\,}{\,-g_{13}g_{11}} {g_{23}g_{21}\,}{\,-g_{23}g_{11}}
\displaybreak[0]\\ & \\
& +4\left({\frac{\overline g_{22}}{g_{33}^2}}\; \frac{\partial H}{\partial
x}- { \frac{\overline g_{12}}{g_{33}^2}}\;
\frac{\partial H}{\partial y}\right)
\, {\frac 1{\overline g_{33}}}\,\dot\pi
\matc{-g_{13}g_{22}\,}{\;g_{13}g_{12}} {-g_{23}g_{22}\,}{\;g_{23}g_{12}}
\displaybreak[0]\\ &\\ =&-4\left(\frac{\partial H}{\partial
x}\, \dot \pi( P )- \frac{\partial H}{\partial y}\, \dot \pi(Q)\right),
\end{align*} where
\begin{align*}
P=&{\frac{\overline g_{21}}{g_{33}|g_{33}|^2}}
\matc{g_{13}g_{21}\,}{\,-g_{13}g_{11}} {g_{23}g_{21}\,}{\,-g_{23}g_{11}}
+{\frac{\overline g_{22}}{g_{33}|g_{33}|^2}}
\matc{g_{13}g_{22}\,}{\,-g_{13}g_{12}}{g_{23}g_{22}\,}{\,-g_{23}g_{12}}
\displaybreak[0]
\\ &\\ =& \frac{1}{g_{33}|g_{33}|^2}
\begin{pmatrix}
{(|g_{21}|^2+ |g_{22}|^2)g_{13}}& {-(g_{11}\overline g_{21}+g_{12}\overline
g_{22})g_{13}} \\ &\\
{(|g_{21}|^2+ |g_{22}|^2)g_{23}}& {-(g_{11}\overline g_{21}+g_{12}\overline
g_{22})g_{23}}
\end{pmatrix}\displaybreak[0]\\ &\\
=& \frac{1}{g_{33}|g_{33}|^2}
\begin{pmatrix}
{(1-|g_{23}|^2)g_{13}}& { \,g_{13}^2\overline g_{23}} \\ &\\
{(1-|g_{23}|^2)g_{23}}& { \,g_{13}|g_{23}|^2} \end{pmatrix}
=\begin{pmatrix} x(1+|x|^2)& x^2\overline y\\ &\\ y(1+|x|^2)& x|y|^2
\end{pmatrix}.
\end{align*}

\smallskip
In the last two steps we have used that $g\in \SU(3)$ and that
$x=g_{13}/g_{33}$, $y=g_{23}/g_{33}$.
If we proceed in the same way with $Q$ we complete the proof of the
proposition.
\end{proof}

\begin{prop4.6}
For $H\in C^\infty(\CC^2)\otimes \End(V_\pi)$ we have
\begin{align*}
E_1(H)&(x,y)=\\
&(1+|x|^2+|y|^2) \left[ \vphantom{\dot\pi
\matc{-(1+|x|^2)}{-3x\overline y}{0}{2(1+|x|^2)}} \, (H_{x_1x_1}
+H_{x_2x_2})\;\dot\pi \matc{-(1+|x|^2)}{-3x\overline
y}{0}{2(1+|x|^2)}\right.\\
&\\&\left. +(H_{y_1y_1}+H_{y_2y_2})\;\dot\pi
\matc{2(1+|y|^2)}{0}{-3\overline x y}{-(1+|y|^2)}\right.
\displaybreak[0]\\
&\\ & \left.
+(H_{y_1x_1}+H_{y_2x_2})\;\dot\pi \matc{2x\overline y-\overline
xy}{-3(1+|y|^2)}
{-3(1+|x|^2)}{2\overline x y-x \overline y}\right.\\ &\\ &\left.
+i\,(H_{x_2y_1}-H_{x_1y_2})\;\dot\pi \matc{\,-(2x\overline y
+\overline xy)\;}{\;-3\,(1+|y|^2)} {3\,(1+|x|^2)\;}{\;\,2\overline x y+x
\overline y} \right].
\end{align*}
\end{prop4.6}
\begin{proof}
We have
\begin{equation}\label{E1def}
\begin{split}
E_1(H)&=-4 \textstyle
\left(X_{-\beta}X_{\beta}(H)\right)\,\Ph_\pi\, \dot\pi(\tilde
H_1)\Ph_\pi^{-1} -4 \left(X_{-\gamma}X_{\gamma}(H)\right)\,
\Ph_\pi\,\dot\pi(\tilde H_2)\Ph_\pi^{-1} \\ & +12
\left(X_{-\beta}X_{\g}(H)\right)\,\Ph_\pi \dot\pi(X_{-\al})\Ph_\pi^{-1}
+12\left(X_{-\gamma}X_{\beta}(H)\right)\,\Ph_\pi\dot\pi(X_\al)\Ph_\pi^{-1}.
\end{split}
\end{equation}
In \eqref{X-beXbe} and \eqref{X-gaXga} we have calculated
$-4\left(X_{-\beta}X_{\beta}(H)\right)=(Y_5^2+Y_6^2)(H)$ and
$-4\left(X_{-\g}X_{\g}(H)\right)=(Y_3^2+Y_4^2)(H)$. Now we shall compute
\begin{equation}\label{X-bXgdef}
-4X_{-\beta}X_{\g}(H)= \left(Y_5Y_3+Y_6Y_4\right)(H)+i\, (
Y_6Y_3-Y_5Y_4)(H).
\end{equation}
\noindent For $U,V\in \lieg$ and $g\in \mathcal A$ we have
$$UV(H)(g)= \left(\frac{d}{ds}\,\frac{d}{dt} H \left(p(g\exp sU\,\exp
tV)\right) \right)_{s=t=0}.$$
If $|s|,|t|$ are small we can write
$$ p(g\exp sU \,\exp tV)= (x(s,t)\, , \, y(s,t)\, , \, 1).$$
Also if we put $x_1=\re(x)$, $x_2=\im(x)$,
$y_1=\re(y)$ and $y_2=\im(y)$, the chain rule gives
\begin{equation*} \begin{split}
UV(H)&(g)=\sum_{i,j=1}^2 H_{x_ix_j}\, \frac{\partial x_i}{\partial
s}\,\frac{\partial x_j}{\partial t}
+\sum_{i,j=1}^2 H_{y_iy_j}\, \frac{\partial x_i}{\partial
s}\,\frac{\partial x_j}{\partial t}\displaybreak[0]\\& +
\sum_{i,j=1}^2 H_{x_iy_j}\,
\left(\frac{\partial x_i}{\partial s}\,\frac{\partial y_j}{\partial
t}+\frac{\partial x_i}{\partial t}\,
\frac{\partial y_j}{\partial s}\right) +\sum_{i=1}^2
H_{x_i}\frac{\partial^2 x_i}{\partial t\partial s}
+H_{y_i}\frac{\partial^2 y_i}{\partial t\partial s}. \end{split}
\end{equation*} For $U=Y_5$ and $V=Y_3$ we have
$$x(s,t)= \frac{g_{11}\tan t +g_{12}\sin s+g_{13}\cos s} {g_{31}\tan t
+g_{32}\sin s+g_{33}\cos s}$$ and
$$y(s,t)=\frac{g_{21}\tan t +g_{22}\sin s+g_{23}\cos s} {g_{31}\tan t
+g_{32}\sin s+g_{33}\cos s}.$$
Moreover at $s=t=0$ we have
$$\frac{\partial x}{\partial s}=-\frac{\overline g_{21}}{g_{33}^2}\,,\qquad
\frac{\partial x}{\partial t}
=\frac{\overline g_{22}}{g_{33}^2}\, ,\qquad \frac{\partial y}{\partial
s}=\frac{\overline g_{11}}{g_{33}^2}\, ,\qquad
\frac{\partial y}{\partial t}=-\frac{\overline g_{12}}{g_{33}^2}.$$

\noindent For $U=Y_6$ and $V=Y_4$ we have
$$ x(s,t)= \frac{i\,g_{11}\tan t +i\,g_{12}\sin s+g_{13}\cos s}
{i\,g_{31}\tan t +i\,g_{32}\sin s+g_{33}\cos s}$$ and $$
y(s,t)=\frac{i\,g_{21}\tan t +i\,g_{22}\sin s+g_{23}\cos s} {i\,g_{31}\tan
t +i\,g_{32}\sin s+g_{33}\cos s}.$$
Moreover at
$s=t=0$ we have $$\frac{\partial x}{\partial s}=-\frac{i\, \overline
g_{21}}{g_{33}^2}\,,\qquad
\frac{\partial x}{\partial t}=\frac{i\,\overline g_{22}}{g_{33}^2}\,
,\qquad \frac{\partial y}{\partial s}
=\frac{i\,\overline g_{11}}{g_{33}^2}\, ,\qquad \frac{\partial y}{\partial
t}=-\frac{i\,\overline g_{12}}{g_{33}^2}.$$

\noindent For $U=Y_6$ and $V=Y_3$ we have
$$x(s,t)= \frac{g_{11}\tan t +i\,g_{12}\sin s+g_{13}\cos s} {g_{31}\tan t
+i\,g_{32}\sin s+g_{33}\cos s}$$ and
$$y(s,t)=\frac{g_{21}\tan t +i\,g_{22}\sin s+g_{23}\cos s} {g_{31}\tan t
+i\,g_{32}\sin s+g_{33}\cos s}.$$
Moreover at $s=t=0$
we have $$\frac{\partial x}{\partial s}=-\frac{i\,\overline
g_{21}}{g_{33}^2} \,,\qquad
\frac{\partial x}{\partial t}=\frac{\overline g_{22}}{g_{33}^2}\, ,\qquad
\frac{\partial y}{\partial s}
=\frac{i\,\overline g_{11}}{g_{33}^2}\, ,\qquad \frac{\partial y}{\partial
t}=-\frac{\overline g_{12}}{g_{33}^2}.$$

\noindent For $U=Y_5$ and $V=Y_4$ we have
$$x(s,t)= \frac{i\,g_{11}\tan t +g_{12}\sin s+g_{13}\cos s} {i\,g_{31}\tan
t +g_{32}\sin s+g_{33}\cos s}$$ and
$$y(s,t)=\frac{i\,g_{21}\tan t +g_{22}\sin s+g_{23}\cos s} {i\,g_{31}\tan t
+g_{32}\sin s+g_{33}\cos s}.$$
Moreover at $s=t=0$
we have $$\frac{\partial x}{\partial s}=-\frac{\overline
g_{21}}{g_{33}^2}\,,\qquad \frac{\partial x}{\partial t}
=\frac{i\,\overline g_{22}}{g_{33}^2}\, ,\qquad \frac{\partial y}{\partial
s}=\frac{\overline g_{11}}{g_{33}^2}\, ,\qquad
\frac{\partial y}{\partial t}=-\frac{i\,\overline g_{12}}{g_{33}^2}.$$

Now if we add the four terms in the right hand side of \eqref{X-bXgdef}
and observe that taking the real or imaginary part of a complex valued
function commutes with taking a derivative, we obtain
\begin{align*}
&-4X_{-\beta}X_{\g}(H)(g)=
-\frac{g_{21}\overline g_{22}}{|g_{33}|^4}
\left(H_{x_1x_1}+H_{x_2x_2}\right) -\frac{g_{11}\overline
g_{12}}{|g_{33}|^4}
\left(H_{y_1y_1}+H_{y_2y_2}\right)\displaybreak[0]\\
&+\frac{(g_{11}\overline g_{22}+g_{21}\overline g_{12})}{|g_{33}|^4}
\left(H_{x_1y_1}+H_{x_2y_2}\right)-\frac{i\,(g_{11}\overline
g_{22}-g_{21}\overline g_{12})}{|g_{33}|^4}
\left(H_{x_2y_1}-H_{x_1y_2}\right).
\end{align*}

\noindent In order to compute $X_{-\g}X_{\beta}$ we observe that
\begin{equation*}\begin{split}
-4X_{-\g}X_{\beta}(H)&= \left(Y_3Y_5+Y_4Y_6\right)(H)+i\, (
Y_4Y_5-Y_3Y_6)(H)\\
&=\left(Y_5Y_3+Y_6Y_4\right)(H)- i\, ( Y_6Y_3-Y_5Y_4)(H),
\end{split}
\end{equation*}
since if $X\in \liek$ then $X(H)=0$. Therefore it is easy to verify that
\begin{align*}
&-4X_{-\g}X_{\beta}(H)(g)=
-\frac{g_{22}\overline g_{21}}{|g_{33}|^4}
\left(H_{x_1x_1}+H_{x_2x_2}\right) -\frac{g_{12} \overline
g_{11}}{|g_{33}|^4}
\left(H_{y_1y_1}+H_{y_2y_2}\right)\displaybreak[0] \\ &
+\frac{(g_{22}\overline g_{11}+g_{12}\overline g_{21})}{|g_{33}|^4}
\left(H_{x_1y_1}+H_{x_2y_2}\right)
-\frac{i\,(g_{22}\overline g_{11}-g_{12}\overline g_{21})}{|g_{33}|^4}
\left(H_{x_2y_1}-H_{x_1y_2}\right).
\end{align*}

\smallskip
\noindent On the other hand we have
\begin{equation}\label{Phpipunto}
\begin{split}
\Ph_\pi(g)\dot\pi(\tilde H_1)\Ph_\pi(g)^{-1}&= \dot\pi\left(A(g) \tilde H_1
A(g)^{-1}\right)\displaybreak[0] \\
&=\frac{1}{\overline g_{33}}
\dot\pi \matc{2g_{11}g_{22}+g_{12}g_{21}}{-3g_{11}g_{12}}{3g_{21}g_{22}}
{-g_{11}g_{22}-2g_{12}g_{21}},\displaybreak[0]\\
&\\ \Ph_\pi(g)\dot\pi(\tilde H_2)\Ph_\pi(g)^{-1}&=\frac{1}{\overline
g_{33}} \dot\pi \matc{-g_{11}g_{22}-2g_{12}g_{21}}
{3g_{11}g_{12}}{-3g_{21}g_{22}}
{2g_{11}g_{22}+g_{12}g_{21}},\displaybreak[0]\\ &\\
\Ph_\pi(g)\dot\pi(X_{-\al})\Ph_\pi(g)^{-1}&=\frac{1}{\overline g_{33}}
\dot\pi \matc{g_{12}g_{22}}{-g_{12}^2}{g_{22}^2}
{-g_{12}g_{22}},\displaybreak[0]\\ &\\
\Ph_\pi(g)\dot\pi(X_\al)\Ph_\pi(g)^{-1}&=\frac{1}{\overline g_{33}} \dot\pi
\matc{-g_{11}g_{21}}{g_{11}^2}{-g_{21}^2}
{g_{11}g_{21}}.
\end{split}
\end{equation}

Now, by \eqref{E1def}, the function multiplying $H_{x_1x_1}+H_{x_2x_2}$ is
a matrix $\dot\pi(A)$ where $A$ is given by:
\begin{align*}
A=&\frac{|g_{21}|^2}{|g_{33}|^4\overline g_{33}}\,
\matc{2g_{11}g_{22}+g_{12}g_{21}}{-3g_{11}g_{12}}{3g_{21}g_{22}}
{-g_{11}g_{22}-2g_{12}g_{21}}\\ & +
\frac{|g_{22}|^2}{|g_{33}|^4\overline g_{33}}\,
\matc{-g_{11}g_{22}-2g_{12}g_{21}}{3g_{11}g_{12}}{-3g_{21}g_{22}}
{2g_{11}g_{22}+g_{12}g_{21}}\\ &+\frac{3g_{21}\overline
g_{22}}{|g_{33}|^4\overline
g_{33}}\matc{g_{12}g_{22}}{-g_{12}^2}{g_{22}^2}{-g_{12}g_{22}} +
\frac{3\overline g_{21} g_{22}}{|g_{33}|^4\overline
g_{33}}\matc{-g_{11}g_{21}}{g_{11}^2}{-g_{21}^2} {g_{11}g_{21}}.
\end{align*}

\medskip
If we recall that $g\in \SU(3)$, $x=g_{13}/g_{33}$ and $y=g_{23}/g_{33}$
then the coefficient $A_{11}$ of $A$ can be written
in the following way, \begin{align*} A_{11}=&\frac{1}{|g_{33}|^4\overline
g_{33}}\, \left(|g_{21}|^2
(-g_{11}g_{22}+g_{12}g_{21}) + |g_{22}|^2 (-g_{11}g_{22} +g_{12}g_{21})
\right)\displaybreak[0]\\ =& \frac{-\left(
|g_{21}|^2+|g_{22}|^2\right)}{|g_{33}|^4}=\frac{-\left(
1-|g_{23}|^2\right)}{|g_{33}|^4}\displaybreak[0] \\ =&
-(1+|x|^2)(1+|x|^2+|y|^2).
\end{align*}
Similarly,
\begin{align*}
A_{22}=&\frac{1}{|g_{33}|^4\overline g_{33}} \, \left( |g_{21}|^2+
|g_{22}|^2\right) (-2g_{21}g_{12}+2g_{11}g_{22})
=\frac{2}{|g_{33}|^4}\,{\left( 1-|g_{23}|^2\right)}\displaybreak[0]\\  =&
2(1+|x|^2)(1+|x|^2+|y|^2).
\end{align*}
We also have $A_{21}=0$ and
\begin{align*}
A_{12}=&\frac{1}{|g_{33}|^4\overline g_{33}} \, \left( 3g_{11}\overline
g_{21 }(-g_{21}g_{12}+g_{11}g_{22})
+ 3g_{12}\overline g_{22 }(g_{11}g_{22}-g_{21}g_{12})
\right)\displaybreak[0]\\ =& \frac{3}{|g_{33}|^4}\,(g_{11} \overline
g_{21}+g_{12}\overline g_{22})= \frac{-3g_{13}\overline
g_{23}}{|g_{33}|^4}\displaybreak[0]\\ =& -3x\overline y
(1+|x|^2+|y|^2).
\end{align*}

\noindent Analogously we calculate the functions multiplying
$H_{y_1y_1}+H_{y_2y_2}$, $H_{y_1x_1}+H_{y_2x_2}$,
$H_{x_2y_1}-H_{x_1y_2}$
and complete the proof of the proposition.
\end{proof}

\begin{prop4.7}
For $H\in C^\infty(\CC^2)\otimes \End(V_\pi)$ we have
\begin{align*} &E_2(H)(x,y)=\\
& 4\frac{\partial H}{\partial x}\!\!\left(\dot\pi\!\matc{0}{x}{0}{y}
\dot\pi\!
\matc{-2x\overline y }{0}{3(1+|x|^2)}{\;x\overline y}
+ \dot\pi\!\matc{x}{0}{y}{0} \dot\pi\!\matc{1+|x|^2}{3x\overline
y}{0}{-2(1+|x|^2)}\!\right)
\displaybreak[0]\\ &\\
& \!\!\!+4\frac{\partial H}{\partial y}\!\!\left(\dot\pi\!\matc{0}{x}{0}{y}
\dot\pi\!\matc{-2(1+|y|^2)}{0} {3y\overline x}
{1+|y|^2}+ \dot\pi\!\matc{x}{0}{y}{0}\dot\pi\!\matc{y\overline x
}{\,3(1+|y|^2)}{0}{-2y\overline x}\!\right). \end{align*}
\end{prop4.7}
\begin{proof} We have
\begin{align*}
E_2(H)&= -4X_{\beta}(H)X_{-\beta}(\Ph_\pi)\, \dot\pi (\tilde
H_1)\Ph_\pi^{-1}-4X_{\gamma}(H)X_{-\gamma}(\Ph_\pi) \,
\dot\pi(\tilde H_2)\Ph_\pi^{-1}\displaybreak[0]\\ & +12
X_{\g}(H)X_{-\beta}(\Ph_\pi)\dot\pi(X_{-\al})\Ph_\pi^{-1} +12
X_{\beta}(H)X_{-\gamma}(\Ph_\pi) \dot\pi(X_\al)\Ph_\pi^{-1}.
\end{align*}
For any $g\in {\mathcal A}$ we have, by Lemma \ref{XbetaH}
$$(X_\beta H)(g)=-{\frac{\overline g_{21}}{g_{33}^2}}\;
\frac{\partial H}{\partial x}+ { \frac{\overline g_{11}}{g_{33}^2}}\;
\frac{\partial H}{\partial y}\, , \qquad
(X_\g H)(g)={\frac{\overline g_{22}}{g_{33}^2}}\;  \frac{\partial
H}{\partial x}- { \frac{\overline g_{12}}{g_{33}^2}}\;
\frac{\partial H}{\partial y},$$
and by Lemma \ref{X-beta}
$$(X_{-\beta} \Phi_\pi)(g)\Phi_\pi(g)^{-1}= \frac{1}{\overline g_{33}}\,
\dot \pi
\matc{-g_{13}g_{21}\;}{g_{13}g_{11}}{-g_{23}g_{21}\;}{g_{23}g_{11}}, $$
$$(X_{-\g} \Phi_\pi)(g)\Phi_\pi(g)^{-1}=
\frac{1}{\overline g_{33}}\, \dot \pi
\matc{g_{13}g_{22}\;}{-g_{13}g_{12}}{g_{23}g_{22}\;}{-g_{23}g_{12}}.$$

\noindent In \eqref{Phpipunto} we computed
$\Ph_\pi(g)\dot\pi(A)\Ph_\pi(g)^{-1}$, for $A=\tilde H_1, \tilde H_2,
X_{-\al}$ and $X_{\al}$. Therefore,
$E_2(H)=4\big(\frac{\partial H}{\partial x} \, P+\frac{\partial H}{\partial
y} Q\big)$
with
\begin{align*}
P=& -\frac{\overline g_{21}}{|g_{33}|^4}\,\dot\pi\matc{g_{13}g_{21}\;}
{-g_{13}g_{11}}{g_{23}g_{21}\;}{-g_{23}g_{11}} \,
\dot\pi \matc{2g_{11}g_{22}+g_{12}g_{21}}{-3g_{11}g_{12}}{3g_{21}g_{22}}
{-g_{11}g_{22}-2g_{12}g_{21}}
\displaybreak[0]\\
&+\frac{\overline g_{22}}{|g_{33}|^4}\, \dot \pi
\matc{-g_{13}g_{22}\;}{g_{13}g_{12}}{-g_{23}g_{22}\;}{g_{23}g_{12}}\,\dot\pi
\matc{-g_{11}g_{22}-2g_{12}g_{21}}{3g_{11}g_{12}}{-3g_{21}g_{22}}
{2g_{11}g_{22}+g_{12}g_{21}}
\displaybreak[0]\\
&-\frac{3\overline g_{22}}{|g_{33}|^4}\,\dot\pi\matc{g_{13}g_{21}\;}
{-g_{13}g_{11}}{g_{23}g_{21}\;}{-g_{23}g_{11}} \,
\dot\pi \matc{g_{12}g_{22}}{-g_{12}^2}{g_{22}^2}{-g_{12}g_{22}}
\displaybreak[0]\\ & +\frac{3\overline g_{21}}{|g_{33}|^4}\,
\dot \pi
\matc{-g_{13}g_{22}\;}{g_{13}g_{12}}{-g_{23}g_{22}\;}{g_{23}g_{12}}\, \dot
\pi\matc{-g_{11}g_{21}}{g_{11}^2}{-g_{21}^2} {g_{11}g_{21}}.
\end{align*}
If we combine the first and third terms and the second and fourth terms,
and use Lemma \ref{menores} several times we obtain
\begin{align*}
P=\,& \frac{1}{|g_{33}|^4}\,\dot\pi\matc{g_{13}g_{21}\;}
{-g_{13}g_{11}}{g_{23}g_{21}\;}{-g_{23}g_{11}}
\displaybreak[0]\\
&\times \dot\pi
\matc{2g_{22}g_{13}\overline g_{23} -g_{12}(1-|g_{23}|^2)}
{-3g_{12}g_{13}\overline
g_{23}}{-3g_{22}(1-|g_{23}|^2)}{2g_{12}(1-|g_{23}|^2)-g_{22}g_{13}\overline
g_{23}}\\
&+\frac{1}{|g_{33}|^4}\,\dot\pi\matc{-g_{13}g_{22}\;}{g_{13}g_{12}}{-g_{23}g
_{
 22}\;}{g_{23}g_{12}} \displaybreak[0] \\
&\times\dot\pi \matc{2g_{21}g_{13}\overline g_{23}
-g_{11}(1-|g_{23}|^2)} {-3g_{11}g_{13}\overline
g_{23}}{-3g_{21}(1-|g_{23}|^2)}{2g_{11}(1-|g_{23}|^2)-g_{21}g_{13}\overline
g_{23}}.
\end{align*}

\smallskip
\noindent Hence using \eqref{expxy} we get
\begin{align*}
P&=\frac{1}{|g_{33}|^2}\,\dot\pi\matc{g_{13}g_{21}\;}
{-g_{13}g_{11}}{g_{23}g_{21}\;}{-g_{23}g_{11}} \displaybreak[0]\\
& \times\left( \dot\pi
\matc{2g_{22}x\overline y} {0}{-3g_{22}(1+|x|^2)\;}{\;-g_{22}x\overline y}
+\dot\pi \matc{-g_{12}(1+|x|^2)}
{-3g_{12}x\overline y}{0}{2g_{12}(1+|x|^2)}\right) \\
&+\frac{1}{|g_{33}|^2}\,\dot\pi\matc{-g_{13}g_{22}\;}{g_{13}g_{12}}{-g_{23}g
_{
 22}\;}{g_{23}g_{12}} \displaybreak[0] \\
&\times\left( \dot\pi \matc{2g_{21}x\overline y}
{0}{-3g_{21}(1+|x|^2)\;}{\;-g_{21}x\overline y} +\dot\pi
\matc{-g_{11}(1+|x|^2)} {-3g_{11}x\overline y}{0}{2g_{11}(1+|x|^2)}\right).
\end{align*}

\noindent Thus,
\begin{align*}
P=&\,\frac{g_{22}}{|g_{33}|^2}\,\dot\pi\matc{g_{13}g_{21}\;}
{-g_{13}g_{11}}{g_{23}g_{21}\;}{-g_{23}g_{11}} \dot\pi
\matc{2x\overline y} {0}{-3(1+|x|^2)\;}{\;-x\overline y}\displaybreak[0]\\
&\;+\frac{g_{21}}{|g_{33}|^2}\,
\dot\pi\matc{-g_{13}g_{22}\;}{g_{13}g_{12}}{- g_{23}g_{
22}\;}{g_{23}g_{12}} \dot\pi\matc{2x\overline y} {0}
{-3(1+|x|^2)\;}{\;-x\overline y}\displaybreak[0]\\ &\;
+\frac{g_{12}}{|g_{33}|^2}\,\dot\pi\matc{g_{13}g_{21}\;}
{-g_{13}g_{11}}{g_{23}g_{21}\;} {-g_{23}g_{11}}\dot\pi \matc{-(1+|x|^2)}
{-3x\overline y}{0}{2(1+|x|^2)}\displaybreak[0]
\\&\; + \frac{g_{11}}{|g_{33}|^2}\,\dot\pi\matc{-g_{13}g_{22}\;}
{g_{13}g_{12}}{-g_{23}g_{22}\;}{g_{23}g_{12}}
\dot\pi\matc{-(1+|x|^2)} {-3x\overline y}{0}{2(1+|x|^2)}.
\end{align*}

\smallskip
\noindent Now, using Lemma \ref{menores} we obtain
\begin{align*}
P=&\,\frac{1}{|g_{33}|^2}\,\dot\pi\matc{0\;}{\;-g_{13}\overline
g_{33}}{0\;}{\;-g_{23}\overline g_{33}}
\dot\pi \matc{2x\overline y} {0}{-3(1+|x|^2)\;}{\;-x\overline
y}\displaybreak[0] \\ &\;+
\frac{1}{|g_{33}|^2}\,\dot\pi\matc{-g_{13}\overline
g_{33}\;}{\;0}{-g_{23}\overline g_{33}\;}{\;0} \dot\pi
\matc{-(1+|x|^2)}{-3x\overline y}{0}{2(1+|x|^2)}\displaybreak[0]
\\ =&\, \dot\pi\matc{0}{\;x}{0}{\,y}
\dot\pi \matc{-2x\overline y }{\,0}{3(1+|x|^2)}{\;\;x\overline y}+
\dot\pi\matc{x}{\,0}{y}{\,0}
\dot\pi \matc{1+|x|^2}{3x\overline y}{0}{-2(1+|x|^2)}.
\end{align*}

\smallskip
\noindent If we proceed in the same way with $Q$ we complete the proof of
the proposition.
\end{proof}

\smallskip

Now we shall compute all first and second order partial derivatives that we
used
in Section \ref{secc1var}.  We start with the following immediate results.

\begin{lem5.3} At $(r,0)\in \CC^2$ We have
$$H_{x_1}(r,0)=\frac{d\tilde H}{dr}(r) \quad \text{ and }\quad
H_{x_1x_1}(r,0)=\frac{d^2\tilde H}{dr^2}(r).$$
\end{lem5.3}

As a preparation to other computations, given $(x,y)\in \CC^2$,
$x=x_1+ix_2$ and $y=y_1+i y_2$ we need to choose an element in $K$ which
carries the point $(x,y)$ to the meridian $\{(r,0)\, : \, r>0\}$.
We take
$$A(x,y)=\frac 1{s(x,y)}
\begin{pmatrix}x& -\bar y \\ y& \bar x \end{pmatrix} \in
\SU(2) \; \text{ where }\;
s(x,y)=\sqrt{|x|^2+|y|^2}. $$
Then $(x,y)=A(x,y)\,(s(x,y),0)$ and
$$H(x,y)=\pi(A(x,y))\tilde H(s(x,y))\pi(A(x,y)^{-1}).$$

For simplicity we shall use the notation
$(u_1,u_2,u_3,u_4)=(x_1,x_2,y_1,y_2)$.
As usual we denote $[A,B]=AB-BA$.

\begin{prop}\label{der2H}
At $(r,0)\in \CC^2$ we have
$$\frac{\partial H}{\partial u_j}(r,0)= \left[ \dot\pi\Bigl(\frac{\partial
A}{\partial u_j}\Bigr), \tilde H(r)\right] +\delta_{1j} \frac{d\tilde
H}{dr}.$$
and
\begin{equation*}
\begin{split}
\frac{\partial^2 H}{\partial u_i\partial u_j}&(r,0)=
\frac{\partial^2 (\pi\circ A)}{\partial u_i\partial u_j} \tilde H(r) +
\tilde H(r)  \frac{\partial^2 (\pi\circ A^{-1})}{\partial u_i\partial u_j}
+\delta_{1i}\delta_{1j} \frac{d^2\tilde H}{dr^2}  \displaybreak[0] \\
&+\frac 1r
\delta_{ij}(1-\delta_{1j}) \frac{d\tilde H}{dr}
+\delta_{1j}\Bigl[
 \dot\pi\Bigl(\frac{\partial A}{\partial u_i}\Bigr),
\frac{d\tilde H}{dr} \Bigr]
 +\delta_{1i}\Bigl[ \dot\pi\Bigl( \frac{\partial A}{\partial u_j}\Bigr),
\frac{d\tilde H}{dr} \Bigr] \displaybreak[0] \\
& -\dot\pi\Bigl( \frac{\partial A}{\partial u_j}\Bigr) \tilde H(r)
\dot\pi\Bigl(\frac{\partial A}{\partial u_i}\Bigr)
-\dot\pi\Bigl( \frac{\partial A}{\partial u_i}\Bigr) \tilde H(r)
\dot\pi\Bigl(\frac{\partial A}{\partial u_j}\Bigr).
\end{split}
\end{equation*}
\end{prop}

\begin{proof}
 We have  $$H(x,y)=\pi(A(x,y))\tilde H(s(x,y))\pi(A(x,y)^{-1}),$$
then
\begin{equation}\label{der1}
\begin{split}
\displaystyle\frac{\partial H}{\partial u_j} =&
 \frac{\partial (\pi\circ A)}{\partial u_j}  \,( \tilde H \circ s )\,
(\pi\circ A^{-1})
+ (\pi \circ A)\, \frac{\partial (\tilde H\circ s)}{\partial u_j}\,
(\pi\circ A^{-1}) \displaybreak[0] \\
&+ (\pi\circ A )\, (\tilde H\circ s )\, \frac{\partial (\pi\circ
A^{-1})}{\partial u_j}.
\end{split}
\end{equation}
Moreover
$$ \frac{\partial (\tilde H\circ s)}{\partial
u_j}=\Bigl(\frac{d\tilde H}{dr}\circ s\Bigr) \frac{\partial s}{\partial
u_j}.$$

\noindent We also note that $A(r,0)=I$ and that at
$(r,0)$ we have
$$\frac{\partial s}{\partial u_j}=\delta_{1j}\, ,\qquad
 \frac{\partial (\pi\circ A)}{\partial u_j} =
\dot\pi \Bigl(\frac{\partial  A}{\partial u_j} \Bigr)\, , \qquad
\frac{\partial (\pi\circ A^{-1})}{\partial u_j} =
-\dot\pi \Bigl(\frac{\partial  A}{\partial u_j} \Bigr).$$
Now replacing in \eqref{der1} the first statement follows.

To compute the second order derivative we start observing that at $(r,0)$
we have
\begin{align*}
 \frac{\partial ^2 (\tilde H \circ s)}{\partial u_i\partial u_j}
&=\frac{d^2 \tilde H}{dr^2}\frac{\partial s}{\partial u_i}\frac{\partial
s}{\partial u_j}+\frac{d\tilde H}{dr}\frac{\partial ^2 s}{\partial
u_i\partial u_j}
 =\delta_{1i}\delta_{1j} \frac{d^2\tilde H}{dr^2} +\frac 1r
\delta_{ij}(1-\delta_{1j})\frac{d\tilde H}{dr}.
\end{align*}

\noindent Now we differentiate the expression given in \eqref{der1} with
respect to $u_i$ and evaluate it at $(r,0)\in \CC^2$.
\begin{equation*}
\begin{split}
\displaystyle\frac{\partial^2  H}{\partial u_i\partial u_j}(r,0) =&
 \frac{\partial^2 (\pi\circ A)}{\partial u_i\partial u_j}  \, \tilde H(r)
+\delta_{1i}\dot\pi\Bigl(\frac{\partial A}{\partial u_j}
\Bigr)\frac{d\tilde H}{dr}
- \dot\pi\Bigl(\frac{\partial A}{\partial u_j}\Bigr) \tilde H(r)
\dot\pi\Bigl(\frac{\partial A}{\partial u_i}\Bigr)\displaybreak[0]\\
&+\delta_{1j}\dot\pi\Bigl(\frac{\partial A}{\partial u_i}
\Bigr)\frac{d\tilde H}{dr}
+ \delta_{1i}\delta_{1j} \frac{d^2\tilde H}{dr^2}
+\frac 1r \delta_{ij}(1-\delta_{1j})\frac{d\tilde H}{dr}\displaybreak[0]\\
&- \delta_{1j}\frac{d\tilde H}{dr}
\dot\pi\Bigl(\frac{\partial A}{\partial u_i}\Bigr)
-\dot\pi\Bigl(\frac{\partial A}{\partial u_i}\Bigr) \tilde H(r)
\dot\pi\Bigl(\frac{\partial A}{\partial u_j}\Bigr)
\displaybreak[0]\\
& - \delta_{1i}\frac{d\tilde H}{dr}
\dot\pi\Bigl(\frac{\partial A}{\partial u_j}\Bigr)
+ \tilde H(r)
\frac{\partial^2 (\pi\circ A^{-1})}{\partial u_i\partial u_j}.
\end{split}
\end{equation*}
Now rearranging the right hand side of the above expression the proof is
finished.
\end{proof}

\begin{prop}\label{derprincipal}
At $(r,0)\in \CC^2$ we have
\begin{equation}\label{der2piA}
\begin{split}
\frac{\partial^2 (\pi\circ A)}{\partial u_i\partial u_j} =& \,
\dot\pi\Bigl(\frac{\partial^2  A}{\partial u_i\partial u_j} \Bigr)
-\tfrac 12 \dot\pi\Bigl(\frac{\partial A}{\partial u_i} \frac{\partial
A}{\partial u_j}\Bigr)
-\tfrac 12 \dot\pi\Bigl(\frac{\partial A}{\partial u_j} \frac{\partial
A}{\partial u_i}\Bigr)\displaybreak[0]\\
& +\tfrac 12 \dot\pi\Bigl(\frac{\partial A}{\partial
u_j}\Bigr)\dot\pi\Bigl( \frac{\partial A}{\partial u_i}\Bigr)
+\tfrac 12 \dot\pi\Bigl(\frac{\partial A}{\partial
u_i}\Bigr)\dot\pi\Bigl(
\frac{\partial A}{\partial u_j}\Bigr)
\end{split}
\end{equation}
and
\begin{equation}\label{der2piA-1}
\begin{split}
\frac{\partial^2 (\pi\circ A^{-1})}{\partial u_i\partial u_j} =&
-\dot\pi\Bigl(\frac{\partial^2  A}{\partial u_i\partial u_j} \Bigr)
+\tfrac 12 \dot\pi\Bigl(\frac{\partial A}{\partial u_i} \frac{\partial
A}{\partial u_j}\Bigr)
+\tfrac 12 \dot\pi\Bigl(\frac{\partial A}{\partial u_j} \frac{\partial
A}{\partial u_i}\Bigr)\displaybreak[0]\\
& +\tfrac 12 \dot\pi\Bigl(\frac{\partial A}{\partial
u_j}\Bigr)\dot\pi\Bigl( \frac{\partial A}{\partial u_i}\Bigr)
+\tfrac 12 \dot\pi\Bigl(\frac{\partial A}{\partial
u_i}\Bigr)\dot\pi\Bigl(
\frac{\partial A}{\partial u_j}\Bigr).
\end{split}
\end{equation}
\end{prop}

\begin{proof}
For $|x|$ and $|y|$ sufficiently small we consider
\begin{equation}\label{logA}
X(x,y)=\log(A(x,y))= B(x,y)-\frac{B(x,y)^2}2+\frac{B(x,y)^3}3-\cdots,
\end{equation}
where $B(x,y)=A(x,y)-I$. Then
$$\pi(A(x,y))=\pi(\text{exp } X(x,y))=\text{exp }\dot\pi
(X(x,y))=\displaystyle\sum_{j\geq 0}\frac{\dot\pi (X(x,y))^j}{j!}.$$
Now we differentiate with respect to $u_j$ to obtain
\begin{equation}\label{derpi}
\begin{split}
\frac{\partial \, (\pi\circ A)}{\partial u_j}= &\, \dot\pi\Bigl(
\frac{\partial X}{\partial u_j}\Bigr)
+\tfrac 1{2!}\dot\pi\Bigl( \frac{\partial X}{\partial u_j}\Bigl)
\dot\pi(X)
+\tfrac 1{2!}\dot\pi(X) \dot\pi\Bigl( \frac{\partial
X}{\partial u_j}\Bigr) \displaybreak[0] \\
&+ \tfrac 1{3!}\dot\pi\Bigl( \frac{\partial X}{\partial u_j}\Bigr)
\dot\pi(X)^2
+ \tfrac 1{3}\dot\pi (X) \dot\pi\Bigl( \frac{\partial
X}{\partial u_j}\Bigr) \dot\pi\Bigl(X\Bigr) \displaybreak[0] \\
&+ \tfrac 1{3!}\dot\pi(X)^2
\dot\pi\Bigl( \frac{\partial X}{\partial u_j}\Bigr) +\cdots .
\end{split}
\end{equation}

\noindent Since $X(r,0)=0$, if we differentiate \eqref{derpi} and evaluate
at $(r,0)$ we obtain
\begin{equation}\label{der2piAaux}
\frac{\partial^2 \,(\pi\circ A)}{\partial u_i\, \partial u_j}  = \dot
\pi\Bigl(\frac{\partial^2 X}{\partial
u_i \partial u_j}\Bigr)
+\tfrac 12 \dot\pi\Bigl(\frac{\partial X}{\partial u_j}\Bigr)\,
\dot\pi\Bigl(\frac{\partial X}{\partial u_i}\Bigr)
+\tfrac 12 \dot\pi\Bigl(\frac{\partial X}{\partial u_i}\Bigr)
\dot\pi\Bigl(\frac{\partial X}{\partial u_j}\Bigr).
\end{equation}

\noindent To compute $\dfrac{\partial X}{\partial u_j}$ and
$\dfrac{\partial^2 X}{\partial u_i \partial u_j}$ we differentiate
\eqref{logA} and we get

\begin{align*}
\frac{\partial X}{\partial u_j}=& \frac{\partial B}{\partial u_j}
-\tfrac 1{2} \Bigl( \frac{\partial B}{\partial u_j}\Bigr)B- \tfrac 12 B
\Bigl(\frac{\partial B}{\partial u_j}\Bigr)+ \tfrac 13 \Bigl(\frac{\partial
B}{\partial u_j}\Bigr)B^2  \\ &
+\tfrac 13  B\Bigl(\frac{\partial B}{\partial
u_j}\Bigr) B +\tfrac 13 B^2 \Bigl(\frac{\partial B}{\partial
u_j}\Bigr)+\cdots .
\end{align*}

\noindent Since $B(r,0)=0$ we have
$$\frac{\partial X}{\partial u_j}(r,0)= \frac{\partial B}{\partial
u_j}(r,0)=\frac{\partial A}{\partial u_j}(r,0).$$

\noindent We also have at $(r ,0)\in \CC^2$
$$\frac{\partial^2 X}{\partial u_i \partial u_j} =
\frac{\partial^2 A}{\partial u_i \partial u_j}
- \tfrac 12 \frac{\partial A}{\partial u_j }\,\frac{\partial A}{\partial
u_i}
- \tfrac 12 \frac{\partial A}{\partial u_i }\,\frac{\partial A}{\partial u_j
}.$$
Replacing in \eqref{der2piAaux} the first statement in the lemma follows.

\smallskip
Now we compute $\dfrac{\partial^2 (\pi\circ A^{-1})}{\partial
u_i\partial u_j}$. We observe that

$$\frac{\partial (\pi\circ A^{-1})}{\partial u_j}= -(\pi\circ A^{-1})\,
\frac{\partial (\pi\circ A)}{\partial u_j}\, (\pi\circ A^{-1}).$$

\noindent Therefore at $(r,0)\in \CC^2$ we have,
\begin{align*}
\frac{\partial^2 (\pi\circ A^{-1})}{\partial u_i \partial u_j} =
&- \frac{\partial (\pi\circ A^{-1})}{\partial u_i}\, \frac{\partial
(\pi\circ A)}{\partial u_j}
-\frac{\partial^2 (\pi\circ A)}{\partial u_i \partial u_j}
-\frac{\partial (\pi\circ A)}{\partial u_j}\frac{\partial (\pi\circ
A^{-1})}{\partial u_i}.
\end{align*}

\noindent Then the statement in \eqref{der2piA-1} follows from
\eqref{der2piA} and
$$ \frac{\partial (\pi\circ A)}{\partial u_j} =
\dot\pi \Bigl(\frac{\partial  A}{\partial u_j} \Bigr)\, , \qquad
\frac{\partial (\pi\circ A^{-1})}{\partial u_j} =
-\dot\pi \Bigl(\frac{\partial  A}{\partial u_j} \Bigr).$$
\end{proof}

\

Now we specialize the above propositions in the  different cases we need.

\begin{lem5.4}
At $(r,0)\in \CC^2$ we have
\begin{equation*}
H_{y_1}(r,0)= -{\frac 1r}\Bigl( \dot\pi(J)\tilde H(r)-\tilde
H(r)\dot\pi(J)\Bigr)
\end{equation*}
and
\begin{equation*}
\begin{split}
H_{y_1y_1}(r,0)&= {\frac 1r} \frac{d \tilde H}{d r} +
\frac{1}{r^2}\left(\dot\pi(J)^2\, \tilde H(r)+ \tilde H(r)\,
\dot\pi(J)^2\right)-\frac 2{r^2}\dot\pi(J)\tilde H(r)\, \dot\pi(J),
\end{split}
\end{equation*}
where $J=\matc{0}{1}{-1}{0}$.
\end{lem5.4}
\begin{proof} To compute ${\partial A}/{\partial u_3}$ at $(r,y_1)\in
\CC^2$ we consider
 $$A(x,y)=\frac 1{\sqrt{r^2+y_1^2}}
\begin{pmatrix}r& - y_1 \\ y_1& r \end{pmatrix}.$$
Therefore
$$\frac{\partial  A}{\partial u_3}(r,0)=-{\frac 1{r}} J\quad{\textstyle and
}
\quad\frac{\partial^2  A}{\partial u_3^2}(r,0)=-{\frac 1{r^2}} I.$$

\noindent By Proposition \ref{derprincipal} we have at $(r,0)\in \CC^2$
$$\frac{\partial^2(\pi\circ A)}{\partial u_3^2}=\frac{\partial^2(\pi\circ
A^{-1})}{\partial u_3^2}= \frac 1{r^2} \dot\pi(J)^2 .$$
Now from Proposition  \ref{der2H} we obtain
\begin{equation*}
\frac{\partial H}{\partial y_1}(r,0)= -{\frac 1r}\Bigl( \dot\pi(J)\tilde
H(r)-\tilde
H(r)\dot\pi(J)\Bigr),
\end{equation*}
\begin{equation*}
\begin{split}
\frac{\partial^2 H}{\partial y_1^2}&(r,0)=
\frac 1{r^2} \dot\pi(J)^2  \tilde H(r) +  \frac 1{r^2} \tilde H(r)
\dot\pi(J)^2+\frac 1r \frac{d\tilde H}{dr}-\frac 2{r^2} \dot\pi(J) \tilde
H(r)\dot\pi(J).
\end{split}
\end{equation*}
This completes the proof of the lemma.
\end{proof}

\begin{lem5.5}
At $(r,0)\in \CC^2$ we have
\begin{equation*}\label{dery2}
H_{y_2}(r,0)= \frac{i}{r}\left( \dot\pi(T)\tilde H(r)-\tilde
H(r)\dot\pi(T)\right)
\end{equation*}
and
\begin{equation*}
\begin{split}
H_{y_2y_2}(r,0)&= {\frac 1r} \frac{d \tilde H}{d r}
-\frac{1}{r^2}\left(\dot\pi(T)^2\, \tilde H(r)+ \tilde H(r)\,
\dot\pi(T)^2\right)+\frac 2{r^2}\dot\pi(T)\tilde
H(r)\, \dot\pi(T),
\end{split}
\end{equation*}
where $T=\matc{0}{1}{1}{0}$.
\end{lem5.5}

\begin{proof} In this case we take $A(x,y)=\frac 1{\sqrt{r^2+y_2^2}}
\begin{pmatrix}r& iy_2 \\ iy_2& r \end{pmatrix}.$
Therefore
$$\frac{\partial  A}{\partial u_4}(r,0)= {\frac i{r}} T\quad{\textstyle and
}
\quad\frac{\partial^2  A}{\partial u_4^2}(r,0)=-{\frac 1{r^2}} I.$$

\noindent By Proposition \ref{derprincipal} we have at $(r,0)\in \CC^2$
$$\frac{\partial^2(\pi\circ A)}{\partial u_4^2}=\frac{\partial^2(\pi\circ
A^{-1})}{\partial u_4^2}= -\frac 1{r^2} \dot\pi(T)^2 .$$
Now from Proposition  \ref{der2H} we get
$$\frac{\partial H}{\partial y_2}(r,0)= \frac{i}{r}\left( \dot\pi(T)\tilde
H(r)-\tilde
H(r)\dot\pi(T)\right),$$
\begin{equation*}
\begin{split}
\frac{\partial^2 H}{\partial y_2^2}&(r,0)=
-\frac 1{r^2} \dot\pi(T)^2  \tilde H(r) -  \frac 1{r^2} \tilde H(r)
\dot\pi(T)^2+\frac 1r \frac{d\tilde H}{dr}+\frac 2{r^2} \dot\pi(T) \tilde
H(r)\dot\pi(T).
\end{split}
\end{equation*}
The proof of the lemma is finished.
\end{proof}

\begin{lem5.6}
At $(r,0)\in \CC^2$ we have
\begin{equation*}
H_{x_2}(r,0)= {\frac ir}\left( \dot\pi(H_\alpha)\tilde H(r)-\tilde
H(r)\dot\pi(H_\alpha)\right)
\end{equation*}
and
\begin{equation*}\begin{split}
H_{x_2x_2}(r,0)&= {\frac 1r} \frac{d \tilde H}{d r} -
\frac{1}{r^2}\left(\dot\pi(H_\alpha)^2\, \tilde H(r)+ \tilde H(r)\,
\dot\pi(H_\alpha)^2\right)
\displaybreak[0]\\ & \quad +\frac
2{r^2}\dot\pi(H_\alpha)\,\tilde H(r)\, \dot\pi(H_\alpha).
\end{split}
\end{equation*}
\end{lem5.6}

\begin{proof}
In this case we have $x=r+ix_2$, $y=0$, and $$A(x,y)=\frac 1{\sqrt{r^2+y^2}}
\begin{pmatrix}r+ix& 0 \\0& r-ix \end{pmatrix}.$$
Therefore
$$\frac{\partial  A}{\partial u_2}(r,0)= \frac 1r\matc{i}{0}{0}{-i}={\frac
i{r}}
H_\alpha\quad{\textstyle and }
\quad\frac{\partial^2  A}{\partial u_2^2}(r,0)=-{\frac1{r^2}} I.$$

\noindent By Proposition \ref{derprincipal} we have at $(r,0)\in \CC^2$
$$\frac{\partial^2(\pi\circ A)}{\partial u_2^2}=\frac{\partial^2(\pi\circ
A^{-1})}{\partial u_2^2}= -\frac 1{r^2} \dot\pi(H_\alpha)^2 .$$
Now from Proposition  \ref{der2H} we have
\begin{equation*}
\frac{\partial H}{\partial x_2}(r,0)= {\frac ir}\left(
\dot\pi(H_\alpha)\tilde H(r)-\tilde
H(r)\dot\pi(H_\alpha)\right)
\end{equation*}
\begin{equation*}
\begin{split}
\frac{\partial^2 H}{\partial x_2^2}(r,0)=&
\frac 1{r^2} \dot\pi(H_\alpha)^2  \tilde H(r) +  \frac 1{r^2} \tilde H(r)
\dot\pi(H_\alpha)^2+\frac 1r \frac{d\tilde H}{dr}\\ &-\frac 2{r^2}
\dot\pi(H_\alpha) \tilde H(r)\dot\pi(H_\alpha).
\end{split}
\end{equation*}
This completes the proof of the lemma.
\end{proof}

\begin{lem5.7}
At $(r,0)\in \CC^2$ we have
\begin{equation*}
H_{x_1y_1}(r,0)=\! \frac {-1}r \Big( \dot\pi(J)\frac{d \tilde H}{d
r}\!- \!\frac{d \tilde H}{d r}\dot\pi(J)\Big)\!
 + \!\frac{1}{r^2}\left( \dot\pi(J)\tilde H(r)-\tilde H(r)\dot\pi(J)\right),
\end{equation*}
\begin{equation*}
H_{x_1y_2}(r,0)= {\frac ir}\Bigl( \dot\pi(T)\frac{d \tilde H}{d
r}-\frac{d \tilde H}{d r}\dot\pi(T)\Bigr) -
\frac{i}{r^2}\left( \dot\pi(T)\tilde H(r)-\tilde H(r)\dot\pi(T)\right).
\end{equation*}
\end{lem5.7}
\begin{proof} It follows from Lemmas \ref{ddery1} and \ref{ddery2} by
differentiating with res\-pect to $r$.
\end{proof}

\begin{lem5.8}
At $(r,0)\in \CC^2$ we have
\begin{equation*}
\begin{split}
H_{y_1x_2}(r,0)=& -
\frac{i}{2r^2}\Big(\dot\pi(H_\alpha)\,\dot\pi(J)
+\dot\pi(J)\,\dot\pi(H_\alpha)\Big)\tilde H(r)\displaybreak[0] \\
& - \frac{i}{2r^2}\tilde H(r)\,\Big(\dot\pi(H_\alpha)\,\dot\pi(J)
+\dot\pi(J)\,\dot\pi(H_\alpha)\Big)\displaybreak[0]\\ &
+\frac{i}{r^2}\left(\dot\pi(H_\alpha)\,\tilde H(r)\,\dot\pi(J)
+\dot\pi(J)\,\tilde H(r)\,\dot\pi(H_\alpha)\right).
\end{split}
\end{equation*}
\end{lem5.8}
\begin{proof}
In this case we have $x=r+ix_2$, $y=y_1$, and
$$A(x,y)=\frac 1{\sqrt{r^2+x^2+y^2}}
\begin{pmatrix}r+ix& -y_1 \\y_1& r-ix \end{pmatrix}.$$
Therefore
$$\frac{\partial  A}{\partial u_2}(r,0)={\frac i{r}} H_\alpha,\qquad
\frac{\partial  A}{\partial u_3}(r,0)=-{\frac 1{r}} J, \qquad
\frac{\partial^2  A}{\partial u_3\partial u_2}(r,0)=0.$$
By Proposition \ref{derprincipal} we have at $(r,0)\in \CC^2$
$$\frac{\partial^2(\pi\circ A)}{\partial u_3\partial
u_2}=\frac{\partial^2(\pi\circ A^{-1})}{\partial u_3\partial u_2}
= -\frac{i}{2r^2}\Big(\dot\pi(H_\alpha)\dot\pi(J)
+\dot\pi(J)\dot\pi(H_\alpha)\Big).$$
Now the lemma follows from Proposition  \ref{der2H}.
\end{proof}

\begin{lem5.9}
At $(r,0)\in \CC^2$ we have
\begin{equation*}
\begin{split}
H_{y_2x_2}(r,0)=& -
\frac{1}{2r^2}\Big(\dot\pi(H_\alpha)\,\dot\pi(T)
+\dot\pi(T)\,\dot\pi(H_\alpha)\Big)\tilde H(r)\displaybreak[0] \\ & -
\frac{1}{2r^2}\tilde H(r)\,\Big(\dot\pi(H_\alpha)\,\dot\pi(T)
+\dot\pi(T)\,\dot\pi(H_\alpha)\Big)\displaybreak[0]\\ &
+\frac{1}{r^2}\left(\dot\pi(H_\alpha)\,\tilde H(r)\,\dot\pi(T)
+\dot\pi(T)\,\tilde H(r)\,\dot\pi(H_\alpha)\right).
\end{split}
\end{equation*}
\end{lem5.9}
\begin{proof}
In this case we have $x=r+ix_2$, $y=y_2$, and
$$A(x,y)=\frac 1{\sqrt{r^2+x^2+y^2}}
\begin{pmatrix}r+ix& iy_2 \\y_2& r-ix \end{pmatrix}.$$
Therefore
$$\frac{\partial  A}{\partial u_2}(r,0)={\frac i{r}} H_\alpha,\qquad
\frac{\partial  A}{\partial u_4}(r,0)={\frac i{r}} T, \qquad
\frac{\partial^2  A}{\partial u_4\partial u_2}(r,0)=0.$$
By Proposition \ref{derprincipal} we have at $(r,0)\in \CC^2$
$$\frac{\partial^2(\pi\circ A)}{\partial u_4\partial
u_2}=\frac{\partial^2(\pi\circ A^{-1})}{\partial u_4\partial u_2}
= -\frac{1}{2r^2}\left(\dot\pi(H_\alpha)\dot\pi(T)
+\dot\pi(T)\dot\pi(H_\alpha)\right).$$
Now the lemma follows from  Proposition  \ref{der2H}.
\end{proof}

\newcommand\bibit[5]{\bibitem
{#2}#3, {\em #4,\!\! } #5}

\end{document}